\documentclass[11pt,letter]{amsart}
\usepackage{amscd}
\usepackage{amssymb}
\usepackage[centering,text={15.5cm,22cm}]{geometry}
\usepackage{graphicx}
\usepackage{color}
\usepackage[all]{xy}
\usepackage{mathrsfs}
\usepackage{marvosym}
\usepackage{stmaryrd}
\usepackage{srcltx}
\usepackage{hyperref}

\definecolor{shadecolor}{rgb}{1,0.9,0.7}

\newtheorem{theorem}{Theorem}[section]
\newtheorem{lemma}[theorem]{Lemma}
\newtheorem{proposition}[theorem]{Proposition}
\newtheorem{corollary}[theorem]{Corollary}

\theoremstyle{definition}
\newtheorem{definition}[theorem]{Definition}
\newtheorem{construction}[theorem]{Construction}
\newtheorem{discussion}[theorem]{Discussion}
\newtheorem{example}[theorem]{Example}
\newtheorem{examples}[theorem]{Examples}

\theoremstyle{remark}
\newtheorem{remark}[theorem]{Remark}

\numberwithin{equation}{section}
\numberwithin{figure}{section}

\newcommand {\lfor} {\llbracket}
\newcommand {\rfor} {\rrbracket}

\newcommand{\NN} {\mathbb{N}}
\newcommand{\ZZ} {\mathbb{Z}}
\newcommand{\QQ} {\mathbb{Q}}
\newcommand{\RR} {\mathbb{R}}
\newcommand{\CC} {\mathbb{C}}

\newcommand{\PP} {\mathbb{P}}
\renewcommand{\AA} {\mathbb{A}}
\newcommand{\GG} {\mathbb{G}}

\newcommand {\shF}  {\mathcal{F}}

\newcommand {\shJ}  {\mathcal{J}}

\newcommand {\shL}  {\mathcal{L}}
\newcommand {\shM}  {\mathcal{M}}

\newcommand {\shW}  {\mathcal{W}}

\newcommand {\bSpec}  {\operatorname{\mathbf{Spec}}}
\newcommand {\GS} {{\underline{\mathrm{GS}}}}

\newcommand {\cl}  {\operatorname{cl}}

\newcommand {\Ext}  {\operatorname{Ext}}

\newcommand {\Gm} {\GG_m}
\newcommand {\gp}  {{\operatorname{gp}}}

\newcommand {\Hom}  {\operatorname{Hom}}

\newcommand {\id}  {\operatorname{id}}
\newcommand {\im}  {\operatorname{im}}

\newcommand {\Int}  {\operatorname{Int}}

\newcommand {\Isom}  {\operatorname{Isom}}

\renewcommand {\ker } {\operatorname{ker}}
\newcommand {\kk} {\Bbbk}

\newcommand {\Log} {\mathcal{L}og\hspace{1pt}}

\newcommand {\lra}  {\longrightarrow}
\newcommand {\ls}  {\dagger}

\newcommand {\M} {\mathcal{M}}

\newcommand {\Mbf} {\mathbf{M}}

\newcommand {\maxid} {\mathfrak{m}}

\newcommand {\MM} {\mathscr{M}}

\newcommand {\N}  {\operatorname{N}}

\renewcommand{\O}  {\mathcal{O}}

\newcommand {\Ob}  {\operatorname{Ob}}
\newcommand {\ol} {\overline}
\newcommand {\ord}  {\operatorname{ord}}

\newcommand {\Pic}  {\operatorname{Pic}}

\newcommand {\pr}  {\operatorname{pr}}

\newcommand {\scrM}  {\mathscr{M}}

\newcommand {\sat}  {{\operatorname{sat}}}

\newcommand {\Sch} {(\mathrm{Sch})}
\newcommand {\SchS} {(\mathrm{Sch}/\underline{S})}

\newcommand {\Spec} {\operatorname{Spec}}

\newcommand {\Trop}  {\operatorname{Trop}}

\newcommand {\ul} {\underline}

\newcommand {\U} {\mathscr{U}}
\newcommand {\V} {\mathscr{V}}

\def\mapright#1{\smash{
  \mathop{\longrightarrow}\limits^{#1}}}

\def\mydate{\ifcase\month \or January\or February\or March\or
April\or May\or June\or July\or August\or September\or October\or 
November\or December\fi \space\number\day,\space\number\year}

\begin{document}

\title[Log GW-invariants]{Logarithmic Gromov-Witten Invariants}

\author{Mark Gross} \address{UCSD Mathematics,
9500 Gilman Drive, La Jolla, CA 92093-0112, USA}
\email{mgross@math.ucsd.edu}
\thanks{This work was partially supported by NSF grants 0505325 and
0805328.}

\author{Bernd Siebert} \address{FB Mathematik,
Universit\"at Hamburg, Bundesstra\ss e~55, 20146 Hamburg,
Germany}
\email{bernd.siebert@math.uni-hamburg.de}

\date{\today}
\maketitle

\tableofcontents

\section*{Introduction.}

The purpose of this paper is the development of a general theory of
Gromov-Witten invariants in logarithmically smooth situations.
Concrete examples of such situations are non-singular projective
varieties with a normal crossing divisor, central fibres of
semi-stable degenerations, or toroidal pairs. They occur naturally
in Gromov-Witten theory in imposing constraints, in dimensional
reduction and in degeneration situations. The easiest, and already
quite useful special case is Gromov-Witten invariants relative to a
smooth divisor. The first treatments used symplectic techniques
\cite{LiRuan},\cite{IonelParker}. Algebraically a direct approach
for very ample divisors is possible \cite{Gathmann}, the much more
complicated general case is due to Jun Li \cite{JunLi1}, \cite{JunLi2}.
Any of the general approaches use a geometrically beautiful, but
technically unpleasant change of target space with the purpose of
avoiding irreducible components mapping into the divisor. As an
application these authors also developed decomposition formulas of
absolute Gromov-Witten invariants into sums of relative
Gromov-Witten invariants under a semistable degeneration of the
target variety into two components intersecting along a smooth
divisor (symplectically, a symplectic cut).

While Jun Li's theory was under development the second named author
suggested that a far more general and potentially technically simpler
theory should be possible using abstract log geometry \cite{Talk}.
This paper is the late realization of this proposal. While several
problems could be solved back then, notably stable reduction and the
construction of virtual fundamental classes, there remained a
fundamental problem of selecting the natural (``basic'') log
structure on the base scheme of a family of stable log maps. Without
a notion of basicness it is virtually impossible to prove
algebraicity and quasi-compactness of the relevant moduli stack.
With hindsight one can say that the suggestion of log Gromov-Witten
theory was premature in 2001 because two major ingredients were only
just developing, tropical geometry on the one hand and more powerful
techniques for the treatment of log moduli problems on the other
hand. In fact, the essential insight for the notion of basicness
came in a discussion of the authors in August~2008 and was guided by
our understanding of moduli of tropical curves and their relation to
log geometry. As for log moduli problems Olsson's algebraic stack of
fine log structures \cite{OlssonCC} is crucial to prove algebraicity
of our stack of basic stable log maps, and Olsson's many other works
on problems and applications of log geometry served as a blueprint
at many stages of this work.

The main concept is that of a \emph{basic stable log map}
(Definition~\ref{Def: Stable log map} and Definition~\ref{Def:
Basicness}). A stable log map is just a stable map with all arrows
defined in the category of log schemes and the condition on the
domain to be pre-stable replaced by log smoothness. We show that
stable log maps to a fixed log scheme $X=(\ul X,\M_X)$ form an
algebraic stack $\tilde\MM(X)$ (Theorem~\ref{Thm: algebraicity of
stack}). The main ingredients in this proof are Olsson's algebraic
stack $\Log$ of fine log schemes, the understanding of pre-stable
curves from the log point of view \cite{Kato 2000}, and a
representability result for spaces of log morphisms in relative
situations (Proposition~\ref{LMor is representable}). The stack
$\tilde \MM(X)$ is far too large because it allows arbitrary log
structures on the base of a family of stable log maps. The notion of
basicness (Definition~\ref{Def: Basicness}) selects a universal
choice that interacts well with geometry. Basicness only depends on
the log morphism on the level of ghost sheaves. The relation to
tropical geometry comes by pulling back to the standard log point
(\S\ref{par: Log maps over std log pt}). In situations where
tropical geometry applies, the moduli space of associated tropical
curves is the dual of the basic monoid (the stalk of the ghost sheaf
of the log structure on the base). In general, this picture is a bit
problematic, but there is still a way to visualize stable log maps
by (families of) tropical curves with values in a topological,
piecewise polyhedral space $\Trop(X)$ that is canonically associated
to $X$ (Appendix~\ref{App: Tropicalization}). In any case, basicness
selects an open substack $\MM(X)$ of $\tilde \MM(X)$ which now also
has a separated diagonal. It also comes with a natural log
structure. This is the first main result (stated as
Proposition~\ref{representability over M x M(X)} and
Corollary~\ref{MM(X) is algebraic log stack} in the text). We work
over a fixed fine saturated log scheme $S= (\ul S,\M_S)$.

\begin{theorem}
Let $X=(\ul X,\M_X)$ be a fine saturated log scheme. Then the stack
$\MM(X/S)$ of basic stable log maps to $X$ over $S$ is an algebraic
log stack locally of finite type over $\ul S$. The forgetful
morphism $\MM(X/S)\to {\mathbf M}(\ul X/\ul S)$ to the stack of
ordinary stable maps over $\ul S$ is representable.
\end{theorem}

A basic problem, which we could not solve in complete generality, is
to identify quasi-compact substacks of $\MM(X/S)$. Natural
conditions concern the classical conditions on genus, number of
marked points and homology class of the underlying ordinary stable
map, plus logarithmic conditions at the marked points similar to the
orders of tangency with the divisor in the classical situation.
Denoting $\MM(X/S,\beta)$ the substack of $\MM(X/S)$ fulfilling such
a set of conditions, a stratawise approach reduces the question of
boundednesss to a generally subtle question in convex geometry
(Theorem~\ref{Thm: Boundedness}), summarized in the concept of
\emph{combinatorial finiteness} of a class $\beta$
(Definition~\ref{Def: Combinatorial finiteness}). We have complete
results provided the pull-back of $\ol\M_X^\gp\otimes_{\ul \ZZ}{\ul
\QQ}$ to any stable log map is globally generated, that is, if any
germ of a section is the restriction of a global section. Then any
$\beta$ is combinatorially finite (Theorem~\ref{quasi-generated
boundedness}). This criterion holds in the important special cases
of genus $0$, of simple normal crossing divisors and in toric
situations. We also have boundedness without further conditions as
long as the dual intersection graph of the domain has genus at most
one (Theorem~\ref{Thm: genus one}) and in certain favourable
stratawise situations (Theorem~\ref{lastboundednesstheorem}). In
concrete situations one can check combinatorial finiteness by
examining a finite list of ordinary stable maps.

Arguably the most difficult single result in this paper is stable
reduction (Theorem~\ref{Thm: Stable reduction}), which requires a
solid understanding of the interaction between the basicness
condition and geometry. Taken together with boundedness we obtain
the second main result (Corollary~\ref{properness of the stack of
stable log maps} in the text).

\begin{theorem}\label{Main Thm II}
For $\ul X$ projective over $\ul S$ and $\beta$ a combinatorially
finite class of stable log maps, $\MM(X/S,\beta)$ is proper over
$\ul S$.
\end{theorem}

The construction of a \emph{virtual fundamental class} on
$\MM(X/S,\beta)$ now is completely straightforward via the approach
of \cite{behrendfantechi} using Olsson's log cotangent complex
\cite{OlssonCC}, as already worked out by Kim \cite{kim}. The details
are given in Section~\ref{Sect: VFC}. One can then define log
Gromov-Witten invariants in the usual way by pairing with cohomology
classes on $\ul X$ via evaluation at the marked points. This
construction then leads to our third main result.

\begin{theorem}
Assume that $X$ is smooth over $S$. Then under the assumptions of
Theorem~\ref{Main Thm II} there exists a virtual fundamental class
$\lfor \MM(X/S,\beta)\rfor$ leading to log Gromov-Witten invariants
with the expected properties.
\qed
\end{theorem}

Among the expected properties are compatibility with base change,
equality with the ordinary virtual fundamental class for trivial log
structures and with the usual fundamental class in unobstructed
situations, and many more that follow by strict analogy with the
ordinary case.

While this work was in progress we learnt of ongoing work of
Abramovich, Chen and coworkers aiming at carrying out our program in
the special case that there exists a surjection $\ul \NN^r\to
\ol\M_X$ (\cite{chen},\cite{abramovich etal},\cite{AbramovichChen}).
\medskip

\emph{Conventions.} We work in the categories of schemes of finite
type over a field $\kk$ of characteristic~$0$ and of fine saturated
log schemes \cite{Kato 1989} over a fixed base log scheme $S$,
itself fine and saturated. Our standard notation for log schemes is
$X=(\ul X,\M_X)$, but we do not underline ordinary schemes that we
do not want to endow with a log structure unless there is a chance
of confusion. Similarly, for a morphism of log schemes the notation
is $f=(\ul f,f^\flat)$. The \emph{ghost sheaf} (also called
\emph{characteristic}) of a log structure $\M$ is denoted
$\ol\M:=\M/\O^\times$. We use multiplicative notation for $\shM$ and
additive notation for $\overline\M$. Throughout $X$ is a log scheme
over $S$ such that $\M_S$, $\M_X$ and the structure morphism $(\ul
X,\M_X)\to (\ul S,\M_S)$ are defined in the Zariski topology. To
have a good theory of ordinary stable maps at our disposal we
require $\ul X$ is quasi-projective over $\ul S$. If $Y$ is a scheme
(algebraic space, stack) then $|Y|$ denotes the set of geometric
points of $Y$ endowed with the Zariski topology, see \cite{knutson},
II.6 or \cite{champs}, Ch.5. We use overlined symbols to denote
elements of $|Y|$; in particular, if $y$ is a scheme-theoretic point
of $Y$ then $\ol y$ is a geometric point with image $y$ . A
\emph{toric monoid} is a fine, saturated, torsion-free monoid
without non-trivial invertibles. For a monoid $P$ we write
$P^\vee:=\Hom (P,\NN)$ for the dual in the category of monoids,
$P^\gp$ for the associated abelian group and $P^*:=\Hom(P,\ZZ)$.
\medskip

\emph{Acknowledgements}. We thank M.~Olsson and R.~Pandharipande
for valuable discussions, and D.~Abramovich, Q.~Chen and  B.~Kim for their
interest in this work and for keeping us updated on the progress of
their work. We also thank one of the referees for a very careful
reading of the paper and several valuable comments.

\section{Stable log maps}
\label{Sect: Stable log maps}

\subsection{Log smooth curves}
\label{par: Stable log curves}
In this subsection we work in the absolute situation $S=\Spec\kk$
with the trivial log structure. A logarithmic version of the theory
of (pre-) stable curves of Deligne, Mumford and Knudsen
\cite{delignemumford}, \cite{knudsen} has been developed by F.Kato
\cite{Kato 2000}, see also \cite{OlssonThesis}. The starting point
is the following result on the structure of log smooth curves.

\begin{theorem}\cite[p.222]{Kato 2000}
\label{Thm: structure of log curves}
Let $\pi:C\to W$ be a smooth and integral morphism of fine saturated
log schemes such that every geometric fibre is a reduced curve.
Assume that $\ul W=\Spec A$ for $(A,\maxid)$ a strictly
Henselian local ring. Let $0\in \ul W$ be the closed point,
$Q=\ol{\M}_{W,0}$ and $\sigma: Q\to A$ a chart for the log
structure on $W$. Then \'etale locally $C$ is isomorphic to one of
the following log schemes $V$ over $W$.
\begin{enumerate}
\item[(i)]
$\Spec(A[z])$ with the log structure induced from the homomorphism
\[
Q\lra \O_V,\quad q\longmapsto \sigma(q).
\]
\item[(ii)]
$\Spec(A[z])$ with the log structure induced from the homomorphism
\[
Q\oplus\NN\lra \O_V,\quad (q,a)\longmapsto z^a \sigma(q).
\]
\item[(iii)]
$\Spec(A[z,w]/(zw-t))$ with $t\in\maxid$ and with the log structure
induced from the homomorphism
\[
Q\oplus_\NN \NN^2\lra \O_V,\quad
\big(q,(a,b)\big)\longmapsto \sigma(q) z^a w^b.
\]
Here $\NN\to\NN^2$ is the diagonal embedding and $\NN\to Q$,
$1\mapsto\rho_q$ is some homomorphism uniquely defined by $C\to W$.
Moreover, $\rho_q\neq 0$.
\end{enumerate}
In this list, the morphism $C\to W$ is represented by the
canonical maps of charts $Q\to Q$, $Q\to Q\oplus\NN$
and $Q\to Q\oplus_\NN \NN^2$, respectively. 
\qed
\end{theorem}

Case (i) deals with smooth points of the central fibre where the log
structure comes entirely from the base while (iii) covers the
situation in a neighbourhood of the degeneracy locus of $\pi$. To
interpret case (ii) observe that the preimage of $z=0$ defines a
section of $\pi$ on an \'etale neighbourhood of $x$. This should be
viewed as the section defining a marked point. Let $\Gamma\subset
\ul C$ be the image of this section and $\iota: \ul C\setminus\Gamma
\hookrightarrow \ul C$ the inclusion. Then the log structure of $C$
near $x$ is the sum of the log structure of the base and the log
structure $\iota_*(\O_{\ul C\setminus\Gamma}^\times)\cap \O_{\ul C}$
associated to $\Gamma$.

\begin{remark}\label{Rem: ghost sheaves of log curves}
(1)\ The monoid $Q$ together with $\rho_q\in Q$ for the nodes $q\in
\ul C$ determine the log structure on the closed fibre $\ul C_0$ and the
morphism to the log point $(\Spec(A/\maxid), Q)$ on the level of
ghost sheaves. In fact, let $\ol\eta$ be the generic point of the
branch of $\ul C_0$ defined by $w=0$ in (iii). Then at
$\ol\eta$ the element $z^a$ becomes invertible. Hence the
compatibility of the charts in~(i) and (iii) implies that the
generization map $\ol\M_{C_0,\ol q}\lra
\ol\M_{C_0,\ol\eta}$ equals
\[
Q\oplus_\NN \NN^2\lra Q,\quad
(m,(a,b))\longmapsto m+b\cdot\rho_q. 
\]
These generization maps together with one copy of $\NN$ at each of
the special points in~(ii) define $\ol\M_{C_0}$ and the structure
homomorphism $Q\to\Gamma(\ul C_0, \ol\M_{C_0})$
uniquely up to unique isomorphism.

Note that if $\ul C$ has non-normal irreducible components then
$\M_C$ is only defined over the \'etale site. Still, $\ol\M_C$
can be described completely by generization maps, the only
difference now being that at a node $q\in \ul C$ in the closure of
only one generic point $\eta$ there are two generization maps
$\ol\M_{C,\ol q}\to \ol\M_{C,\ol\eta}$.\\[1ex]
2)\ At a node $q\in\ul C_0$ the two generization maps
$Q\oplus_\NN\NN^2\to Q$ define the homomorphism
\[
\iota:Q\oplus_\NN\NN^2\lra Q\times Q,\quad
(m,(a,b))\longmapsto (m+a\cdot\rho_q,m+b\cdot\rho_q).
\]
If $\iota\big( (m,(a,b))\big)=0$ then $m+a\rho_q=m+b\rho_q=0$, and
hence $a=b$ since $\rho_q\neq0$. In view of the definition of
$Q\oplus_\NN \NN^2$ this implies
\[
(m,(a,b))=(m,(a,a))=\big(m+a\rho_q,(0,0)\big)=0.
\]
Thus $\iota$ exhibits $Q\oplus_\NN\NN^2$ as a submonoid of $Q\times
Q$:
\begin{equation}\label{Q oplus_NN NN^2}
Q\oplus_\NN\NN^2\simeq \big\{(m_1,m_2) \in Q\times Q\,\big|\, 
m_1-m_2 \in \ZZ\rho_q\text{ in }Q^\gp\big\}\subset Q\times Q.
\end{equation}
In fact, if $m_1-m_2=\alpha\rho_q$ then
\[
(m_1,m_2)= \begin{cases}
\iota\big( m_2,(\alpha,0)\big),& \alpha\ge 0\\
\iota\big( m_1,(0,-\alpha)\big),&\alpha<0.
\end{cases}
\]
\end{remark}

The theorem suggests the following generalization of the notion of
marked (pre-) stable curves to log geometry.

\begin{definition}
\label{Def: Prestable log-curve}
A \emph{pre-stable (marked) log curve} over $W$ is a pair $(C/W,
\mathbf{x})$ consisting of a proper log smooth and integral morphism
$\pi: C\to W$ of fine saturated log schemes over $S$ together with a
tuple of sections $\mathbf{x}=(x_1,\ldots,x_l)$ of $\ul \pi$, such that
every geometric fibre of $\pi$ is a reduced and connected curve, and
if $U\subset \ul C$ is the non-critical locus of $\pi$ then
$\ol \M_C|_U\simeq \ul\pi^* \ol\M_W\oplus \bigoplus_i
{x_i}_* \NN_{\ul W}$.

A pre-stable log curve is \emph{stable} if forgetting the log
structure leads to an ordinary stable curve.
\end{definition}

\begin{remark}
1)\ \ The underlying morphism of a pre-stable log curve is flat and its
geometric fibres have at most ordinary double points. The
underlying morphism of schemes of a pre-stable log curve is hence an
ordinary pre-stable curve.\\[1ex]
2)\ \ The condition on $\ol\M_C$ says that the sections $x_i$ label
precisely the special non-nodal points occurring in
Theorem~\ref{Thm: structure of log curves},(ii). A straightforward
generalization would only label a subset of the special non-nodal
points.
\end{remark}

We sometimes deal with diagrams of spaces endowed with fine
sheaves of monoids that look like they are the diagrams of ghost
sheaves of a pre-stable log curve. For later reference we cast this
situation into the following definition.

\begin{definition}
\label{Def: pre-stable at ghost sheaf level}
Let $(\ul \pi:\ul C\to\ul W, \mathbf{x})$ be an ordinary pre-stable
curve over a scheme $\ul W$ and let $\ul\pi^*\ol\M_W\to\ol\M_C$ be a
morphism of fine sheaves of monoids on $\ul C$ and $\ul W$,
respectively. We say that $\big( (\ul C,\ol\M_C)/(\ul
W,\ol\M_W),\mathbf{x})$ has the structure of a pre-stable log curve
\emph{on the level of ghost sheaves} if for every geometric point
$\ol w:\Spec\kappa\to \ul W$ the situation is described by
Remark~\ref{Rem: ghost sheaves of log curves},1. In particular,
writing $Q=\ol\M_{W,\ol w}$, then for each node $q\in\ul C_{\ol w}$
there is an element $\rho_q\in Q$ such that $\ol M_W|_{\ul C_{\ol
w}}$ is isomorphic to the sheaf with stalks $Q\oplus_\NN \NN^2$
(with $1\in\NN$ mapping to $\rho_q\in Q$ and to $(1,1)\in\NN^2$) at
$q$, with stalks $Q\oplus\NN$ at the marked points, and with stalks
$Q$ at all other points; the generization maps and the morphism
$\ol\M_{W,\ol w}\to \ol\M_C|_{\ul C_{\ol w}}$ are as described in
Remark~\ref{Rem: ghost sheaves of log curves},1.
\end{definition}

F.~Kato also introduces the notion of basic log structure of a marked
\emph{stable} curve (\cite{Kato 2000} p.227f). A simple way to think
about this concept is as follows. Let $(\ul\pi: \ul C\to \ul W, \mathbf{x}=
(x_1,\ldots,x_k))$ be a stable marked curve. Locally with respect to
the base, $(\ul C/\ul W,\mathbf{x})$ is the pull-back by a morphism
$\ul W\to \ul T$ of a pre-stable marked curve $(\ul q:\ul U\to \ul T,
\mathbf{y})$ that is formally versal at any point of $\ul T$.
Versality implies that the image of the subspace of $\ul U$ defined
by the first Fitting ideal of $\Omega^1_{\ul U/\ul T}$ is a normal
crossings divisor $D\subset \ul T$. This divisor is the scheme
theoretic version of the subset of $\ul T$ parametrizing singular
curves. We endow $\ul T$ and $\ul U$ with the log structures
associated to the divisors $D$ in $\ul T$ and $q^{-1}(D)$ and the
divisor of marked points in $\ul U$, respectively. The basic log
structure on $\ul C\to \ul W$ is then obtained by pull-back via $\ul
W\to \ul T$. Globally one obtains an \'etale descent datum for
$\M_C$, $\M_W$ and for the morphism $\pi^*\M_W\to \M_C$. Note also
that at a geometric point $\ol w$ of $\ul W$ there is an isomorphism
$\ol{\M}_{W,\ol w} \simeq\NN^r$ where $r$ is the number of double
points of $\ul C_{\ol w}$, which by versality equals the number of
branches of $D$ at the image of $\ol w$ in $\ul T$.

The existence of basic log structures on stable curves makes it
possible to endow the stacks $\Mbf_{g,k}$ of $l$-marked
stable curves with a logarithmic structure \cite[p.230f]{Kato 2000}.
This means (\cite{OlssonENS}, Definition~5.1 and Corollary~5.8) that
there exists a factorization
\begin{equation}\label{Mgk is a log stack}
\Mbf_{g,k}\lra (\mathrm{Log})\lra \Sch
\end{equation}
of the functor defining the stack $\Mbf_{g,k}$. Here
$(\mathrm{Log})$ is the category of fine saturated log schemes with
strict morphisms. The first arrow maps a family $(\ul C/ \ul
W,\ul{\mathbf{x}})$ of $k$-marked stable curves of genus $g$ to its
base scheme $\ul W$ endowed with the basic log structure. The
factorization \eqref{Mgk is a log stack} endows $\Mbf_{g,k}$ with a
log structure, defining a log stack $\scrM_{g,k}$. Of course, not
every stable log curve carries the basic log structure of the
underlying stable marked curve. The log structure can rather always
be obtained from the basic log structure by a unique base change
inducing the identity on the underlying spaces (\cite{Kato 2000},
Proposition~2.1). See also \cite{OlssonThesis}, Ch.5, for an
extended treatment.

Analogous statements hold for pre-stable curves, leading to the log
algebraic stack $\scrM$, an Artin stack (see Appendix~\ref{App:
Prestable curves}).

\subsection{Stable log maps}
\label{par: Stable log maps}
We now turn to the main concept of this paper, a logarithmic
version of the notion of stable map. Recall that $X$ is a log scheme
over $S$, with log structures defined on the Zariski sites.

\begin{definition}\label{Def: Stable log map}
A \emph{log curve over $X$} with \emph{base} $W$ is a pre-stable
marked log curve $(C/W,\mathbf{x})$ (Definition~\ref{Def:
Prestable log-curve}) together with a morphism $f:C
\to X$ fitting into a commutative diagram of log schemes
\begin{equation}\label{Fig: Stable log map}
\begin{CD}
C @>f>> X\\
@V\pi VV @VVV\\
W@>>> S\
\end{CD}
\end{equation}
A log curve over $X$ is a \emph{stable log map} if for every
geometric point $\ol w\to \ul W$ the restriction of $\ul f$ to
the underlying pre-stable marked curve belonging to $\ul
C_{\ol w}\to \ol w$ is an ordinary stable map.
The notation is $(C/W,\mathbf{x},f)$ with the morphisms $\pi:C\to W$
and $W\to S$ usually understood.

A morphism of stable log maps
\[
\Phi: (C_1/W_1, \mathbf{x}_1,f_1)\lra
(C_2/W_2, \mathbf{x}_2,f_2)
\]
is a cartesian diagram of log curves $\Phi:C_1/W_1\to
C_2/W_2$ over $S$ with $W_1\to W_2$ \emph{strict} and such that
$f_1=f_2\circ\Phi$ and $\mathbf{x}_2= \ul \Phi\circ \mathbf{x}_1$. 

The category of stable log maps thus obtained is denoted
$\tilde\MM(X/S)$ or just $\tilde\MM(X)$.
\end{definition}

Let $\kappa$ be a field. For a toric monoid $Q$ with
$Q^\times=\{0\}$ denote by $ (\Spec \kappa,Q)$ the associated
logarithmic point, that is, $\Spec \kappa$ with log structure
\[
Q\times \kappa^\times\lra \kappa,\quad
(q,a)\longmapsto \begin{cases}0,&q\neq0\\ a,&q=0.\end{cases}
\]

Now let $\ul C$ be a pre-stable curve over $\kappa$, and assume
given a fine saturated log structure $\alpha:\M\to \O_{\ul C}$ on
the Zariski site of $\ul C$. For the application to stable log maps
we will only consider $\M=\ul f^*\M_X$ for some $\ul
f:\ul C\to \ul X$, but this is not relevant for the following
discussion. In any case, this log structure $\M$ can be rather
arbitrary and certainly does not need to be smooth for any log
structure on $\Spec \kappa$. We need to understand diagrams
\begin{equation}\label{Diag: Universal diagram}
\parbox{0pt}{\xymatrix{
(\ul C,\M_C)\ar[d]\ar[r] &(\ul C,\M)\\
(\Spec \kappa,Q)
}}
\end{equation}
with left-hand vertical arrow smooth, for some toric monoid $Q$. In
the simplest case of trivial $\M$ we are in the situation recalled
in Section~\ref{par: Stable log curves}. Thus in this case there is
a universal such diagram with $Q=\NN^r$ and $r$ the number of nodes
of $C$, in the sense that any other diagram is obtained by unique
pull-back. We will see that there is a similarly universal object in
complete generality. The essential step is the characterization of a
universal such diagram on the level of ghost sheaves. This is the
object of the next subsection.

\subsection{The category {\protect $\mbox{\underline{\rm{GS}}}$}}
\label{par: GS}

Let $\ol\M$ be a fine saturated sheaf on a pre-stable curve
$(\ul\pi:\ul C\to \ul W,\mathbf x)$. We consider the following category.

\begin{definition}\label{Def: GS}
Let $\GS(\ol\M)$ (for ``ghost sheaves'') be the category with
objects
\[
(\ol\M_W, \ol\M_C, \psi: \ul\pi^*\ol\M_W\to\ol\M_C,
\varphi: \ol\M\to \ol \M_C),
\]
where $\ol\M_W$ and $\ol\M_C$ are fine sheaves on $\ul W$ and $\ul
C$, respectively, and $\psi$ and $\varphi$ are local\footnote{A
homomorphism of monoids $\varphi: P\to Q$ is called \emph{local} if
$\varphi^{-1}(Q^\times) =P^\times$ or, equivalently,
$\varphi^{-1}(Q\setminus Q^\times)= P\setminus P^\times$. Thus for
the case of toric monoids this means $\varphi^{-1}(0)=0$.}
homomorphisms of fine
(sheaves of) monoids. We require that $\psi$  endows $(\ul C/\ul W,
\mathbf x)$ with the structure of a pre-stable log curve over $(\ul
W,\ol\M_W)$ on the level of ghost sheaves (Definition~\ref{Def:
pre-stable at ghost sheaf level}). Thus the objects of $\GS(\ol\M)$
can alternatively be taken as diagrams
\begin{equation}\label{Diag: category BC}
\parbox{0pt}{\xymatrix{
\ol \M_C&\ar[l]_{\varphi} \ol\M\\
\ar[u]^\psi \ul\pi^*\ol\M_W
}}
\end{equation}
of fine saturated sheaves on $\ul C$ with local homomorphisms. A
morphism
\[
(\ol\M_{W,1}, \ol\M_{C,1}, \psi_1,\varphi_1)\lra
(\ol\M_{W,2}, \ol\M_{C,2}, \psi_2, \varphi_2)
\]
in $\GS(\ol\M)$ is given by a pair of homomorphisms $\ol\M_{W,1}\to
\ol\M_{W,2}$ and $\ol\M_{C,1}\to \ol\M_{C,2}$ with the obvious
compatibilities with $\psi_i$ and $\varphi_i$, $i=1,2$. 
\end{definition}

For the remainder of this subsection we restrict to the case $\ul
W=\Spec\kappa$ for a field $\kappa$. In this case we write $Q$ for
the stalk of $\ol\M_W$ and $\psi:Q\to \Gamma(\ul C,\ol\M_C)$. One
central insight in this paper is the characterization of the
connected components of $\GS(\ol\M)$ in this case by what we call
the \emph{type} of a stable log map. Even more usefully, for each
type we construct a universal object of the corresponding connected
component (see Proposition~\ref{Prop: Universal property of
basicness} below). Saying that $\GS(\ol\M)$ has a universal object
essentially means that any two log enhancements of an ordinary
stable map fit into one family, at least on the level of ghost
sheaves.

\begin{discussion}
\label{Disc: Stalkwise}
To introduce the concept of type let us reformulate
Diagram~\eqref{Diag: category BC} on the level of stalks.
By the structure of stable log curves over $(\Spec\kappa, Q)$ there are
three types of points $x$ on $\ul C$, depending on the stalks of
$\ol \M_C$, as follows. Write $P_x:=
\ol\M_x$.\footnote{Recall that we assumed $\ol\M$ is a
sheaf on the Zariski site, so $\ol\M_x= \ol\M_{\bar
x}$.}
\begin{enumerate}
\item[(i)]
$x=\eta$ is a generic point or a general closed point.
Then $\ol \M_{C,\ol\eta}=Q$ and $\varphi$ defines a
homomorphism
\[
\varphi_{\ol \eta}: P_\eta\lra Q.
\]
\item[(ii)]
$x=p$ is a marked point. Then $\ol \M_{C,\ol p} =Q\oplus\NN$
with $\psi_{\ol p}$ inducing the inclusion of $Q$ as first factor.
Then $\varphi_{\ol p}$ is determined by $\varphi_{\ol \eta}$ for
$\ol\eta$ the generic point of the irreducible component containing
$p$ together with
\begin{equation}\label{u_p}
u_p:=\pr_2\circ \varphi_{\ol p}: P_p\lra \NN.
\end{equation}
\item[(iii)]
$x=q$ is a node. Then $\ol\M_{C,\ol q} \simeq
Q\oplus_\NN\NN^2$ with $\NN\to Q$, $1\mapsto\rho_q$ and
$\NN\to\NN^2$, $1\mapsto(1,1)$. Let $\ol\eta_1$, $\ol\eta_2$ be
the geometric generic points of the branches of $\ul C$ at $\ol q$.
We have a commutative diagram
\[
\xymatrix@C=20pt
{&P_{\eta_1}\ar[rr]^{\varphi_{\ol\eta_1}}&&Q\\
P_q\ar[rr]^{\varphi_{\ol q}}\ar[ru]^{\chi_1}\ar[rd]_{\chi_2}
&&Q\oplus_\NN\NN^2 \ar[ru]\ar[rd]\ar[r]^\iota
&\ar[u]_{\pr_1} \ar[d]^{\pr_2} Q\times Q\\
&P_{\eta_2}\ar[rr]^{\varphi_{\ol\eta_2}}&&Q
}
\]
where the diagonal arrows are generization maps. Recall from
\eqref{Q oplus_NN NN^2} that $\iota$ is injective with $(m_1,m_2)\in
Q\times Q$ in the image of $Q\oplus_\NN \NN^2$ iff
$m_1-m_2\in\ZZ\rho_q$, viewed as equation in $Q^\gp$. Thus
$\varphi_{\ol q}$ is determined uniquely by $\varphi_{\ol\eta_1}$,
$\varphi_{\ol\eta_2}$, and there exists a homomorphism
\begin{equation}\label{u_q}
u_q: P_q\lra \ZZ
\end{equation}
defined (since $\rho_q\neq 0$) by the equation
\begin{equation}\label{eqn at node}
\varphi_{\ol\eta_2}\big(\chi_2(m)\big)- 
\varphi_{\ol\eta_1}\big(\chi_1(m)\big)=
u_q(m)\cdot\rho_q.
\end{equation}
\end{enumerate}

Note that the definition of $u_q$ in (iii) depends on the choice of
an ordering of the adjacent branches of $\ul C$. Such a choice can be
implemented by orienting the \emph{dual intersection graph}
$\Gamma_{\ul C}$ of $\ul C$. The graph $\Gamma_{\ul C}$ has one
vertex $v_\eta$ for each generic point $\eta\in \ul C$ and an edge
$E_q$ for each node $q$ joining $v_{\eta_1}$ and $v_{\eta_2}$ for
$\eta_1,\eta_2$ the generic points of the two adjacent branches at
$q$. Note that we can have $\eta_1=\eta_2$ if the component is not
normal. In addition there is an unbounded edge (a flag) $E_p$ for
each marked point $p$, with adjacent vertex $v_\eta$ for the unique
generic point $\eta$ with $p\in\cl(\eta)$, the topological closure
of $\eta$.
\end{discussion}

\begin{remark}
\label{Rem: structure of M}
In the above discussion, we only needed to use the stalks $P_x$ of
$\overline\M$ at points of $\ul C$ which were either generic, marked
or double points. In fact, for any closed point $x\in \ul C$ which
is neither a double point or marked point, we must have
$P_x=P_{\eta}$, where $\eta$ is the generic point of the component
of $C$ containing $x$. Indeed, we have a surjective generization map
$\chi:P_x\rightarrow P_{\eta}$. Since $\ol\M$ is the ghost sheaf of
a fine log structure, $\chi$ is given by localization at a face of
$P_x$ followed by dividing out the submonoid of invertible elements.
Thus if $\chi$ is not an isomorphism we can always find $0\not= m\in
P_x$ with $\chi(m)=0$. But $\overline\M_{C,\bar
x}=\overline\M_{C,\ol\eta}=Q$, and since generization is compatible
with the map $\varphi$ on stalks, we have
$\varphi_{\ol\eta}\circ\chi = \varphi_{\ol x}$, hence
$\varphi_{\ol x}(m)=0\in Q$. This contradicts the property that
$\varphi_{\ol x}$ is local. 
\end{remark}

\begin{definition}\label{Def: types in GS}
1)\ The \emph{type} of an object $(Q,\ol\M_C, \psi,\varphi)$ of the
category $\GS(\ol\M)$ is the set $\mathbf u:=\big\{ u_p\in
P_p^\vee,u_q \in P_q^*\big\}$, where $u_p$ and $u_q$ are defined in
\eqref{u_p} and \eqref{u_q}, respectively. Here $p$ and $q$ run over the
marked and nodal points of $\ul C$, respectively. Given $\mathbf u$,
the full subcategory of $\GS(\ol\M)$ with objects of type $\mathbf
u$ is denoted $\GS(\ol\M,\mathbf u)$. \\[1ex]
2)\ The \emph{type} $(\Gamma_{\ul C},\mathbf u)$ of a stable log map
over a log point $(C/(\Spec\kappa,Q), \mathbf{x}, f)$ is the dual
intersection graph $\Gamma_{\ul C}$ of $\ul C$ together with the
type $\mathbf u$ of the corresponding object of $\GS(\ol\M)$.
\end{definition}

The type is compatible with generization:

\begin{lemma}\label{Lem: type under generization}
Let $\mathbf u$, $\mathbf u'$ be the types of a stable log map
$(C/W,\mathbf x,f)$ at two geometric points $\ol w\to \ul W$, $\ol
w'\to \ul W$  with $\ol w\in\cl(\ol w')$. For $x\in \ul C_{\ol w}$,
$x'\in \ul C_{\ol w'}$ and $x\in\cl(x')$ let $\chi_{x'x}: P_x \to
P_{x'}$ be the generization map of the stalks of $\ul f^* \ol\M_X$.
Then for marked or nodal points $x,x'$ with $x\in\cl(x')$ it holds
\[
u_x= u_{x'}\circ \chi_{x'x}.
\]
\end{lemma}

\proof
For marked points $p,p'$ this follows readily from compatibility of
\eqref{u_p} with generization. For nodal points $q,q'$ comparing the
generization of \eqref{eqn at node} for $q$ with the equation for $q'$
yields
\[
u_q\cdot \kappa(\rho_q)= (u_{q'}\circ\chi_{q'q})\cdot\rho_{q'}.
\]
where $\kappa: \ol\M_{W,\ol w}\to \ol\M_{W,\ol w'}$.
The claimed equation now follows from $\kappa(\rho_q)=\rho_{q'}\neq0$.
\qed

\subsection{The standard log point and tropical curves}
\label{par: Log maps over std log pt}
An interesting special case is stable log maps over standard log
points $(\Spec\kappa,\NN)$. This provides the connection to tropical
geometry. In toric situations this connection has previously been
discussed in \cite{nisi} and in \cite{gross_seattle},\S10.
Another motivation is that this case suffices to
characterize universal stable log maps. To explain this we
consider the situation of Diagram~\ref{Diag: Universal
diagram} of a pre-stable curve $\ul C/\Spec\kappa$ and a fine
saturated log structure $\M$ on $\ul C$. We think of the case
$\M= \ul f^* \M_X$ for an ordinary stable
map $(\ul C/\Spec\kappa,\mathbf{x}, \ul f)$. Now if
Diagram~\ref{Diag: Universal diagram} is universal (for a fixed
type) then diagrams of the same form over the standard log point are
given by morphisms
\[
(\Spec\kappa,\NN)\lra (\Spec\kappa,Q).
\]
Moreover, two such morphisms lead to isomorphic log maps if and only
if they differ by a homomorphism $Q\to \kappa^\times$.
Now a morphism of log structures $Q\times\kappa^\times\lra
\NN\times\kappa^\times$ has the form
\[
(m,a)\longmapsto \big(\varphi(m),h(m)\cdot a\big)
\]
for some $\varphi\in\Hom(Q,\NN)$ with $\varphi^{-1}(0) =\{0\}$
and $h\in\Hom(Q,\kappa^\times)$. Composing with the automorphism
\[
(m,a)\longmapsto \big(m,h(m)\cdot a\big)
\]
of $(\Spec\kappa, Q)$ we may assume $h=1$. Hence the set of
isomorphism classes of stable log maps over $(\Spec \kappa,\NN)$
obtained from (\ref{Diag: Universal diagram}) by base change is in
one-to-one correspondence with
\[
\big\{ \varphi\in \Hom(Q,\NN) \,\big|\, \varphi^{-1}(0)=\{0\}\big\}.
\]

The upshot of this discussion is that $\Int(Q^\vee)$ is equal to the
set of isomorphism classes of Diagrams~\ref{Diag: Universal diagram}
with $Q=\NN$. Note that while this is a discrete set, the set of
isomorphism classes of Diagrams~\ref{Diag: Universal diagram}
relative to a \emph{fixed} log point $(\Spec\kappa, Q)$ is
fibred over this discrete set with fibres $\Hom(Q,\kappa^\times)$.

Let us now specialize Discussion~\ref{Disc: Stalkwise} to $Q=\NN$.
Then for a node $q$ the element $\rho_q\in Q\setminus\{0\}$ is a
number $e_q\in\NN\setminus\{0\}$, and $Q\oplus_\NN \NN^2$ is
isomorphic to the submonoid $S_{e_q}$ of $\NN^2$ generated by
$(e_q,0)$, $(0,e_q)$, $(1,1)$.\footnote{In terms of
generators and relations we have $S_{e_q}=\langle a_1,a_2,a_3\, |\,
a_1+a_2=e_q\cdot a_3\rangle$, where $a_1=(e_q,0)$, $a_2=(0,e_q)$,
$a_3=(1,1)$.} Thus specifying the sheaf $\ol\M_C$
together with a homomorphism $\ol\M\to\ol\M_C$ is equivalent to the
following data, for points $x\in\ul C$:
\begin{enumerate}
\item[(i)]
$x=\eta$ is a generic point. Then
$\varphi_{\ol \eta}: P_\eta\lra\NN$ defines an element $V_\eta\in
P_\eta^\vee$.
\item[(ii)]
$x=p$ is a marked point. As in the general case, $\varphi_p$ is
determined by (i) and by $u_p\in P_p^\vee$, fixed by the type.
\item[(iii)]
$x=q$ is a node and $\ol\eta_1$, $\ol\eta_2$ are the geometric
generic points of the adjacent branches of $\ul C$. Letting
$i_{q,\eta_i}: P_{\eta_i}^\vee\to P_q^\vee$ be the inclusion induced
by the generization maps, Equation~\ref{eqn at node} now reads
\begin{equation}\label{Eq: basic relation}
i_{q,\eta_2}(V_{\eta_2})- i_{q,\eta_1} (V_{\eta_1})= e_q u_q.
\end{equation}
\end{enumerate}
Thus apart from the type $\mathbf u$ a stable log map over a
standard log point defines points $V_\eta\in P_\eta^\vee$ and
$e_q\in\NN\setminus\{0\}$. We call the tuple $\big((V_\eta)_\eta,
(e_q)_q\big)$ \emph{tropical data} of a stable log map over a
standard log point. Similarly we can talk of tropical data for an
object $(Q,\ol\M_C,\psi, \varphi)$ of $\GS(\ol\M)$ with $Q=\NN$.
\medskip

To discuss the relationship with tropical geometry we recall the
basic definition of a tropical curve. By abuse of notation we
confuse a graph and its topological realization.

\begin{definition}\label{Def: tropical curve}
Let $\ol\Gamma$ be a connected graph and let $\Gamma$ be the
topological space obtained by removing from $\ol\Gamma$ a set of
univalent vertices of $\ol\Gamma$, so that $\Gamma$ has
both compact and non-compact edges.\footnote{The set of univalent
vertices to be removed is part of the data. A traditional tropical
curve in $\RR^n$ is the image of $\Gamma$, and $\overline\Gamma$ is
obtained by attaching a univalent vertex at each unbounded edge.}
We assume that $\Gamma$
has at least one vertex. Let $N\simeq\ZZ^n$ be a lattice,
$N_{\RR}:=N\otimes_{\ZZ}\RR$. A \emph{tropical curve in $N_{\RR}$
with domain $\Gamma$} consists of the following data:
\begin{enumerate}
\item[(i)]
For each flag $(v,E)$ of $\Gamma$, where $v$ is a vertex of $\Gamma$
and $E$ an edge containing $v$, we are given a \emph{weight vector}
$u_{(v,E)}\in N$. If $E$ has two vertices $v_1$ and $v_2$, then
$u_{(v_1,E)}=-u_{(v_2,E)}$, and if $E$ is a loop, then
$u_{(v,E)}=0$. The graph $\Gamma$ along with the weight vectors is
called the \emph{type} of the tropical curve.
\item[(ii)]
 A map $h:\Gamma\to N_{\RR}$ with the following properties:
\begin{enumerate}
\item
For any edge $E$ of $\Gamma$ with vertex $v$, $h|_E$ is constant
if $u_{(v,E)}=0$, and otherwise $h|_E$ is proper and identifies $E$ with
an affine line segment or ray. Furthermore, $u_{(v,E)}$ is a tangent
vector to $h(E)$ pointing away from $h(v)$.
\item
For each vertex $v$, we have the \emph{balancing condition}
\[
\sum_E u_{(v,E)}=0,
\]
where the sum is over all edges $E$ with vertex $v$.
\end{enumerate}
\end{enumerate}
\end{definition}

\begin{discussion}
\label{Disc: Tropical curve}
The term ``tropical data'' is motivated by the case that $\ol\M^\gp$
is globally generated. In this case the tropical data gives rise to
a generalized tropical curve with domain $|\Gamma_{\ul C}|$, the
geometric realization of $\Gamma_{\ul C}$, and image in $N_\RR$ for
\[
N:=\Hom(\Gamma(\ul C,\ol\M^\gp),\ZZ),
\]
as follows. The generalization concerns the balancing condition, as to be
discussed. The restriction maps $\Gamma(\ul C,\ol\M^\gp)\to \ol\M_{\ol
x}^\gp$ induce injections
\[
P_x^\vee\lra N.
\]
Denote by $\tilde V_\eta, \tilde u_p, \tilde u_q$ the images of
$V_\eta, u_p, u_q$ in $N$. Define a continuous map
\[
h:|\Gamma_{\ul C}|\lra N_\RR
\]
by sending $v_\eta$ to $\tilde V_\eta$, an edge $E_q$ with adjacent
vertices $v_{\eta_1}, v_{\eta_2}$ to the line segment connecting
$\tilde V_{\eta_1}$ and $\tilde V_{\eta_2}$, and an unbounded edge
$E_p$ with adjacent vertex $v_\eta$ to the ray $\tilde V_\eta+
\RR_{\ge0} \tilde u_p$. The weight vectors are given by $u_{(v_\eta,
E_p)}:=\tilde u_p$ for the unbounded edges and $u_{(v_{\eta_i},
E_q)}:= \pm \tilde u_q$, with the sign chosen so that
$u_{(v_{\eta_i}, E_q)}$ points away from $\tilde V_{\eta_i}$. Note
that by \eqref{Eq: basic relation} for an edge $E_q$ with vertices
$v_{\eta_1}, v_{\eta_2}$ it holds $h(v_{\eta_2}) -h(v_{\eta_1})= \pm
e_q \tilde u_q$, so $e_q$ is the integral length of the
corresponding edge of the tropical curve, as a multiple of the
weight vector, at least for $\tilde u_q\neq 0$. 

As is, this does not in general fulfill the balancing condition
(Definition~\ref{Def: tropical curve},b). However, in
Proposition~\ref{Prop: Balancing condition} below we formulate a
modified balancing condition involving a correction term. The
correction term turns out to depend only on the given log structure
$\M$ on $\ul C$, that is, on the underlying ordinary stable map. The
balancing condition holds unmodified if for any irreducible component
$D\subset \ul C$ and $m\in\Gamma (D,\ol\M|_D)$ the degree of the
corresponding $\O_D^\times$-torsor $L_m\subset\M|_D$
vanishes\footnote{Recall that if $\M$ is a log structure on a space
$X$ then any section $m\in\Gamma(X,\ol\M)$ gives rise to an
$\O_X^\times$ torsor $\kappa^{-1}(m)\subset\M$ for
$\kappa:\M\to\ol\M$ the quotient homomorphism.}.
In general one can add one more unbounded edge at each vertex
$v_\eta$ with the weight vector $\tau_\eta^X$ derived from $\ul f$
via \eqref{tau_eta} to obtain an honest tropical curve in $\N_\RR$.

In the toric degenerations of toric varieties of \cite{nisi} already
$\ol\M_X^\gp$ is globally generated. It is then appropriate to
consider the composition with
\[
N\lra N':=\Hom\big(\Gamma(\ul X,\ol\M_X^\gp),\ZZ\big).
\]
Moreover, since in \cite{nisi} we work relative $(\Spec\kk,\NN)$
the images of the generator of $\ol \M_{\Spec\kk}=\NN$ under $X\to
(\Spec\kk,\NN)$ and under $(\Spec\kappa,\NN)\to (\Spec\kk,\NN)$
define a global section $\rho$ of $\ol\M_X$ and an element $b\in\NN=
\ol \M_{(\Spec\kappa,\NN)}$, respectively. Commutativity of
\eqref{Fig: Stable log map} at a generic point $\eta$ now implies
$V_\eta (\rho_{\ul f(\eta)}) =b$. Thus $h$ maps the vertices of $\Gamma$ to
the affine hyperplane $H\subset N'_\RR$ defined by
$\langle\rho,.\rangle=b$. Finally, from the structure of log smooth
curves at a marked point, it holds $\langle\rho, \tilde u_p\rangle
=u_p(\rho_p)=0$, and hence $\im(h)\subset H$. It is the map to
$H$ that traditionally and in \cite{nisi} is called a tropical
curve. This ends our discussion of tropical curves in this context.
\end{discussion}

As suggested by the tropical curve interpretation of
Discussion~\ref{Disc: Tropical curve} there should be a balancing
condition at each vertex $v_\eta$ of $\Gamma_{\ul C}$ imposing
restrictions on $u_p$, $u_q$ for the adjacent edges $E_p$, $E_q$.
Denote $\M:=\ul f^*\M_X$. For a generic point $\eta\in \ul C$ let
$D:=\cl(\eta)$ and $g:\tilde D\to \ul C$ the normalization of
$D$. This gives rise to maps which are compositions
\begin{equation}\label{tau_eta}
\begin{aligned}
\tau_\eta^X:&&\Gamma(\tilde D,g^*\ol{\M})&&\lra&&\Pic
\tilde D&&\mapright{\deg}&& \ZZ\\
\tau_\eta^C:&&\Gamma(\tilde D,g^*\ol{\M}_C)&&\lra&&\Pic
\tilde D&&\mapright{\deg} &&\ZZ&.
\end{aligned}
\end{equation}
The first map associates to a section of $g^*\ol{\M}$ or
$g^*\ol\M_C$ the corresponding $\O_D^{\times}$-torsor,
and the second map is the degree homomorphism.

The balancing condition is due to the basic fact that $f^\flat$ must
induce isomorphisms of torsors, so that the pull-back $\varphi:
g^*\ol \M\to g^*\ol\M_C$ of $\ol{f^\flat}$ to $\tilde D$ fits
into the commutative diagram
\begin{equation} \label{taucomm}
\xymatrix@C=30pt
{
\Gamma(\tilde D,g^*\ol{\M})\ar[r]^{\varphi}\ar[rd]_{\tau^X_\eta}
&\Gamma(\tilde D,g^*\ol{\M}_C)\ar[d]^{\tau^C_\eta}\\
&\ZZ
}
\end{equation}

The map $\tau^X_\eta$ is given by $\ul f$ and $\M$, so we have no
control over it, except that if $\ul f$ contracts $D$, then
$\tau^X_\eta=0$. Otherwise there is nothing general we can say.

Similarly, $\tau^C_\eta$ is determined by $\M_C$. Explicitly,
for $q\in D$ identify $S_{e_q}\subset \NN^2$ with the
submonoid generated by $(0,e_q), (e_q,0), (1,1)$ so that the
generization map $\chi_{q}:\ol{\M}_{C,\ol q}=S_{e_q}
\to \ol{\M}_{D,\ol \eta}=\NN$ is the projection onto the
second coordinate: $\chi_{q}(a,b)=b$. We have the identification
\[
\Gamma\big(\tilde D,g^*\ol{\M}_C\big)=
\big\{(n_q)_{q\in  \tilde D}\,\big|\, \hbox{$n_q\in S_{e_q}$
and $\chi_{q}(n_q)=\chi_{q'}(n_{q'})$ for $q,q'\in  \tilde D$}\big\}
\oplus \bigoplus_{p\in  \tilde D} \NN.
\]
In particular, if
$\big((a_q,b_q),(n_p)\big)$ represents an element of
$\Gamma\big(\tilde D,g^*\ol{\M}_C\big)$ then all second entries
$b_q$ agree.

\begin{lemma}\label{tau_eta^C}
$\displaystyle\tau_\eta^C\big(((a_q,b)_{q\in \tilde D}, (n_p)_{p\in \tilde D})\big)=
-\sum_{p\in \tilde D} n_p+\sum_{q\in \tilde D} \frac{b-a_q}{e_q}$.
\end{lemma}
\proof
By log smoothness the element $\big((b,b)_q,(0)_p\big) \in
\Gamma(\tilde D,g^*\ol\M_C)$ with all $a_q=b$, $n_p=0$, maps
to the trivial $\O_{\tilde D}^\times$-torsor. Thus it suffices to consider
sections of $g^*\ol\M_C$ of the form $\big(
(a_q,0),(n_p)\big)$ with $a_q/e_q\in\ZZ$. Let $L\subset g^*\M_C$ be
the corresponding $\O_{\tilde D}^\times$-torsor.  The structure map
$L\to g^*\M_C\to \O_{\tilde D}$ identifies $L$ with the sheaf of
regular functions on $\tilde D$ with zeros of order $a_q/e_q$ at $q$
and of order $n_p$ at $p$. In fact, if the log structure at $q$ is
induced from the toric model $\Spec\kappa[x,y,t]/ (xy-t^{e_q})$ with
$\tilde D$ corresponding to $V(y)$, then $x$ defines an element of
$\M_{C,\ol q}$ mapping to $(e_q,0)\in S_{e_q}$. Hence
\[
\deg(L)=\deg \O_{\tilde D}\Big(-\sum_q \frac{a_q}{e_q} q-\sum_p n_p p\Big)
=-\sum_q \frac{a_q}{e_q}-\sum_p n_p. 
\]
This is the claimed formula.
\qed
\medskip

The equation $\tau_\eta^X= \tau_\eta^C\circ \varphi$ is a
formula in $N_D:=\Gamma(\tilde D,g^*\ol\M^\gp)^*$, which is the
inductive limit of abelian groups $(P_x^\vee)^\gp= (P_x^\gp)^*$ with
respect to the homomorphisms $\iota_{x,\eta}: (P_\eta^\vee)^\gp \to
(P_x^\vee)^\gp$, $x\in\tilde D$. Here $P_x$ for $x\in \tilde D$
means $P_{g(x)}$. More explicitly, if $\Sigma\subset \tilde D$ is
the set of special points $p,q$, that is, mapping to a special point
in $D$,
\begin{equation}\label{Eq: Balancing condition}
N_D=\lim_{\text{\shortstack{$\lra$\\[-0.5ex] $x\in \tilde D$}}}
(P_x^\vee)^\gp = \Big(\bigoplus_{x\in \Sigma}
(P_x^\vee)^\gp\Big)
\Big/\sim,
\end{equation}
where for any $a\in (P_\eta^\vee)^\gp$ and $x,x'\in \Sigma$,
\begin{equation}\label{Eq: Monoid relation}
(0,\ldots,0,\iota_{x,\eta}(a),0,\ldots,0)
\sim (0,\ldots,0,\iota_{x',\eta}(a),0,\ldots,0).
\end{equation}
We may thus represent an element of $N_D$ as a tuple $(a_x)_{x\in
\Sigma}$, but keep in mind the relations~(\ref{Eq: Monoid
relation}).

With this representation of $N_D$ we are now in position to write
down the balancing condition.

\begin{proposition}\label{Prop: Balancing condition}
Consider a Diagram~\ref{Diag: Universal diagram} with $Q=\NN$ and
$D\subset \ul C$ an irreducible component with generic point $\eta$
and $\Sigma\subset\tilde D$ the preimage of the set of special
points. If $\tau_\eta^X$ defined in~(\ref{tau_eta}) is
represented by $(\tau_x)_{x\in\Sigma}$ then
\[
(u_x)_{x\in\Sigma} + (\tau_x)_{x\in \Sigma}=0
\]
in $N_D=\Gamma(\tilde D,\ol \M^\gp)^*$.
\end{proposition}

\proof
Let $m\in \Gamma(\tilde D,g^*\ol\M)$. Recall that $\varphi$
denoted the pull-back of $\ol{f^\flat}$ by $g$. In view of
Lemma~\ref{tau_eta^C} it holds
\begin{eqnarray*}
\tau_\eta^C\big(\varphi(m)\big)
&=&\tau_\eta^C \Big(\big(\langle V_{\eta_q},m\rangle,
\langle V_\eta,m\rangle\big)_{q\in \tilde D},
\big(\langle u_p,m\rangle\big)_{p\in D}\Big)\\
&=&\sum_{q\in \tilde D} \frac{1}{e_q}
\Big(\langle V_{\eta},m\rangle-\langle V_{\eta_q},m\rangle\Big)
-\sum_{p\in \tilde D} \langle u_p,m\rangle.
\end{eqnarray*}
Thus since $u_q(m)= \frac{1}{e_q} \big(\langle V_{\eta_q},m\rangle
-\langle V_\eta,m\rangle\big)$,
\[
\tau_\eta^C\circ \varphi= 
\big((-u_q)_{q\in \tilde D}, (-u_p)_{p\in \tilde D}\big),
\]
and the claimed formula follows from $\tau_\eta^C\circ
\varphi= \tau_\eta^X$, the commutativity of \eqref{taucomm}.
\qed

\subsection{The basicness condition}
\label{Par: basicness}
For a fine saturated sheaf $\ol \M$ over a pre-stable curve
over a field $\ul C/\Spec\kappa$, in Definition~\ref{Def: GS}
we introduced the category $\GS(\ol\M)$. Given a type $\mathbf
u$ for objects of $\GS(\ol\M)$ (Definition~\ref{Def: types in
GS}) we are now in position to construct a universal object for the
full subcategory $\GS(\ol\M, \mathbf u)$.

\begin{construction}\label{Q and basic log structure}
Let $\mathbf u=\big\{(u_p)_p,(u_q)_q \big\}$ be a type for $\GS
(\ol\M)$ and assume $\GS (\ol\M)\neq\emptyset$. For a node $q\in \ul
C$ denote by $\chi_{\eta_i,q}: P_q\to P_{\eta_i}$ the two
generization maps, ordered as in the definition of $u_q$ in
\eqref{eqn at node}. Then if $m\in P_q$ let
\[
a_q(m):= \big( (\ldots,\chi_{\eta_1,q}(m),\ldots,-\chi_{\eta_2,q}(m),\ldots),
(\ldots,u_q(m),\ldots)\big) \in \Big(\prod_\eta
P_\eta\times\prod_q\NN\Big)^\gp
\]
be the element with all entries vanishing except the indicated ones
at places $\eta_1$, $\eta_2$ and $q$. Let $R\subset \big(\prod_\eta
P_\eta\times\prod_q \NN\big)^\gp$ be the saturated subgroup
generated by the $a_q(m)$ for all nodes $q\in\ul C$ and $m\in
P_q$. Now define the \emph{basic monoid} $Q$ as the saturation of
the quotient by $R$:
\begin{equation}\label{Def: Q}
Q:=\Bigg[\iota\Big(\prod_{\eta\in\ul C} P_\eta \times\prod_{q\in \ul
C}\NN\Big)\Big/ R\Bigg]^{\sat}.
\end{equation}
Here $\iota$ denotes the inclusion of $\prod_\eta P_\eta \times
\prod_q\NN$ into its associated group. By the very definition $Q$ is
fine and saturated. Taking the saturation of $R$ amounts to dividing
out any torsion of the associated group, so $Q$ is also torsion-free.
But note that at this point there is no reason to infer
$Q^\times=\{0\}$, and in fact, this is not true in general.

The inclusion of the various factors into $\prod_\eta
P_\eta\times\prod_q\NN$ composed with the surjection to $Q$ defines
homomorphisms
\begin{align*}
&&&& \varphi_{\ol \eta}:P_\eta&\lra
\prod_\eta P_\eta\times\prod_q\NN \lra Q,\\
&&&&\NN&\lra \prod_\eta P_\eta\times\prod_q\NN\lra Q,&
1\longmapsto \rho_q.&&&&
\end{align*}
Since the difference of the two sides of Equation~\eqref{eqn at node}
is nothing but $a_q(m)$, we have the following equality of maps
$P_q\to Q$ (where $\chi_i=\chi_{\eta_i,q}$):
\[
\varphi_{\ol\eta_2}\circ\chi_2 -\varphi_{\ol\eta_1}\circ\chi_1
=  u_q\cdot \rho_q.
\]
From Discussion~\ref{Disc: Stalkwise}, (i)--(iii) the data $Q$,
$\rho_q$ and $\varphi_{\ol\eta}$ thus define a distinguished object
$(Q,\ol\M_C, \psi, \varphi)$ of $\GS (\ol\M,\mathbf u)$,
except that we do not know at this point that $Q^\times=0$ and that
all morphisms are local. These properties will be established in
Proposition~\ref{Prop: Universal property of basicness} provided
$\GS (\ol\M, \mathbf u)\neq\emptyset$.
\end{construction}

\begin{example}
Let us illustrate some features of this definition with some simple
examples, notably concerning the locality property of homoomorphisms
and saturation issues.
\begin{enumerate}
\item
Consider a curve $\ul C$ with two irreducible components
mutually intersecting in two nodes $q_1$, $q_2$.
Assume that the whole curve maps to a standard log point, and hence
$P_x=\NN$ for all $x$ and all generization maps are isomorphisms. Choose
$u_{q_1}(1)=0$, $u_{q_2}(1)=1$. Then in $\prod_\eta P_\eta^\gp\times
\prod_q \ZZ=\ZZ^4$ we have
\[
a_{q_1}(1)=(1,-1,0,0),\quad a_{q_2}(1)=(1,-1,0,1).
\]
Thus $(0,0,0,1)\in R$ and hence $\rho_{q_2}=0$. Thus at $q_2$ the
morphism of monoids $\NN\to \ol\M_{\ul C,q_2}$ mapping $1$ to
$\rho_q$ is not local. In particular, $((\ul
C,\ol\M_C)/(\Spec,\kappa/\NN),\mathbf{x})$ is not a pre-stable log
curve on the level of ghost sheaves (Definition~\ref{Def: pre-stable
at ghost sheaf level}).

For a geometric interpretation of the situation note that $u_q$
compares the lifts of elements of the log structure of $X$ to the
two branches of $\ul C$ at $q$. In the present situation the result
at the two nodes has to agree since $\ul C$ has only two irreducible
components. This can be viewed as a manifestation in log geometry of
the impossibility of a deformation situation where one of the two
nodes smooths while the other stays.
\item
Considering the same situation as in (1) but with $u_{q_1}(1)=2$,
$u_{q_2}(1)=3$, leads to a non-saturated image of $\prod_\eta
P_\eta\times \prod_q\NN$. In fact, $R$ is now generated by
\[
a_{q_1}(1)=(1,-1,2,0),\quad a_{q_2}(1)=(1,-1,0,3),
\]
and the map
\[
\left(\begin{matrix}1&1&0&0\\
-6&0&3&2\end{matrix}\right):
\ZZ^4\lra \ZZ^2
\]
describes the quotient of $\prod_\eta P_\eta^\gp\times \prod_q \ZZ=
\ZZ^4$ by $R$. The image of $\prod_\eta P_\eta\times \prod_q
\NN=\NN^4$ is then generated by the columns of the matrix, and hence
is not saturated for it contains $(0,2)$ and $(0,3)$ but not
$(0,1)$. Thus we need to saturate in the definition of $Q$ to stay
in the category of fine saturated log schemes.\footnote{We realized
this saturation issue only after studying \cite{chen} in greater
detail.}
\item
Consider again a curve as in (1), but now mapping to a pre-stable
log curve with one node in such a way that
\[
P_{q_1}=\NN^2,\quad P_{\eta_1}=P_{q_2}=P_{\eta_2}=\NN,
\]
with only non-trivial generization maps $P_{q_1}\to P_{\eta_i}$ the
two projections $\NN^2\to \NN$. Take $u_{q_1}(a,b)= a+b$,
$u_{q_2}(c)= 2c$. Then
\[
a_{q_1}(1,0)= (1,0,1,0),\quad
a_{q_1}(0,1)=(0,-1,1,0),\quad
a_{q_2}(1)= (1,-1,0,2).
\]
The subgroup of $\prod_\eta P_\eta^\gp\times \prod_q\ZZ= \ZZ^4$
generated by these elements contains $(0,0,2,-2)= a_{q_1}(1,1)-
a_{q_2}(1)$, but it does not contain $(0,0,1,-1)$. Hence saturation
in the definition of $R$ is necessary to make $(\prod_\eta
P_\eta^\gp\times \prod_q\ZZ)/R$ torsion-free.
\end{enumerate}
\end{example}

\begin{remark}\label{Rem: Q^vee}
Another useful way of thinking about $Q$ is in the dual space:
\[
Q^\vee= \Big\{ \big((V_\eta)_\eta,(e_q)_q\big)\in
{\textstyle \bigoplus_\eta
P_\eta^\vee \oplus\bigoplus_q \NN} \,\Big|\,
\forall q: V_{\eta_2}-V_{\eta_1}=e_q u_q\Big\}.
\]
In particular, this avoids any saturation issues. Note that
following the discussion in \S\ref{par: Log maps over std log pt}
any $((V_\eta),(e_q))\in Q^\vee$ that does not lie in any proper
face gives rise to an object in the category $\GS(\ol\M)$ of the
form $(\NN,\ol\M_C, \psi,\phi)$. In fact, this is how the authors
first found $Q$.
\end{remark}

\begin{proposition}
\label{Prop: Universal property of basicness}
If $\GS(\ol\M,\mathbf u)\neq \emptyset$ then it has as initial
object the tuple $(Q,\ol\M_C, \psi, \varphi)$ from
Construction~\ref{Q and basic log structure}.
\end{proposition}

\proof
We first prove the universal property in an enlarged category
without the assumption that $Q^\times=\{0\}$ and the morphisms of
monoids local. Let $(Q',\ol\M'_C, \psi',\varphi')\in \Ob\big(
\GS(\ol\M,\mathbf u)\big)$. Let $\rho'_q\in Q'$ be the element
defining $\ol\M'_{C,\ol q}$. Then because of \eqref{eqn at
node} for $(Q',\ol\M'_C, \psi',\varphi')$ and since $Q'$ is
torsion-free and saturated the map
\begin{equation}\label{universal stanard lifting}
\xymatrix@C=60pt{
\prod_\eta P_\eta\times
\prod_q\NN \hspace{-5pt}
\quad\ar[r]^{\mbox{}\hspace{40pt}
\prod_\eta \varphi'_{\ol\eta}\times \prod_q\rho'_q}
&Q'}
\end{equation}
factors over $Q$. Tracing the image of the generator of the $q$-th
copy of $\NN$ shows that the induced map $Q\to Q'$ maps $\rho_q$ to
$\rho'_q$. Moreover, by the very definition of this factorization it
is compatible with $\varphi_{\ol \eta}:P_\eta \to Q$ and
$\varphi'_{\ol\eta}:P_\eta\to Q'$. This proves existence of a
morphism $(Q,\ol\M_C, \psi, \varphi)\to (Q',\ol\M'_C,
\psi',\varphi')$.

For uniqueness note that by compatibility of $\rho_q$, $\rho'_q$ and of
$\varphi_{\ol \eta}$, $\varphi'_{\ol \eta}$ any such morphism would
lift to the homomorphism stated in \eqref{universal stanard
lifting}.

It remains to show that $(Q,\ol\M_C, \psi, \varphi)\in
\GS(\ol\M, \mathbf u)$, that is, that $Q^\times=0$
and all sheaf homomorphisms are local. Since $\GS(\ol \M,\mathbf
u)\neq \emptyset$ there is at least one morphism $(Q,\ol\M_C,
\psi, \varphi)\to (Q',\ol\M'_C, \psi',\varphi')$ as
constructed above. Now because $(Q',\ol\M'_C,
\psi',\varphi')\in \Ob(\GS(\ol \M, \mathbf u))$, for any
$\eta$ and $q$ the compositions
\[
P_\eta\stackrel{\varphi_{\ol\eta}}{\lra}Q\lra Q',\quad
\NN\stackrel{\cdot\rho_q}{\lra} Q\lra Q'
\]
are local. This proves that $\psi$ and $\varphi$ are
indeed local homomorphisms. Thus also the composition
\[
\prod_\eta P_\eta\times\prod_q\NN\lra Q\lra Q'
\]
is local and hence, by surjectivity of the first arrow up to
saturation, also $Q\to Q'$ is local and $Q^\times=\{0\}$. Finally
$\ol\M_C \to\ol\M'_C$ is local because at a node the homomorphism is
defined by the product of the homomorphisms at the adjacent generic
points.
\qed

\begin{definition}\label{Def: Basicness}
A stable log map $(C/W, \mathbf{x}, f)$ is called \emph{basic} if for any
geometric point $\ol w\to \ul W$ the induced object of $\GS(\ul
f_{\ol w}^*\ol \M_X, \mathbf u)$ is universal. Here $\mathbf
u$ is the type of $(C/W, \mathbf{x}, f)$ at $\ol w$.
\end{definition}

\begin{remark}
Following up on Remark~\ref{Rem: Q^vee} the tropical interpretation
of the basicness condition is as follows. If $\ol w:\Spec \kappa\to
\ul W$ is a geometric point then $\Int {\ol\M}_{W,\ol w}^\vee$ can
be identified with the set of tropical data of the pull-backs to
$(\Spec \kappa,\NN)$ via a log enhancement of $\ol w$. Moreover, in the
proper tropical situations of Discussion~\ref{Disc: Tropical curve}
the real cone generated by $\ol\M_{W,\ol w}^\vee$ is canonically
isomorphic to
\[
\Hom(Q,\RR_{\ge 0})= \Big\{ \big((V_\eta)_\eta,(l_q)_q\big)\in
{\textstyle \bigoplus_\eta
\Hom(P_\eta,\RR_{\ge 0}) \oplus\bigoplus_q \RR_{\ge0}} \,\Big|\,
\forall q: V_{\eta_2}-V_{\eta_1}=l_q u_q\Big\}.
\]
Each point in this space defines a tropical curve with vertices
$V_\eta$ and interior edges $q$ mapping to a translation of the line
segment connecting $0$ with $l_q u_q$. This can be interpreted as
the moduli space of tropical curves of the given type. The proper
faces of $\Hom(Q,\RR_{\ge0})$ parametrize degenerate curves for
which there may be some edge of length~$0$ or some vertex
$V_{\eta}$ which maps to the boundary of $\Hom(P_{\eta},\RR_{\ge
0})$.
\end{remark}

As a first property we show that basicness is an open condition:

\begin{proposition}\label{Prop: Basicness is open}
Let $(C/W, \mathbf x, f)$ be a stable log map to a log
scheme $X$. Then
\[
\Omega:=\big\{\ol w\in |\ul W|\,\big|\, {\{\ol w\}} \times_{\ul W}
(C/W, \mathbf x,f)
\text{ is basic}\big\}
\]
is an open subset of $|\ul W|$.\footnote{In writing ${\{\ol w\}}
\times_{\ul W} (C/W, \mathbf x,f)$ we view $(C/W, \mathbf
x,f)$ as an object over $\ul W$. In particular, this fibre
product is a log curve over $\Spec\kappa(\ol w)$ endowed with the log
structure making the inclusion into $W$ strict.}
\end{proposition}

\proof
Since basicness is a condition on morphisms of fine sheaves,
$\Omega$ is constructible.\footnote{Alternatively, the following
arguments indeed show that basicness holds on subsets of $\ul W$
which admit a chart $Q\to \M_W$ inducing an isomorphism $Q\simeq
\ol\M_{W,\ol w}$ for some $\ol w\in |\Omega|$ (cf.\ \cite{chen},
Proposition~3.5.2).} It remains to show that $\Omega$ is closed
under generization. So let $\ol w_1\in \Omega$, $\ol w_2\in |\ul W|$
and $\ol w_1\in\cl(\ol w_2)$. We need to show $\ol w_2\in\Omega$.
Since basicness is stable under strict base change we may first
replace $\ul W$ by $\Spec\big( \O_{\ul W,{\ol w_1}}\big)$ and then
by $\cl(\ol w_2)$ with the induced reduced scheme structure, to
reduce to the case $\ul W=\Spec R$ for a strictly Henselian local
domain $R$, and with $\ol w_1$ and $\ol w_2$ the closed point $0$
and generic point $\Spec K$, $K$ the quotient field of $R$. Denote
by $\kappa=R/\maxid$ the residue field of $R$, and endow
$\Spec\kappa$ and $\Spec K$ with the log structures induced by the
embeddings into $W= (\Spec R,\M_R)$.

Now we have two relevant stable log maps over fields, the closed fibre
\[
(C_0/(\Spec \kappa,Q), \mathbf x_0,f_0) := 
\Spec\kappa \times_{\ul W}(C/W, \mathbf x,f),
\]
which is basic by assumption, and the generic fibre
\[
(C_K/(\Spec K,Q_K), \mathbf x_K,f_K) := 
\Spec K \times_{\ul W}(C/W, \mathbf x,f),
\]
a stable log map over some log point $(\Spec K,Q_K)$.
We use the standard notations for the points and the monoids
of the closed fibre, while for the generic fibre we add hats. Note
that apart from the usual generization maps between points on the
same fibre we also have generization homomorphisms from the closed
to the generic fibre,
\[
P_\eta\lra P_{\hat\eta},\quad
P_p\lra P_{\hat p},\quad
P_q\lra P_{\hat q},
\]
where $P_q\lra P_{\hat q}$ only exists for those nodes $q\in
\ul C_0$ that are contained in the closure of a node of
$\ul C_K$.

For each generic point $\eta\in \ul C_0$ with generization
$\hat\eta\in \ul C_K$ the homomorphism $\ul
f^*\ol\M_X\to \ol \M_C$ defines a commutative square
\[
\begin{CD}
P_\eta @>>> Q\\
@VVV@VVV\\
P_{\hat\eta}@>>> Q_K.
\end{CD}
\]
Moreover, if $q\in\cl(\hat q)$ then $\rho_q$ maps to $\rho_{\hat
q}$. We thus obtain a commutative diagram
\begin{equation}\label{Fig: Generization and basicness}
\begin{CD}
P:=\prod_\eta P_\eta \times\prod_q\NN @>{\theta}>> Q\\
@V{\tilde\psi}VV @VV{\psi}V\\
\hat P:=\prod_{\hat \eta} P_{\hat \eta} \times\prod_{\hat q}\NN
@>{\hat\theta}>> Q_K,
\end{CD}
\end{equation}
with $\theta$ the map defining $Q$ up to saturation and the
vertical arrows the generization epimorphisms. If $q$ is a node
not in the closure of some $\hat q$ then $\tilde\psi$ maps
this copy of $\NN$ to $0$.

As a generization homomorphism, $\psi$ induces an isomorphism
$S^{-1}Q/(S^{-1}Q)^\times\to Q_K$ with $S=\psi^{-1}(0)$ a face of
$Q$. Similarly, $\tilde\psi$ induces an isomorphism $\tilde S^{-1}
P/(\tilde S^{-1} P)^\times\to \hat P$. As the morphisms $\theta$ and
$\hat\theta$ are local we have the relation
\begin{equation}\label{Eq: S and tilde S}
\tilde S={\tilde\psi}^{-1}(0)=
{\tilde\psi}^{-1} \big({\hat\theta}^{-1}(0)\big)=
\theta^{-1}\big(\psi^{-1}(0)\big)=\theta^{-1}(S).
\end{equation}
Recall from Construction~\ref{Q and basic log structure} that
$Q$ is the saturation of the image of $P$ in $P^\gp/R$, where $R$ is
the saturation of the subgroup of $P^\gp$ generated by certain
elements $a_q(m)$, for all nodes $q$ and all $m\in P_q$. By the
compatibility of types with generization (Lemma~\ref{Lem: type under
generization}) $\tilde\psi^\gp$ maps $R$ to the analogous subgroup
$\hat R\subset \hat P$. Since by \eqref{Fig: Generization and
basicness} $\hat \theta$ is surjective up to saturation, it remains
to show that up to saturation $\hat\theta$ is the quotient of $\hat
P$ by $\hat R$. Thus let $\hat m_1,\hat m_2 \in \hat P$ with
$\hat\theta (\hat m_1)= \hat\theta(\hat m_2)$. Let $m_1,m_2$ be
lifts of $\hat m_1,\hat m_2$ to $P$. Then
$\theta(m_1)-\theta(m_2)\in S^\gp$ (viewed in $Q^\gp$) for
\[
\psi(\theta(m_1))=\hat\theta(\hat m_1)=\hat\theta (\hat m_2)
=\psi(\theta(m_2)).
\]
Thus there exist $h_i\in \theta^{-1}(S)$ such that $\theta(m_1+h_1)
= \theta(m_2+h_2)$. Using $\theta^{-1}(S)= \tilde S$ from \eqref{Eq:
S and tilde S} we may replace $m_i$ by $m_i+ h_i$ to achieve
$\theta(m_1)=\theta(m_2)$. But then $m_1- m_2\in R$, and hence
\[
\tilde\psi(m_1)-\tilde\psi(m_2)\in \tilde\psi^\gp(R)=\hat R,
\]
finishing the proof.
\qed
\medskip

The next proposition establishes a universal property for basic
stable log maps. This shows in particular that we do not lose any
generality in imposing basicness. This result is not needed for the
construction of log Gromov-Witten invariants, but is included for
reassurance. It is also referred to in the comparison with Jun Li's
moduli space in Corollary~\ref{Cor: Comparison with Jun Li's}. As an
auxiliary result we first treat the problem on the level of ghost
sheaves.

\begin{lemma}
\label{Lem: GS(M/W)}
Let $(\ul C/\ul W,\mathbf x)$ be a pre-stable curve and $\ol\M$ a
fine saturated sheaf on $\ul C$. For each geometric point $\ol
w:\Spec\kappa(\ol w)\to \ul W$ let be given a type $\mathbf{u}_{\ol
w}$ of an object of $\GS(\ol\M|_{\ul C_{\ol w}})$ such that the
collection $(\mathbf u_{\ol w})_{\ol w}$ of types is compatible with
generization. Then the full subcategory $\GS\big(\ol\M,(\mathbf
u_{\ol w})\big)$ of objects of $\GS(\ol\M)$ that have type $\mathbf
u_{\ol w}$ over the geometric point $\ol w$ has a universal object. 
\end{lemma}

\proof
For geometric points this is the statement of Proposition~\ref{Prop:
Universal property of basicness}. 
In the general case Diagram~\eqref{Fig: Generization and
basicness} in the proof of Proposition~\ref{Prop: Basicness is open}
shows that the fibrewise defined diagrams of ghost sheaves are
compatible with generization. Hence they define the desired initial
object in the category of diagrams of ghost sheaves.
\qed
\medskip

\begin{proposition}
\label{Prop: Unique pull-back from basic stable log map}
Any stable log map arises as the pull-back from a basic stable log map
with the same underlying ordinary stable map. Both the basic stable log
map and the morphism are unique up to unique isomorphism.
\end{proposition}

\proof
Let $(\pi:C\to W,\mathbf x,f)$ be the given stable log map, defining
morphisms of log structures $\pi^\flat:\ul\pi^*\M_W\to\M_C$ and
$f^\flat:\ul f^*\M_X\to\M_C$. We consider the category of morphisms
of log structures $\M'_W\to \M_W$ on $\ul W$ and $\M'_C\to \M_C$,
$\ul\pi^*\M'_W\to \M'_C$, $\ul f^*\M_X\to\M'_C$ on $\ul C$,
compatible with $\pi^\flat$ and $f^\flat$ in the obvious way. The
statement follows once we show that this category has an initial
object, and that this object is a basic stable log map.

If $\M_1\to \M_2$ is a morphism of fine log structures on a scheme $Y$
then from the commutative diagram
\[
\begin{CD}
1@>>> \O_Y^\times @>>> \M_1^\gp @>>> \ol\M_1^\gp @>>>0\\
@.@|@VVV@VVV\\
1@>>> \O_Y^\times @>>> \M_2^\gp @>>> \ol\M_2^\gp @>>>0
\end{CD}
\]
it follows that $\M_1= \M_2\times_{\ol \M_2} \ol\M_1$. Moreover,
a morphism of log structures with target $\shM_2$ lifts to $\shM_1$
if and only if this is true on the level of ghost sheaves. In
particular, the functor defined by going over from a stable log map
to the associated diagram of ghost sheaves is an equivalence from
the present category to $\GS\big(\ol\M, (\mathbf u_{\ol w})\big)$.
Here $\mathbf u_{\ol w}$ is the type of the given stable log map
$(C/W,\mathbf x,f)$ at the geometric point $\ol w$.
The statement now follows from Lemma~\ref{Lem: GS(M/W)}.
\qed
\medskip

Another remarkable property of basic stable log maps is that they do
not admit non-trivial automorphisms that are the identity on the
underlying ordinary stable maps. Stack-theoret\-ically this means
that the forgetful map from the stack of basic stable log maps to
the stack of ordinary stable maps is representable
(Proposition~\ref{Prop: forgetful morphism to M(X)}). The statement
is also useful for checking that the stack of stable log maps has a
separated diagonal (Proposition~\ref{Prop: separated diagonal}).

\begin{proposition}
\label{Prop: no non-trivial log automorphisms}
An automorphism $\varphi: C/W\to C/W$ of a basic
stable log map $(\pi:C\to W, \mathbf{x},f)$ with
$\ul\varphi=\id_{\ul C}$ is trivial.
\end{proposition}

\proof
An automorphism of $C/W$ is an automorphism $(\ul\varphi,
\varphi^\flat)$ of $C=(\ul C,\M_C)$ descending to an automorphism
$\psi=(\id, \psi^\flat)$ of $(\ul W,\M_W)$. It is an
automorphism of $(C/W, \mathbf{x},f)$ if $\ul{\varphi}
(\mathbf{x}) =\mathbf{x}$ and if it commutes with the morphism
$f^\flat: \ul f^*\M_X\to \M_C$. Since $\ul\varphi=\id_C$ the
latter condition means
\begin{equation}\label{varphi^flat eqn}
\varphi^\flat\circ f^\flat=f^\flat.
\end{equation}

We claim that it suffices to show $\psi^\flat=\id_{\M_W}$. In fact,
let $U\subset \ul C$ be the complement of the set of special points
(images of marked points or criticial points of $\ul \pi$). Because 
$\pi^\flat: \pi^*\M_W\to \M_C$ is an isomorphism on $U$,
$\psi^\flat= \id_{\M_W}$ implies $\varphi^\flat|_U = \id$. But by
Theorem~\ref{Thm: structure of log curves}, $\M_C$ has no section
with support on the set of special points. Hence
$\varphi^\flat=\id_{\M_C}$.

Now on the level of ghost sheaves a basic stable log map is
determined by the underlying ordinary stable map and the type.
Moreover, there is no morphism between stable log maps of different
types. This shows that $\ul\varphi =\id_{\ul C}$ implies
$\ol{\varphi^\flat} =\id_{\ol\M_C}$. By strictness of $\pi: C\to W$
away from the special points this imples $\ol{\psi^\flat}=
\id_{\ol\M_W}$. Thus $\psi^\flat$ has the form
\[
\psi^\flat(m)= h(\ol m)\cdot m
\]
for a  homomorphism $h:\ol\M_W\to \O_{\ul W}^\times$. Now if $\ol
w\to \ul W$ is a geometric point, then by basicness, up to
saturation $\ol\M_{W,\ol w}$ is generated by the fibres $P_\eta$ of
$\ul f^*\ol\M_X$ at the generic points $\eta\in \ul C_{\ol w}$, and
by one copy of $\NN$ for each singular point of $\ul C_{\ol w}$.
Strictness of $\pi|_U$ together with \eqref{varphi^flat eqn}
implies that $h$ is trivial on the part of $\ol\M_W$ generated by
$P_\eta$. The factors of $\NN$ are generated by the image of the
basic log structure $\pi^*\M_W^o\to \M_C^o$ on the pre-stable curve
$\ul C/\ul W$. But $\M_C^o$ has no non-trivial automorphism inducing
the identity on $\ul C$ and on $\ol{\M^o_C\!\!}\,\,$. This is a
direct consequence of the existence of basic log structures for
pre-stable curves, see Appendix~\ref{App: Prestable curves}. Hence
$h=1$ also on this part. Finally note that a homomorphism from a
fine monoid to a group is trivial if and only if it is trivial on
its saturation. Hence $h$ is trivial and thus
$\psi^\flat=\id_{\M_W}$ as remained to be shown.
\qed

\section{Algebraicity}

\subsection{The stack of stable log maps}
We are now ready to define the stack of stable log maps. We continue
with the convention that $X$ is a log scheme over $S$, with the log
structures defined on the Zariski sites. Remember also that all our
log structures are defined over $S$, that is, come with a morphism
from the pull-back of $\M_S$ that is compatible with $\M_S\to\O_{\ul
S}$. We endow $\SchS$ with the \'etale Grothendieck topology.

\begin{definition}\label{Def: Stack of basic stable log maps}
The \emph{stack of stable log maps to $X$ (over $S$)} is the
category $\tilde\MM(X)= \tilde\MM(X/S)$ (Definition~\ref{Def: Stable
log map}) together with the forgetful morphism $\tilde\MM(X) \to
(\mathrm{Sch}/\ul S)$ mapping $(C/W,\mathbf x,f)$ to $\ul W$. The full
subcategory of \emph{basic} stable log maps is denoted $\MM(X)$.
\end{definition}

Since the morphisms in $\tilde\MM(X)$ and $\MM(X)$ are given by
cartesian diagrams of log smooth curves over the underlying base
schemes, $\tilde \MM(X)\to (\mathrm{Sch}/\ul S)$ and
$\MM(X)\to (\mathrm{Sch}/\ul S)$ are fibred groupoids. As is
customary, for any $a\in\tilde \MM(X)$ over $\ul W\in\SchS$
and $\varphi: \ul V\to \ul W$ we choose one morphism
in $\tilde \MM(X)$ covering $\varphi$ and denote it by
$\varphi^*a\to a$.

Once we prove that $\tilde\MM(X)$ is an algebraic stack,
Proposition~\ref{Prop: Basicness is open} shows that $\MM(X)$ is
also algebraic, for it is an an open substack of $\tilde \MM(X)$.
We therefore restrict attention to $\tilde \MM(X)$ for most of this
section.

\begin{lemma}
$\tilde\MM(X)$ is a stack.
\end{lemma}
\proof
We verify Axioms~(i) and (ii) in \cite{champs}, Definition~3.1.

(i)~We have to check the sheaf axioms for morphisms between two
objects in $\tilde\MM(X)$ over the same base scheme. This amounts to
the following. Let $a_i\in \tilde\MM(X)_{\ul W}$, $i=1,2$, be
two stable log maps with the same base scheme $\ul W$. Let
$h:\tilde{\ul W}\to \ul W$ be an \'etale cover,
$\tilde\psi: h^*a_1 \to h^*a_2$ a morphism over $\tilde {\ul
W}$ and $\pr_i:\tilde{\ul W}\times_{\ul W} \tilde
{\ul W}\to \tilde {\ul W}$ for $i=1,2$ the projections.
With $\pi:=h\circ\pr_1= h\circ\pr_2$ there are two morphisms
\[
\pr_i^*\tilde\psi:\pi^*a_1\lra \pi^* a_2,\quad i=1,2.
\]
The sheaf axiom says that if these two morphisms agree then there
exists a unique morphism $\psi: a_1\to a_2$ with
$\tilde\psi=h^*\psi$. For the underlying morphisms of schemes this
follows from faithfully flat descent (\cite{SGA1}, VIII,
Theorem~5.2). Since on the domains we work with log structures in
the \'etale topology the refinement to morphisms of log schemes is a
tautology.

(ii)~This axiom deals with descent for objects in $\tilde\MM(X)$. On
the underlying domain of the stable log map this follows by the
sheaf property of the stack of pre-stable curves $\mathbf M$ (see
Appendix~\ref{App: Prestable curves}). Then the underlying morphisms
to $\ul X$ descend as in~(i). Again the refinement to morphisms of
log spaces is a tautology.
\qed
\medskip

The rest of this section is devoted to proving algebraicity of
$\tilde\MM(X)$. Denote by $\scrM= \scrM_S$ the log stack over $S$ of
(not necessarily basic) pre-stable marked log curves. In
Appendix~\ref{App: Prestable curves} we recall the folklore result
that $\scrM$ is an algebraic stack locally of finite type over $\ul
S$ (Proposition~\ref{Prop: MM is algebraic}). There is a forgetful
morphism of stacks
\[
\tilde\MM(X)\lra \scrM,
\]
mapping a stable log map $(C/W, \mathbf{x}, f)$ to the pre-stable
marked log curve $(C/W,\mathbf{x})$. Note that this functor is
faithful, since a morphism of stable log maps is given by a morphism
on the domains. Another forgetful morphism is to the algebraic stack of
ordinary stable maps $\mathbf{M}(\ul X)$ \cite{behrendmanin}:
\[
\tilde\MM(X)\lra \mathbf{M}(\ul X).
\]
In the next subsection we will prove the following.

\begin{proposition}
\label{representability over M x M(X)}
$\tilde\MM(X)\to \scrM\times \mathbf{M}(\ul X)$ is representable and
locally of finite presentation.
\end{proposition}

A direct consequence is the main result of this section.

\begin{theorem}
\label{Thm: algebraicity of stack}
$\tilde\MM(X)$ is an algebraic stack locally of finite type over
$\ul S$.
\end{theorem}

\proof
By \cite{champs}, Proposition~4.5,(ii) algebraicity of $\tilde
\MM(X)$ follows from Proposition~\ref{representability over M x
M(X)}. It is locally of finite type over $\ul S$ since so are $\MM$
and $\mathbf M(\ul X)$.
\qed
\medskip

As is the case with $\scrM$, the stack $\tilde\MM(X)$ is only an
algebraic stack in the sense of \cite{OlssonENS}, p.750. This
definition drops the separatedness of the diagonal morphism
(quasi-separatedness) from the definition in \cite{champs}. In
contrast, the open substack $\MM(X)$ of $\tilde\MM(X)$ does have a
separated diagonal, so is an algebraic stack in the sense of
\cite{champs}:

\begin{proposition}\label{Prop: separated diagonal}
The diagonal morphism $\Delta_{\MM(X)/\ul S}: \MM(X)\lra
\MM(X)\times \MM(X)$ is separated.
\end{proposition}

\proof
The statement amounts to the following: An automorphism of a stable
log map $(C/W,\mathbf x,f)$ over an integral scheme $\ul W$ that is
generically the identity is trivial. This is clearly true for the
underlying ordinary stable map, and lifts to basic stable log maps
by virtue of Proposition~\ref{Prop: no non-trivial log
automorphisms}.
\qed
\medskip

As a corollary of Theorem~\ref{Thm: algebraicity of stack} we obtain
algebraicity of the open substack $\MM(X)$ of $\tilde\MM(X)$ of
basic stable log maps. Moreover, $\MM(X)$ carries a canonical log
structure, that is, a factorization
\[
\MM(X)\lra (\mathrm{Log}/S)\lra \SchS,
\]
where $(\mathrm{Log}/S)$ is the category of fine saturated log
schemes over $S$ with strict morphisms.

\begin{corollary}\label{MM(X) is algebraic log stack}
$\MM(X)$ is an algebraic log stack with separated diagonal and
locally of finite type over $\ul S$.
\end{corollary}

\proof
For the log structure, define the functor $\MM(X)\lra
(\mathrm{Log}/S)$ by sending a stable log map $(C/W,\mathbf x,f)$ to
the logarithmic base $W$.
\qed
\medskip

Assuming Proposition~\ref{representability over M x M(X)} we can
establish at this point some further properties of $\MM(X)$.

\begin{proposition}
\label{Prop: forgetful morphism to M(X)}
The forgetful morphism $\MM(X)\to \mathbf{M}(\ul X)$ of algebraic
stacks is representable.
\end{proposition}

\proof
By Corollaire~8.1.2 in \cite{champs} we have to show that the
diagonal morphism
\[
\MM(X)\lra \MM(X)\times_{\mathbf{M}(\ul X)} \MM(X)
\]
is a monomorphism. This amounts to the statement about automorphisms
of basic stable log maps verified in Proposition~\ref{Prop: no non-trivial
log automorphisms}.
\qed

\begin{corollary}
\label{Cor: C(X^ls) is DM}
The algebraic stack $\MM(X)$ is a Deligne-Mumford stack.
\end{corollary} 

\proof
An algebraic stack representable over a Deligne-Mumford stack is
itself a Deligne-Mumford stack. In fact, if $X\to \mathfrak X$ is an
\'etale presentation of an algebraic stack and $\mathfrak Y\to
\mathfrak X$ is representable, then $Y:= X\times_{\mathfrak X}
\mathfrak Y\to \mathfrak Y$ is also an \'etale surjection, and $Y$
is an algebraic space by representability.
\qed

\subsection{Representability of spaces of log morphisms}
We now prove Proposition~\ref{representability over M x M(X)}. To
avoid excessive underlining, in this subsection we change our
convention and denote schemes or algebraic spaces by unadorned
letters and use a dagger for log schemes, as in $X^\ls=(X,\M_X)$. We
have to show that if $W$ is a scheme\footnote{We follow the usual
convention to confuse a scheme and its associated stack} and $W\to
\scrM\times \mathbf{M}( X)$ is a morphism then the fibre product
$W\times_{\scrM\times \mathbf{M}(X)} \tilde\MM(X)$ is an algebraic
space locally of finite type over $W$. Explicitly, this means the
following. The morphism $W\to \scrM\times \mathbf{M}(X)$ amounts
to giving a pair $(\shW^\ls,f)$ with $\shW^\ls=(C^\ls/W^\ls,\mathbf{x})$ a
pre-stable marked log curve and $f: C\to X$ a morphism of schemes
making $(C/W,\mathbf{x},f)$ an ordinary stable map. For $V\to W$, the
fibre category of $W\times_{\scrM\times \mathbf{M}(X)} \tilde\MM(X)$
over $V$ can be taken as the category with objects morphisms of log
structures
\[
(f^*\M_X)_V\lra (\M_C)_V
\]
over $S^\ls$ and with only the identity as morphisms. Here the index
$V$ means pull-back via the base change morphism $C_V:=V\times_{W}
C\to C$. Triviality of the automorphisms in this fibre category is
due to the fact that $\tilde\MM(X)\to \scrM\times \mathbf{M}(X)$ is
faithful. Thus the question is about the representability of the
functor of morphisms between two given log structures over $S^\ls$
along the fibres of the proper morphism $C\to W$.

Abstracting, now let $\pi: Y\to W$ be a proper morphism of schemes
and let $\alpha_i:\M_i\to\O_Y$, $i=1,2$, be two fine saturated log
structures on $Y$. Consider the functor
\begin{equation}
\label{LMor_Y/W}
\operatorname{LMor}_{Y/W}(\M_1,\M_2): (\mathrm{Sch}/
W)\lra (\operatorname{Sets})
\end{equation}
that on objects is defined by
\[
V\longmapsto \big\{ \varphi:(\M_1)_V\to (\M_2)_V\ \text{morphism of
log structures}\big\}.
\]
Then the statement that $W\times_{\scrM\times
\mathbf{M}(X)} \tilde\MM(X)$ is an algebraic space\footnote{As for
algebraic stacks we have to drop the condition of
quasi-separatedness from the definition of algebraic spaces
(\cite{knutson}, Ch.II, Definition~1.1).}
 essentially is a special case of the following proposition.
 
\begin{proposition}
\label{LMor is representable}
$\operatorname{LMor}_{Y/W} (\M_1,\M_2)$ is represented by an
algebraic space locally of finite type over $W$.
\end{proposition}

The proof is provided, after some preparations, at the end of this
subsection. This then finishes the proof of
Proposition~\ref{representability over M x M(X)} in the case
$S=\Spec \kk$ with the trivial log structure. In the general case we
have in addition two morphisms of log structures $\psi_i:\pi^*
\M_S\to\M_i$, $ \pi: C\to W$ the projection, and need to restrict to
those $\varphi: (\M_1)_V\to (\M_2)_V$ compatible with $\psi_i$. But
by Proposition~\ref{LMor is representable} composition with $\psi_1$
defines a morphism of algebraic spaces
\[
\operatorname{LMor}_{Y/W}(\M_1,\M_2) \lra
\operatorname{LMor}_{Y/W}(\pi^*\M_S,\M_2).
\]
Now $W\times_{\scrM\times \mathbf{M}(X)}
\tilde\MM(X)$ arises as the fibre product with the morphism
\[
W\lra \operatorname{LMor}_{Y/
W}(\pi^*\M_S,\M_2)
\]
defined by $\psi_2$, and is hence represented by an algebraic
space. This finishes the proof of Proposition~\ref{representability
over M x M(X)} also in the general case.

\begin{remark}
\label{Rem: non-sparatedness issue}
One problem in showing representability of
$\operatorname{LMor}_{Y/W}(\M_1, \M_2)$ is that it is
non-separated, essentially because the induced map $\ol\M_1\to
\ol\M_2$ cannot be determined by its restriction to an open dense
subset. As an example (cf.\ \cite{OlssonENS}, Remark~3.12) consider
the log structure $\M$ on $\AA^1=\Spec \kk[x]$ with chart
\[
\NN^2\lra \kk[x],\quad
(a,b)\longmapsto x^{a+b}.
\]
The map $\NN^2\to\NN^2$, $(a,b)\mapsto (b,a)$ induces a non-trivial
automorphism of $\M$ that restricts to the identity on
$\AA^1\setminus\{0\}$.
\end{remark}

To find an \'etale cover of the algebraic space representing
$\operatorname{LMor}_{Y/W} (\M_1,\M_2)$ in
Proposition~\ref{LMor is representable} we
thus first restrict the map $\ol\varphi$. To this end let
$\ol w\to W$ be a geometric point and
\[
\ol\varphi_{\ol w}: (\ol\M_1)_{\ol w}
\to (\ol\M_2)_{\ol w}
\]
be a choice of $\ol\varphi$ over one geometric fibre
$Y_{\ol w}$ of $Y\to W$.
Now since $\ol\M_i$ are fine sheaves the choice at a geometric
point $\ol x\to Y$ determines $\ol\varphi$ at any
generization $\ol y$ of $\ol x$. Moreover, if $\ol
y$ specializes to some other point $\ol z$ such that the
generization map $\ol \M_{i,\ol z} \to \ol
\M_{i,\ol y}$ is an isomorphism then $\ol\varphi$ is
also determined at $\ol z$. Iterating the
generization-specialization process we are lead to the following
definition.

\begin{definition}\label{Def: A_gen}
Let $\ol\M$ be a fine sheaf on a scheme $Y$ and let $A\subset |Y|$
be a set of geometric points. We say $\ol x\in |Y|$ has \emph{property
$(A_{\mathrm{gen}})$ with respect to $\ol\M$} if there exists a
sequence $\ol y_1,\ol z_1,\ldots,\ol y_r, \ol z_r\in|Y|$ for some
$r$ with the following properties.
\[
\begin{array}{ll}
\ol y_1\in A,\quad \ol z_r=\ol x,\\
\ol y_i\in\cl(\ol z_i),&i=1,\ldots,r,\\
\ol y_i\in\cl(\ol z_{i-1})&\text{and}\quad
\ol \M_{\ol y_i}\to \ol\M_{\ol z_{i-1}}
\text{ is an isomorphism},\ i=2,\ldots,r.
\end{array}
\]
\end{definition}
\medskip

Thus the giving of $\ol\varphi$ on a closed subset $A\subset
|Y|$ then determines $\ol\varphi$ also on the subset
\[
U_A:=\big\{ \ol x\in|Y|\,\big|\, \ol x\text{ fulfills } (A_{\mathrm
gen})\big\}
\]
of $|Y|$. Note that by definition $U_A$ is closed under
generization. Since $\ol\M$ is a fine sheaf it is also immediate
that $U_A$ is a constructible subset of $|Y|$, and hence $U_A\subset
|Y|$ is open.

Since the statement of Proposition~\ref{LMor is representable} is
local in $W$ and by properness of $Y\to W$ we
may assume any point of $|Y|$ fulfills $(A_{\mathrm{gen}})$ for
$A=Y_{\ol w}$ with respect both to $\ol\M_1$ and to
$\ol\M_2$, that is, $U_{Y_{\ol w}}= Y$. Then for any
$V\to W$ there is at most one
$\ol\varphi: (\ol\M_1)_V\to (\ol\M_2)_V$
compatible with $\ol\varphi_{\ol w}$ under sequences of
generization maps. Let us call such $\ol\varphi$ (or a lift
$\varphi$ to a morphism of log structures) \emph{compatible with
$\ol\varphi_{\ol w}$}, and similarly for any $A\subset
|Y|$. Note that $\ol\varphi_{\ol w}$ may not extend to
$Y$, but it may do so after certain base changes.

We first treat the representability problem locally on $Y$, that is,
for $Y=W$.

\begin{lemma} 
\label{logmapmodulilemma}
Let $Y=W$ and suppose that there exists a
closed subset $A\subset Y$ such that any $\ol x\in |Y|$ fulfills
$(A_\mathrm{gen})$ with respect to both $\ol\M_i$
(Definition~\ref{Def: A_gen}). Let
\[
\ol\varphi_A: (\ol\M_1)_A\lra  (\ol\M_2)_A
\]
be a homomorphism of sheaves of monoids. Then the functor
\[
\operatorname{LMor}_Y^{\ol\varphi_A}:
(Y'\mapright{f} Y)\longmapsto \big\{\varphi:(Y',f^*\M_2)\to
(Y',f^*\M_1)\,\big|\, \text{$\varphi$ is compatible with
$\ol\varphi_A$}\big\}
\]
is represented by a scheme $\mathrm{LMor}_Y^{\ol\varphi_A}$ 
of finite type and affine over $Y$.
\end{lemma}

\proof
It is sufficient to prove the statement on an \'etale open cover of
$Y$, since we can then use descent for affine morphisms
(\cite{SGA1}, VIII, Theorem~2.1) to obtain a scheme over $Y$. Thus
we can assume that we in fact have charts $\psi_i:P_i\to
\Gamma(Y,\M_i)$ for the two log structures. We can also assume that
$\ol\varphi_A$ is induced by a homomorphism of monoids
$\ol\varphi_A:P_1\to P_2$.

Let $p_1,\ldots,p_n \in P_1$ be a generating set for $P_1$ as a
monoid.  Consider the sheaf of finitely generated $\O_Y$-algebras
\[
\shF_Y:=\O_Y[P_1^{\gp}]/\langle
\alpha_1(\psi_1(p_i))-z^{p_i}\alpha_2(\psi_2(\ol\varphi_A(p_i)))\,|\,
1\le i\le n\rangle.
\]
Then the desired scheme is
$\mathrm{LMor}_Y^{\ol\varphi_A} :=\bSpec \shF_Y$.

To see that this is the correct scheme, suppose $f:Y'\to Y$
is given. We wish to show that giving a commutative diagram of
schemes
\[
\xymatrix@C=30pt
{ Y'\ar[r]^<<<<<g\ar[dr]_f& \bSpec\shF_Y\ar[d]\\
&Y}
\]
is equivalent to giving a log morphism $\varphi:(Y',f^*\M_2)
\to (Y',f^*\M_1)$ which is the identity on $Y'$ and lifts
$f^*(\ol\varphi_A): f^*\ol\M_1\to f^*\ol\M_2$. 
Of course giving $g$ is equivalent to giving a section of
$(\bSpec\shF_Y)\times_Y Y'$ over $Y'$. But 
\begin{align*}
(\bSpec\shF_Y)\times_Y Y'{} =&\\
\bSpec \O_{Y'}&[P_1^{\gp}]
\big/\big\langle f^*\big(\alpha_1(\psi_1(p_i))\big)-z^{p_i}
f^*\big(\alpha_2(\psi_2(\ol\varphi_A(p_i)))\big)
\,\big|\,1\le i\le n\big\rangle,
\end{align*}
and the latter scheme is $\bSpec\shF_{Y'}$ associated to the data
$(Y',f^*\M_1)$, $(Y',f^*\M_2)$ with charts $f^*(\psi_i)=f^\flat\circ
\psi_i:P_i \to \Gamma(Y',f^*\M_i)$. Thus, without loss of
generality, we can assume that $Y=Y'$ and $f$ is the identity.

Now giving $\varphi:(Y,\M_2)\to(Y,\M_1)$ lifting
$\ol\varphi_A$ is equivalent to specifying $\varphi^\flat$. In
order for $\varphi^\flat$ to lift $\ol\varphi_A$, there must be
a map $\eta:P_1\to\Gamma(Y,\O_{Y}^{\times})$ with the
property that for all $p\in P_1$,
\[
\varphi^\flat(\psi_1(p))=\eta(p)\cdot \psi_2(\ol\varphi_A(p)).
\]
Giving $\eta$ completely determines $\varphi^\flat$. In addition,
$\varphi^\flat$ is a homomorphism of monoids if and only if $\eta$
is a homomorphism, and since $\eta$ takes values in the group
$\O_{Y}^{\times}$, specifying $\varphi^\flat$ is equivalent to
specifying a section of $\bSpec \O_Y[P_1^{\gp}]$. Indeed, giving a
section of $\bSpec\O_Y[P_1^{\gp}]$ over $Y$ is the same as giving
a morphism $Y\to \Spec\kk[P_1^\gp]$, which in turn is the same
as giving an element of $\Hom(P_1,\Gamma(Y,\O_Y^{\times}))$.

Secondly, since $\varphi^*=\id$, we must have $\alpha_1=
\alpha_2\circ\varphi^\flat$, so for each $p\in P_1$, we must have
\[
\alpha_1(\psi_1(p))=\alpha_2(\varphi^\flat(\psi_1(p)))=
\eta(p)\cdot \alpha_2(\psi_2(\ol\varphi_A(p))).
\]
If this holds for each $p_i$, it holds for all $p$. Thus a section
of $\bSpec\O_Y[P_1^{\gp}]$ over $Y$ determines a morphism of log
structures if and only if
it lies in the subscheme determined by the equations
\[
\alpha_1(\psi_1(p_i))-z^{p_i}\alpha_2(\psi_2(\ol\varphi_A(p_i))),
\]
demonstrating the result.
\qed

\begin{lemma}
\label{logmapmoduli2}
Let $Y\to W$ be a projective, separated morphism of schemes. Let
$\ol w \to W$ be a geometric point, and assume any $\ol x\in |Y|$
fulfills $(A_\mathrm{gen})$ with respect to both $\ol\M_i$ for
$A=Y_{\ol w}$ (Definition~\ref{Def: A_gen}). Then for a homomorphism
$\ol\varphi_A: (\ol\M_1)_{\ol w} \to (\ol\M_2)_{\ol w}$ of sheaves
of monoids the functor
\[
\operatorname{LMor}_{Y/W}^{\ol\varphi_A}:
(V\to W)\longmapsto
\big\{\varphi: (\M_1)_V
\to (\M_2)_V \,\big|\, \text{$\varphi$ is
compatible with $\ol\varphi_A$} \big\}
\]
is represented by a scheme $\mathrm{LMor}_{Y/
W}^{\ol\varphi_A}$ of finite type over $W$.
\end{lemma}

\proof
Let $Z=\mathrm{LMor}_{Y}^{\ol\varphi_A}$. By Lemma
\ref{logmapmodulilemma}, $\operatorname{LMor}_{Y/
W}^{\ol\varphi_A}$ is isomorphic to the functor
\[
(V\to W)\longmapsto \big\{\hbox{sections
of $Z_V\to Y_V$}\big\}.
\]
This is precisely the functor of sections $\prod_{Y/W}
(Z/Y)$ discussed in \cite{Gr}, p.221-19 (see also \cite{Nt}), and it
is represented by an open subscheme of $\mathrm{Hilb}_{Z/
W}$. Furthermore, if $\shL$ is a relatively ample line bundle on
$Y$, then as $Z\to Y$ is affine, the pull-back of $\shL$ to
$Z$ is also (trivially) relatively ample over $W$. We use
this ample line bundle to define Hilbert polynomials. Any section of
$Z_V \to Y_{V}$ then must have the same Hilbert polynomial
with respect to $\shL$ as the Hilbert polynomial of $\shL$ on
$Y_{V}$, so in fact $\Pi_{Y/W} (Z/Y)$ defines an open
subscheme of $\mathrm{Hilb}_{Z/W}^P$, for $P$ this fixed
Hilbert polynomial. Thus $\Pi_{Y/W} (Z/Y)$ is represented
by a scheme of finite type over $W$.
\qed
\bigskip

We are now in position to prove Proposition~\ref{LMor is
representable}.
\medskip

\noindent
\emph{Proof of Proposition~\ref{LMor is representable}.}
Since $\M_1$, $\M_2$ are sheaves in the \'etale topology,
$\mathrm{LMor}_{Y/W}(\M_1,\M_2)$ is clearly a sheaf in the \'etale
topology. It remains to prove the local representability statement
of \cite{knutson}, Ch.II, Definition~1.1,b, for
$\mathrm{LMor}_{Y/W}(\M_1,\M_2)$. For any geometric fibre $A=Y_{\ol
w}$ take an \'etale neighbourhood $U=U(\ol w)\to W$ of $\ol w$ such
that $Y_U$ fulfills the assumptions of Lemma~\ref{logmapmoduli2}.
Then for any $\ol \varphi_{\ol w}$ we have a scheme
$\mathrm{LMor}_{Y_U/U}^{\ol\varphi_{\ol w}}$ of
finite type over $W$. We claim that the natural functor
\[
Z:=\coprod_{\ol \varphi_{\ol w}}\mathrm{LMor}_{Y_U/
U}^{\ol\varphi_{\ol w}} \lra \mathrm{LMor}_{Y/W}(\M_1,\M_2)
\]
is schematic and an \'etale surjection. This means explicitly that
if $V$ is a scheme and $V\to \mathrm{LMor}_{Y/W}(\M_1,\M_2)$ is a
morphism of stacks then the fibre product in the sense of
$2$-categories $V\times_{\mathrm{LMor}_{Y/W} (\M_1,\M_2)} Z$ is
represented by a scheme, and the projection to $V$ is an \'etale
surjection. We claim that for one
$\mathrm{LMor}_{Y_U/U}^{\ol\varphi_{\ol w}}$ the representing scheme
is an open subset of $V\times_W U$. In fact, the functor $V\to
\mathrm{LMor}_{Y/W}(\M_1,\M_2)$ says that we fix a morphism of log
structures $\varphi: (\M_1)_V\to (\M_2)_V$. Now a functor
\[
\psi: T\lra V\times_{\mathrm{LMor}_{Y/W} (\M_1,\M_2)}
\mathrm{LMor}_{Y_U/ U}^{\ol \varphi_{\ol w}}
\]
from a scheme $T$ is nothing but (i) a morphism $T\to V$ and (ii) a
morphism $T\to \mathrm{LMor}_{Y_U/ U}^{\ol \varphi_{\ol w}}$ such
that (iii) the compositions with the morphisms to
$\mathrm{LMor}_{Y/W} (\M_1,\M_2)$ coincide. Note that (i)~provides a
pull-back $\varphi_T: (\M_1)_T\to (\M_2)_T$ of $\varphi$, (ii)~gives
a morphism of schemes $T\to U$ and a morphism of log structures
$\varphi': (\M_1)_T\to (\M_2)_T$ with $\ol{\varphi'}$ induced by
$\ol\varphi_{\ol w}$, and (iii)~says $\varphi'= \varphi_T$. Thus
$\psi$ is nothing but a factorization of the composition $T\to V\to
W$ through $U$ along with the information that $\ol\varphi_T$ is
induced by $\ol\varphi_{\ol w}$. Given $\varphi$ the latter
condition defines an open subset of $V\times_ W U$. Thus $\psi$ is
canonically identified with a homomorphism from $T$ to an open
subset of $V\times_W U$. This proves the claim. Note that since
$U\to W$ is \'etale so is the projection $V\times_W U\to V$.
Finally, surjectivity of $V\times_{\mathrm{LMor}_{Y/W} (\M_1,\M_2)}
Z\to V$ follows from the fact that we took the union over all
$\ol\varphi_{\ol w}$.
\qed

\section{Boundedness}
\label{Sect: Boundedness}

The aim of this section is to identify parts of $\MM(X)$ that are of
finite type. The main results are Theorem~\ref{quasi-generated
boundedness} in \S\ref{Par: Combinatorial finiteness} and
Theorem~\ref{Thm: Boundedness} in \S\ref{Par: Boundedness}.


\subsection{Finiteness of combinatorial types}
\label{Par: Combinatorial finiteness}
In Definition~\ref{Def: types in GS},(2) we defined the type of a
stable log map over a geometric point. It is given by the dual
intersection graph $\Gamma_{\ul C}$ of the domain and data $\mathbf{u}
=\big\{(u_p)_p,(u_q)_q\big\}$. Morally the $u_p: P_p\to\NN$ tell the
order of contact with the toric divisors in a local chart for the
log structure on $X$. Since by log smoothness these orders stay
locally constant in families of stable log maps, they are part of
the data distinguishing connected components of $\MM(X)$.

\begin{definition}\label{log classes}
A \emph{class} $\beta$ of stable log maps to $X$ consists of the following.
\begin{enumerate}
\item[(i)]
The data $\ul\beta$ of an underlying ordinary stable map,
that is, the genus $g$ of $\ul C$, the number $k$ of marked points,
and data $A$ bounding the degree, e.g.\ as described in
\cite{behrendmanin}, p.12.\footnote{In
\cite{behrendmanin} $A$ is defined by the degree function on the
cone of isomorphism classes of ample invertible sheaves; if
$\kk\subset\CC$ one might prefer prescribing a class in the
singular homology group $H_2(\ul X_\CC,\ZZ)$ of
the associated complex variety $\ul X_\CC$.}
\item[(ii)]
Strict closed embeddings $Z_1,\ldots,Z_k\subset X$, with $\ul Z_i$
carrrying the reduced induced scheme structure, together with
sections $s_i\in\Gamma (\ul Z_i,(\ol \M_{Z_i}^\gp)^*)$. We call $\beta$
\emph{maximal} if none of the $s_i$ extends to any strictly larger
subset of $X$.
\end{enumerate}
A stable log map $(C/W,\mathbf x, f)$ \emph{is of class $\beta$} if
the underlying ordinary stable map is of type $(g,k,A)$, and if for
any $i$ we have $\im(\ul f\circ x_i)\subset \ul Z_i$
and for any geometric marked point $\ol w\to
\ul W$ the map 
\begin{equation}\label{u_p is s_i}
\ol\M_{Z_i,\ul f(x_i(\ol w))}=(\ul f^*\ol\M_X)_{x_i(\ol w)} \stackrel{\ol
f^\flat}{\lra} \ol \M_{C,x_i(\ol w)} =\ol \M_{W,\ul
w}\oplus\NN\stackrel{\pr_2}{\lra} \NN
\end{equation}
equals the germ of $s_i$ at $\ul f(x_i(\ol w))$.

The substack of $\MM(X)$ of stable log maps of class $\beta$ is
denoted $\MM(X,\beta)$.
\end{definition}

Note that the composition \eqref{u_p is s_i} is denoted $u_p$ in
other parts of the text.

\begin{remark}
The matching condition $\im(\ul f\circ x_i)\subset\ul Z_i$ clearly
defines a closed algebraic substack of $\MM(X)$, while the remaining
conditions are open. In particular, $\MM(X,\beta)$ is also an
algebraic stack locally of finite type over $\ul S$.

Moreover, if $\beta$ is maximal then $\MM(X,\beta)$ is an
\emph{open} substack of $\MM(X)$. In fact, the maximality condition
says that if $z\in \ul Z_i$ lies in the closure of $y\in\ul
X\setminus\ul Z_i$ then $s_i$ does not factor over the generization
map
\[
\ol\M_{Z_i,z}^\gp= \ol\M_{X,z}^\gp\lra \ol\M_{X,y}^\gp,
\]
for otherwise $s_i$ extends to $\ul Z_i\cup \cl(y)$.
Thus if $(C/W,\mathbf x,f)$ is a stable log map and the
composition~\eqref{u_p is s_i} equals the germ of $s_i$ at a point
$\ol w\to \ul W$ then the same is true in a whole neighbourhood.
\end{remark}

A necessary condition for boundedness of $\MM(X,\beta)$ is that only
finitely many types of stable log maps to $X$ of class $\beta$
occur. Unfortunately, we have been unable to prove this in complete
generality; so far we have only been able to prove finiteness given
certain assumptions on $X$. On the other hand, we also couldn't
find an $X$ for which this finiteness does not hold. Thus we believe
that the following definition is in fact empty, at least locally
over $\mathbf M(\ul X)$.

\begin{definition}\label{Def: Combinatorial finiteness}
A class $\beta$ of stable log maps is called \emph{combinatorially
finite} if the set of types of stable log maps of class $\beta$ is
finite.
\end{definition}

Thus in general, if one wishes to deal with
log Gromov-Witten invariants for an $X$ for which we do not
prove finiteness below, one will have to check finiteness for
that $X$. However, we think that the cases discussed below will
cover most, if not all, applications of log Gromov-Witten theory.

We fix in this section an ordinary stable map $\ul f:\ul
C\to \ul X$ over $\Spec\kappa$, and we consider all possible
types of liftings of such a map to $f:C\to X$ over the
standard log point $(\Spec\kappa,\NN)$.

We introduce several conditions a log scheme can satisfy which
will be useful for proving boundedness.

\begin{definition}
We say a sheaf of monoids $\ol\M$ on (the Zariski site of) a
scheme $\ul Y$ is \emph{almost generated} if the maps
\[
\Hom(\ol\M_y, \RR_{\ge 0})\to \Hom(\Gamma(\ul
Y,\ol\M),\RR_{\ge 0})
\]
are injective for all $y\in \ul Y$. We say a log scheme $Y$ is
\emph{almost generated} if $\ol\M_Y$ is almost generated.

We say a sheaf of monoids $\ol\M$ on a scheme $\ul Y$ is
\emph{quasi-generated} if
\[
\Hom(\ol\M_y, \RR)\to \Hom\big(\Gamma(\ul Y,\ol\M^\gp),\RR\big)
\]
are injective for all $y\in \ul Y$. We say a log scheme $Y$ (with a
log structure on the Zariski site) is \emph{quasi-generated} if
$\ol\M_y$ is quasi-generated.
\end{definition}

\begin{remark}
1)\ If $\ol\M$ is a sheaf of fine monoids then being almost
generated is equivalent to saying that for any $y\in \ul Y$ the
image of the restriction map $\Gamma(\ul Y,\ol \M)\to\ol\M_y$ spans
${\ol \M}_y^\gp\otimes_\ZZ\RR$. In fact, write $P=\ol\M_y$ and
$Q=\im\big(\Gamma(\ul Y,\ol\M)\to \ol\M_y\big) \subset P$. If
$\ol\M$ is almost generated then $\Hom(P, \RR_{\ge0}) \to \Hom(Q,
\RR_{\ge0})$ is injective. But then, since $P$ and $Q$ are fine
monoids, also the induced map of associated groups
\[
\Hom(P,\RR)=\Hom(P, \RR_{\ge0})^\gp
\lra \Hom(Q, \RR_{\ge0})^\gp= \Hom(Q,\RR)
\]
is injective. Thus $Q$ spans $P^\gp\otimes_\ZZ\RR$.

Conversely, if $\Hom(P, \RR) \to \Hom(Q, \RR)$ is injective then so is
the restriction to $\Hom(P,\RR_{\ge 0})$.\\[1ex]
2)\ Similarly, a sheaf of fine monoids $\ol\M$ is
quasi-generated iff for any $y\in\ul Y$
the image of $\Gamma(\ul Y,\ol\M^\gp)\to \ol\M_y^\gp$ spans
$\ol\M_y^\gp\otimes_\ZZ\RR$.
\end{remark}

\begin{examples}
\label{globalgeneratedexamples}
(1) The condition that $X$ be Deligne-Faltings arises in the work of
Abramovich and Chen \cite{abramovich etal}, \cite{chen}. This means
that there is a surjection $\ul{\NN}^r\to \ol\M_X$ for some
$r$. The condition of $X$ being almost generated is strictly weaker.
For example, suppose $\ul X$ is a surface which is singular only at
a point $p\in \ul X$, where $\ul X$ has an $A_{e-1}$ singularity,
$e\ge 2$. Let $D=D_1\cup D_2$ be a divisor in $\ul X$ with $D_1,D_2$
irreducible and $D_1\cap D_2=\{p,q\}$, where $q\in \ul X$ is a
smooth point of $\ul X$. Assume that near $p$ we have $\ul X$
locally given by the equation $xy=t^e$, and $D_1\cup D_2$ is locally
given by $t=0$. Let $\ul X$ be given the divisorial log structure
induced by $D$. Then one checks easily that $\Gamma(\ul
X,\ol\M_X)\cong \ol\M_{X,p}=S_e$, the submonoid of $\NN^2$ generated
by $(e,0)$, $(0,e)$, and $(1,1)$. In particular, $\Gamma(\ul
X,\ol\M_X) \to \ol\M_{X,q}=\NN^2$ is not surjective. Thus $X$
cannot be Deligne-Faltings. However, one checks easily that $X$ is
almost generated. 

On the other hand, Deligne-Faltings log structures are always almost
generated, as a surjection $P=\NN^r\to \ol\M_x$ yields an
injection $\Hom(\ol\M_x,\RR_{\ge 0})\to \Hom(P,\RR_{\ge
0})$.

(2) It is not clear that being quasi-generated is weaker than being
almost generated. However, in some common situations, it is easier
to check. In fact, it suffices to find a group $M$ and a map
$M\to\Gamma(\ul Y,\ol\M^{\gp})$ such that for all $y\in\ul
Y$ the induced map
\[
\Hom(\ol\M_y, \RR)\to \Hom(M,\RR)
\]
is injective. For example, if $X$ is a toric variety with the
divisorial log structure defined by the toric divisors and $M$ is
the character lattice, there is a natural surjection
$\ul{M}\to \ol\M^{\gp}_X$.

(3) If $X$ is almost generated (quasi-generated), and $\ul f: \ul C
\to \ul X$ is an ordinary stable map of curves, then the
pull-back log structure $\ul f^*\M_X$ on $\ul C$ is almost generated
(quasi-generated).

(4) Suppose we are given an ordinary stable map $\ul f: \ul C\to \ul
X$ and  the dual intersection graph of $\ul C$ is a tree, for
example if $g(C)=0$. Then the pull-back log structure $\ul f^*\M_X$
is almost generated. In fact, even better, the map $\Gamma(\ul C,\ul
f^*\ol\M_X)\to (\ul f^*\ol\M_X)_x$ is surjective for every
$x\in \ul C$, so the pull-back log structure is Deligne-Faltings.

Indeed, with these assumptions, the sheaf $\ul f^*\ol\M_X$ is
entirely determined by the stalks $P_{\eta}$ and $P_x$ for special
points $x$, along with uniquely determined generization maps
$\chi_{\eta,x}: P_x\to P_\eta$ for every distinguished point $x$. To
specify a section of $\ul f^*\ol\M_X$, we just need to specify
elements $s_x\in P_x$ for all $x\in \ul C$ such that
$\chi_{\eta,x}(s_x)=s_{\eta}$ for $x\in D$, $D:=\cl(\eta)\subset \ul
C$. Now picking a point $x\in D$, and $s_x\in P_x$, set
$s_{\eta}=\chi_{\eta,x}(s_x)$. For every other point $x'\in D$,
$x'\not=x$, the generization map $\chi_{\eta,x'}$ is always
surjective, so we can choose $s_{x'}\in P_{x'}$ with $\chi_{\eta,x'}
(s_{x'})=s_{\eta}$. Some of these points $x'$ will be double points,
hence allowing us to define $s_{\eta'}$ for other generic points
$\eta'$. Continuing in this fashion, using the fact there are no
cycles in the dual intersection graph of $\ul C$, gives us a section
of $\ul f^*\ol\M_X$ whose germ at $x$ is the given $s_x$.
\end{examples}

Our arguments hinge on the following finiteness result from tropical
geometry:

\begin{proposition}
\label{NSfiniteness} 
Fixing a graph $\Gamma$, a lattice $N$, and weight vectors
$u_{(v,E)}\in N$ for every non-compact edge $E$ of $\Gamma$, there
are only a finite number of types of tropical curves with target
$N_{\RR}$ (as in Definition~\ref{Def: tropical curve}) with this
given $\Gamma$ and $u_{(v,E)}$.
\end{proposition}

\proof This is a weaker result than \cite{nisi}, Proposition 2.1.
\qed

\begin{theorem}
\label{globalgeneratedfiniteness}
Let $(\ul C/\Spec\kappa,\mathbf x,\ul f)$ be an ordinary stable map and
suppose $\ol\M:=\ul f^*\ol\M_X$ is almost generated. Then there is
only a finite number of types of log curves with the given
underlying ordinary stable map.
\end{theorem}

\proof
Let 
\[
M=\Gamma(\ul C,\ol\M)^{\gp}.
\]
Let $N=\Hom(M,\ZZ)$, $N_{\RR}= N\otimes_{\ZZ} \RR$. Clearly
$M$ is spanned by the submonoid $\Gamma(\ul C,\ol\M)\subset M$. Thus
also the dual submonoid $\Gamma(\ul C,\ol\M)^{\vee}\subset N$ spans
the dual space $N$, and such dual submonids are also sharp. In
particular, $\Gamma(\ul C,\ol\M)^{\vee}$ coincides with $K\cap N$ for
some strictly convex rational polyhedral cone $K$ in $N_{\RR}$.

Observe that since we can pull-back any stable log map to a standard
log point it is enough to bound the types of tropical curves over
standard log points. We thus consider now a stable log map
$(C/(\Spec\kappa,\NN) ,\mathbf x, f)$ over the standard log point. As
in \S\ref{par: Log maps over std log pt} this determines data
$V_\eta, e_q, u_p, u_q$. The intersection graph $\Gamma_{\ul C}$ of
$\ul C$ along with the data $\big\{(u_p),(u_q)\big\}$ is the type of
$(C/(\Spec\kappa,\NN) ,\mathbf x, f)$ by pull-back. Similar to
Discussion~\ref{Disc: Tropical curve} we now look at the associated
tropical curve. By the definition of almost generated, we obtain for
all $x\in \ul C$ inclusions
\[
\Hom(P_x,\RR_{\ge 0})\hookrightarrow K.
\] Thus $V_{\eta}\in P_\eta^\vee$, $u_p\in P_p^\vee$  live naturally
in $K$ and in fact  in the monoid $\Gamma(\ul C,\ol\M)^{\vee}
\subset N$.  We also have, for each irreducible component
$D=\cl(\eta)$ of $\ul C$, a map $\tau^X_\eta:\Gamma(\tilde
D,g^*\ol\M)\to\ZZ$, where $g:\tilde D\to \ul C$ is
the normalization of $D$ followed by inclusion into $\ul C$, see
\eqref{tau_eta}. We then have a composition of the pull-back map on
sections with $\tau^X_\eta$:
\[
M\to M_\eta:=
\Gamma(\tilde D,g^*\ol 
\M)^\gp \mapright{\tau^X_\eta} \ZZ.
\]
This composition determines an element
of $N$ which we also denote by $\tau^X_\eta$.

To build a tropical curve from this data  add to $\Gamma_{\ul C}$ a
number of unbounded edges:  for each vertex $v_{\eta}$  we attach an
unbounded edge, $E_{\eta}$, to $v_{\eta}$.  We then define
$h:\Gamma_{\ul C}\to N_{\RR}$ by
\[
h(v_{\eta})=V_{\eta};
\]
each edge $E_q$ with endpoints $v_{\eta_1}$ and $v_{\eta_2}$ is
mapped to the line segment joining $V_{\eta_1}$ and $V_{\eta_2}$;
and each edge $E_p$ with endpoint $v_{\eta}$ is mapped to the ray
with endpoint $V_{\eta}$ with direction $u_p$. Finally, we map the
ray $E_{\eta}$ to the ray with endpoint $V_{\eta}$ in the direction
defined by $\tau^X_\eta$. (If $\tau^X_\eta=0$, the edge is
contracted).

To give this the structure of a tropical curve, one also needs to
assign integral weight vectors to each flag $(v,E)$ of $\Gamma_{\ul C}$.
If $E=E_q$, we assign  the vector $\pm u_q$, with the sign chosen so
that $\pm u_q$ points away from $h(v)$. If $E=E_p$, we
associate the vector $u_p$, and if $E=E_{\eta}$, we associate the
vector $\tau^X_\eta$. These associated vectors are integral
tangent vectors to the image of the corresponding edge. Then the
tropical curve balancing condition  is just the statement that for a
given $\eta$,
\[
\tau_\eta^X+
\sum_{x} u_x=0
\]
where the sum is over all special points $x\in \cl(\eta)$. But this
is precisely the image of the equation of Proposition~\ref{Prop:
Balancing condition} in $M_\eta^{\vee}$, under the map
$M_\eta^{\vee}\to N$. Thus $h$ defines a balanced
tropical curve in $N_{\RR}$. 

Summing the balancing condition over all vertices gives a global
balancing condition involving all unbounded edges. Namely, in $N$,
we have
\[
\sum_\eta \tau_\eta^X+\sum_p u_p=0,
\]
where $\eta$ runs over all generic points of $\ul C$ and $p$ runs
over all marked points. Now $\sum_\eta \tau_\eta^X\in N$ is given,
completely specified by the original stable map $\ul f: \ul C\to \ul
X$ and independent of $f^\flat$, and $u_p\in \Gamma(\ul
C,\ol\M)^{\vee}$, which as observed at the beginning of the proof,
is the set of integral points of a strictly convex cone in
$N_\RR$. Thus there is only a finite number of possibilities for
writing $-\sum_\eta \tau_\eta^X$ as a sum of such $u_p$'s. This
shows finiteness of the choices of the $u_p$'s.

For given $(u_p)_p$ and $\tau_\eta^X$, Proposition
\ref{NSfiniteness} shows there is only a finite number of possible
combinatorial types of this tropical curve. This means there is
only a finite number of possibilities for the vectors $u_q\in N$. 
On the other hand, these vectors are images of $u_q\in
(P_q^{\gp})^*$. Since $(P_q^{\gp})^*$ injects into $N$, we
conclude we only have a finite number of allowable types.
\qed

\bigskip

We have a slightly weaker result in the quasi-generated case, where
we need to fix the $u_p$'s, that is, the class.

\begin{theorem}\label{quasi-generated boundedness}
In the situation of Theorem~\ref{globalgeneratedfiniteness} suppose
that $\ol\M=\ul f^*\ol\M_X$ is only quasi-generated. Then if the
$u_p$'s are fixed, there is only a finite number of types of log
curves with the given underlying ordinary stable map.
\end{theorem}

\proof
The argument is essentially the same as the proof of the previous
theorem, this time taking $M=\Gamma(\ul C,\ol\M^{\gp})$. Then via
the same construction, one obtains a tropical curve in $N_{\RR}$.
However, one no longer has all $u_p$ living in a strictly convex
cone in $N_{\RR}$, so we cannot use this to control the vectors
$u_p$. However if we assume the $u_p$ are given, then
Proposition~\ref{NSfiniteness} still gives the needed finiteness.
\qed
\bigskip

We next address what can be accomplished in more general situations.
In Appendix~\ref{App: Tropicalization} we introduce the
\emph{tropicalization} $\Trop(X)$ of $X$ as a natural target space
for tropical curves. While this quite generally provides a tropical
curve, it can be difficult to use the balancing condition in the
often strange space $\Trop(X)$ in a useful way to prove finiteness.

Here is a special case that we state without proof, where it is
still possible to prove finiteness. This is likely to be a good
model for the type of application we have in mind, in which one
considers varieties degenerating to unions of reasonably simple
varieties.

\begin{theorem}
\label{lastboundednesstheorem}
Let $\ul X$ be a scheme with $D\subset \ul X$ a divisor inducing a
divisorial log structure on $\ul X$. Suppose that this makes $X$ an fs
log scheme, log smooth over $\Spec\kappa$, and suppose furthermore:
\begin{enumerate}
\item For each irreducible component $\ul Y$ of $D$, the restriction of
the log structure of $X$ to $\ul Y$ is almost generated. Furthermore, for
$ y\in \ul Y$, the inclusion $\Hom(\ol\M_{X,
y},\RR_{\ge 0})\to \Hom(\Gamma(\ul Y, \ol\M_X),\RR_{\ge
0})$ is an inclusion of faces.\footnote{Note this holds, for example, if the
log structure on $Y$ is Deligne-Faltings.}
\item $X$ is monodromy free (Definition~\ref{Def: monodromy free}).
\end{enumerate}
Then given an ordinary stable map $\ul f:\ul C\to \ul X$ and
a collection of $u_p\in P_p^{\vee}$ for $p\in \ul C$ marked points,
there is only a finite number of possible types of log curves with
the given underlying ordinary stable map.
\qed
\end{theorem}

We were also able to prove boundedness without any further
hypotheses in the case of genus one, also stated here without proof.

\begin{theorem}\label{Thm: genus one}
Let $(\ul C/\Spec\kappa,\ul x,\ul f)$ be an ordinary stable map and
suppose the dual intersection graph $\Gamma_{\ul C}$ has genus at
most one. Then there is only a finite number of types of log curves
with the given underlying ordinary stable map.
\qed
\end{theorem}


\subsection{Stable log maps of constant type and boundedness}
\label{Par: Boundedness}

Here is the main result of this section.

\begin{theorem}
\label{Thm: Boundedness}
For any combinatorially finite class $\beta$ of stable log maps to
$X$, the algebraic stack $\MM(X,\beta)$ is of finite type over $\ul S$.
\end{theorem}

By Corollary~\ref{MM(X) is algebraic log stack} we already know that
$\MM(X)\to \ul S$ is locally of finite type. Moreover, the stack of
ordinary stable maps of fixed class $\mathbf M(\ul X,\ul\beta)$ is
of finite type over $\ul S$. So to finish the proof of
Theorem~\ref{Thm: Boundedness} it remains to show that for any
morphism $\ul W\to \mathbf M(\ul X)$ from a quasi-compact scheme
$\ul W$ the fibre product $\ul W\times_{\mathbf M(\ul X)}
\MM(X,\beta)$ is quasi-compact. By \cite{champs}, Corollaire~5.6.3
this is equivalent to showing that the topological space of
geometric points $\big|\ul W\times_{\mathbf M(\ul X)}
\MM(X,\beta)\big|$ is quasi-compact. We do this by a stratawise
approach.

\begin{definition}
A \emph{weak covering} of a topological space $Z$ is a collection of
subsets $\{A_i\}_{i\in I}$ with the following property:
For any $z\in Z$ there exists $i\in I$ with $\cl(z)\cap
A_i\neq\emptyset$. 
\end{definition}

\begin{lemma}\label{Lem: qc criterion}
Let $Z$ be a topological space weakly covered by finitely many
quasi-compact subsets $A_i\subset Z$, $i\in I$. Then $Z$ is
quasi-compact.
\end{lemma}

\proof
Let $U_j\subset Z$, $j\in J$, be open subsets covering $Z$. By
quasi-compactness of the $A_i$, for any $i\in I$ there exist
finitely many $U_j$ covering $A_i$. Hence, since $I$ is finite,
there exists a finite subset $J'\subset J$ with $A_i\subset
\bigcup_{j\in J'} U_j$, for any $i$. We claim $Z=\bigcup_{j\in J'}
U_j$. In fact, let $z\in Z$. Then by the weak covering assumption
there exists $i\in I$ with $\cl(z)\cap A_i\neq\emptyset$. Thus by
the choice of $J'\subset J$ there exists $j\in J'$ with $\cl(z)\cap
U_j\neq\emptyset$, and then $z\in U_j$ for $U_j\subset Z$ is open.
Thus already finitely many of the $U_j$ cover $Z$.
\qed
\medskip

Our strata $A_i$ will be defined by taking locally trivial families
of ordinary stable maps.

\begin{definition}\label{Def: cc}
An ordinary stable map $(\ul\pi: \ul C\to \ul W,\mathbf{x}, \ul f)$
over an integral scheme $\ul W$ is called \emph{combinatorially
constant} if the following conditions are satisfied, where we write
$\M:= \ul f^*\M_X$ as usual.
\begin{enumerate}
\item
If $g:\tilde{\ul C}\to \ul C$ is the normalization, then the
composition $\ul \pi\circ g$ is a smooth map, and there are
pairwise different sections $y_q$ of $\ul\pi$ with $\bigcup_q
\im(y_q)= \operatorname{crit}(\ul\pi)$.
\item
Each irreducible component $\ul C_\eta\subset \ul C$ is geometrically
connected and there is a section $\sigma_\eta$ of $\ul\pi|_{\ul
C_\eta}$ with $\im(\sigma_\eta)$ disjoint from $\im(\mathbf x)$ and
from $\operatorname{crit}(\ul \pi)$.
\item
For any $i,q$ the sheaves $x_i^*\ol\M$ and $y_q^*\ol\M$ are constant,
and there are charts of log structures on $\ul W$
\[
\ul{P_\eta}\lra \sigma_{\eta}^*\M,\quad
\textstyle\prod_q\ul\NN\lra \M_W^0
\]
inducing isomorphisms $\ul{P_\eta} \simeq \sigma_\eta^*\ol\M$ and
$\prod_q \ul\NN\simeq \ol\M_W^0$. Here $\M_W^0$ is the basic log
structure for the pre-stable curve $\ul C/\ul W$.
\end{enumerate}
\end{definition}

Note that (3) defines isomorphisms $ \sigma_\eta^* \M \simeq \ul
P_\eta\times\O_{\ul W}^\times$ and $\M_W^0\simeq\prod_q \ul\NN
\times \O_{\ul W}^\times$. Note also that for a combinatorially
constant stable map $(\ul C/\ul W,\mathbf{x}, \ul f)$ the existence
of the sections $\sigma_\eta$ and Stein factorization imply that the
geometric fibres of $\ul\pi\circ g|_{\ul C_\eta}$ are connected.
Thus all geometric fibres of $\ul \pi$ have the same dual
intersection graph $\Gamma$ as the generic fibre, and $\ul C$ is
obtained by gluing together families of smooth, connected curves
along pairs of sections. Moreover, a type $\mathbf{u}=
\big\{(u_p),(u_q)\big\}$ at the generic fibre induces a type for the
fibre over every geometric point $\ol w\to \ul W$. By abuse of
notation we call $\mathbf{u}$ a \emph{type for $(\ul C/\ul
W,\mathbf{x}, \ul f)$}.
\medskip

Let $\ul W\to \mathbf M(\ul X)$ be the morphism defined by a
combinatorially constant stable map $\mathfrak f= (\ul C/\ul
W,\mathbf{x}, \ul f)$. Then $\ul W\times_{\mathbf M(\ul X)} \MM(X)$
is isomorphic as a stack to the (non-full) subcategory
$\MM(X,\mathfrak f)\subset \MM(X)$ with objects stable log maps with
underlying ordinary stable map obtained by pull-back from $\mathfrak
f$, and with morphisms on the underlying schemes induced by the
identity on $\mathfrak f$. Now the type of a stable log map in
$\M(X,\mathfrak f)$ is locally constant. Hence we have a
decomposition into disjoint open substacks
\begin{equation}\label{decomposition according to types}
\ul W\times_{\mathbf M(\ul X)} \MM(X) \simeq
\coprod_{\mathbf u} \MM(X,\mathfrak f,\mathbf u)
\end{equation}
according to the type $\mathbf u$ for $\mathfrak f$.

The following lemma shows that provided $\ul W$ is reduced, the log
structure on the base of a stable log map in $\MM(X,\mathfrak
f,\mathbf u)$ is locally constant.

\begin{lemma}
\label{reducedtrivlog}
Let $W$ be a reduced scheme with a fine log structure
$\M$ such that $\overline{\M}$ is a constant sheaf.
Then \'etale locally, $\M\simeq O_{\ul W}^{\times}\times\overline{\M}$.
\end{lemma}

\proof Let $Q$ be the stalk of $\overline{\M}$. Then a chart
for $\M$ takes the form, \'etale locally on $\ul W$,
\[
Q\lra \O_{\ul W}\qquad
m\longmapsto \begin{cases} 1&m=0;\\
0&m\not=0.
\end{cases}
\]
Indeed, a chart defining $\ul W$ must take every non-zero element
of $Q$ to a function which is not invertible at any point of its domain.
But as $\ul W$ is reduced, the only function which is non-invertible
at every point of an open set of $\ul W$ is $0$. Thus locally $\M$
takes the given form.
\qed

\begin{proposition}\label{Prop: Quasi-compactness for fixed type}
Let $\mathfrak f= (\ul C/\ul W,\mathbf{x}, \ul f)$ be a
combinatorially constant ordinary stable map over an integral,
qausi-compact scheme $\ul W$. Then for any type $\mathbf u$ for
$\mathfrak f$ the stack $\MM(X,\mathfrak f,\mathbf u)$ is
quasi-compact.
\end{proposition}

\proof
Let $Q$ be the basic monoid defined by $\mathbf u$. Let
$\M_1:=\M_{\ul W}^0$ be the basic log structure for the pre-stable
curve $\ul C/\ul W$, $\M_2:=Q\times \O_{\ul W}^\times$ the constant log
structure and
\[
\ol\varphi:\ol\M_1=\prod_q \ul\NN\lra  \ul Q= \ol\M_2
\]
the homomorphism coming from the definition of $Q$. To take care of
the domains of the relevant stable log maps we now look at $\ul
Y:=\mathrm{LMor}^{\ol\varphi}_{\ul W}$ from
Lemma~\ref{logmapmodulilemma}. By the universal property of $\M_W^0$
this scheme classifies log smooth structures on $\ul C_{\ul V}/\ul
V$ with trivialized base log structure $Q\times\O_{\ul V}^\times$.
Let $C_{\ul Y}/Y$ be the universal log smooth curve.

It remains to lift the pull-back $\ul f_{\ul Y}$ of $\ul f$ to a log
morphism. This is done by $\ul Z:=\mathrm{LMor}^{\ol\psi}_{{\ul
C}_{\ul Y}/\ul Y}$, where now the two log structures are $\M_1:=\ul f_{\ul
Y}^*\M_X$ and $\M_2:=\M_{C_{\ul Y}}$ on ${\ul C}_{\ul Y}$. The map
$\ol\psi: \ol\M_1\to\ol\M_2$ is again fixed by basicness. We
then obtain a universal basic stable log map $(C_{\ul Z}/ Z, \mathbf
x_{\ul Z}, f_{\ul Z})$. Since $\ul Z\to \ul W$ is of finite
type, $\ul Z$ is quasi-compact by quasi-compactness of $\ul W$.

Now the morphism $\ul Z\to \MM(X,\mathfrak f,\mathbf u)$ of
algebraic stacks thus obtained is an epimorphism. In fact, by
Lemma~\ref{reducedtrivlog} locally on $\ul W$ the log structure
$\M_W$ of a stable log map $(C/W,\mathbf x,f)\in
\Ob\big(\MM(X,\mathfrak f,\mathbf u)\big)$ is isomorphic to the
constant log structure $Q\times\O_{\ul W}^\times$, and hence
$(C/W,\mathbf x,f)$ is locally the image of the pull-back of the
universal basic stable log map over $\ul Z$. Thus the induced map of
spaces of geometric points $|\ul Z|\to |\MM(X,\mathfrak f,\mathbf
u)|$ is a continuous surjection, and quasi-compactness of
$|\MM(X,\mathfrak f,\mathbf u)|$ follows from quasi-compactness of
$|\ul Z|$.
\qed
\medskip

\begin{lemma}\label{Lem: cc stratification}
Let $({\ul C}/\ul W,\mathbf{x},\ul f)$ be an ordinary stable map
with $\ul W$ quasi-compact. Then there exist finitely many locally
closed, integral subschemes $\ul W_i\subset \ul W$ weakly covering
$\ul W$ and \'etale surjections $\ul{\tilde W}_i\to \ul W_i$ such
that the pull-back of $({\ul C}/\ul W,\mathbf{x},\ul f)$ to
$\ul{\tilde W}_i$ is combinatorially constant.
\end{lemma}

\proof
By quasi-compactness it suffices to construct the $\ul W_i$ locally.
Standard arguments for families of pre-stable curves provide the
$\ul W_i$ fullfilling Conditions~(1) and (2) in Definition~\ref{Def:
cc}. Then use the fact that $\M$ is a fine sheaf to also achieve~(3).
\qed
\bigskip

We are now in position to prove Theorem~\ref{Thm: Boundedness}.
\medskip

\noindent
\emph{Proof of Theorem~\ref{Thm: Boundedness}.}
By the discussion at the beginning of this subsection we have to
show that if $\ul W\to \MM(X,\beta)$ is a morphism from a
quasi-compact scheme $\ul W$ then $\big|\ul W\times_{\mathbf M(\ul
X)} \MM(X,\beta)\big|$ is quasi-compact. By the quasi-compactness
criterion Lemma~\ref{Lem: qc criterion} together with
Lemma~\ref{Lem: cc stratification} we can assume that $\ul W$ is
integral and that the ordinary stable map $\mathfrak f=({\ul C}/\ul
W,\mathbf{x},\ul f)$ defined by $\ul W\to \mathbf{M}(\ul X)$ is
combinatorially constant. Let $\mathbf{u}_1,\ldots,\mathbf{u}_s$ be
the finite list of types $\mathbf{u}$ of stable log maps in class
$\beta$ according to Definition~\ref{Def: Combinatorial finiteness}.
From \eqref{decomposition according to types} we have the decomposition
\[
\big|\ul W\times_{\mathbf{M}(\ul X)} \MM(X,\beta)\big| =
\coprod_{\mu=1}^s \big|\MM(X,\mathfrak f,\mathbf u_\mu)\big|.
\]
Now according to Proposition~\ref{Prop: Quasi-compactness for fixed
type} each component on the right-hand side is quasi-compact. Hence
$\ul W\times_{\mathbf{M}(\ul X)} \MM(X,\beta)$ is quasi-compact, as
had to be shown.
\qed

\section{Stable reduction}

In this section we prove a stable reduction theorem for basic stable
log maps. Throughout, $R$ is a discrete valuation ring over our base
scheme $\ul S$ with maximal ideal $\maxid$, residue field
$R/\maxid=\kappa$ and quotient field $K$. We assume $K$ is endowed
with a fine, saturated log structure over $\M_S$, hence defining a
log point $(\Spec K, Q_K)$ over $S$ for a toric monoid $Q_K$. The
closed point in $\Spec R$ is denoted $0$.

\begin{theorem}
\label{Thm: Stable reduction}\sloppy
Assume that $\ul X\to \ul S$ is proper. Let $\big(\pi_K: C_K\to
(\Spec K, Q_K), \mathbf{x}_K,f_K\big)$ be a basic stable log map to
$X$ over the log point $(\Spec K, Q_K)$. Then possibly after
replacing $K$ by a finite extension $\tilde K$ and $R$ by its
integral closure in $\tilde K$ and pulling back the stable log map
via $(\Spec\tilde K, Q_K)\to (\Spec K,Q_K)$, the following holds:
There exists a log structure $\M_R$ on $R$ over $S$ together with a
strict morphism $(\Spec K, Q_K)\to (\Spec R,\M_R)$, and a basic
stable log map $(C/(\Spec R,\M_R),\mathbf{x},f)$ to $X$
such that the restriction to $\Spec K$ is isomorphic to
$(C_K/(\Spec K,Q_K), \mathbf{x}_K,f_K)$.

Moreover, such an extension as a basic stable log map is unique up
to unique isomorphism.
\fussy
\end{theorem}

\begin{corollary}
\label{properness of the stack of stable log maps}
For any combinatorially finite class $\beta$ of stable log maps,
$\MM(X,\beta)$ is proper over $\ul S$.
\end{corollary}

\proof
Recall that a morphism of algebraic stacks is proper if it is
separated, of finite type and universally closed (\cite{champs},
Definition~7.11). Separatedness follows by the uniqueness part in
the theorem from the valuative criterion \cite{champs},
Proposition~7.8. Theorem~\ref{Thm: Boundedness} established that
$\MM(X,\beta)\to \ul S$ is of finite type. In view of the existence
part in Theorem~\ref{Thm: Stable reduction} universal closedness
then follows from the valuative criterion \cite{champs},
Theorem~7.10.
\qed
\medskip

The proof of Theorem~\ref{Thm: Stable reduction} is divided into
three steps presented in the following subsections.

\subsection{Extension on the level of ghost sheaves}
\label{par: Ghost sheaves}
By the stable reduction theorem for ordinary stable maps
\cite{fultonpandh} we may assume $(\ul C_K/K,\mathbf{x}_K, \ul f_K)$
is the restriction to $K$ of a marked stable map $(\ul\pi: \ul C\to
\Spec R,\mathbf{x}, \ul f)$ over $R$. This step may involve a base
change, but note that any two such extensions are uniquely
isomorphic. In particular, the central fibre $(\ul C_0/\Spec \kappa,
\mathbf{x}_0, \ul f_0)$ is uniquely determined. The objective in
this subsection is to derive the extension as a stable log map on
the level of ghost sheaves. As we will see in the proof of
Proposition~\ref{Prop: Uniqueness of type} this amounts to
determining the type of the central fibre.

The central fibre being a stable curve over a field, we use our usual
convention that $\eta$, $q$, $p$ denote generic points, nodes and
marked points of $\ul C_0$, respectively. For the corresponding
objects on $\ul C_K$ we use an index $K$. Thus the sheaf $\ul
f^*\ol \M_X$ is determined by the usual generization maps
$\chi_{\eta,p}:P_p\lra P_\eta$, $\chi_{\eta,q}:P_q\lra P_\eta$ on
$\ul C_0$, together with the following generization maps from $\ul
C_0$ to $\ul C_K$:
\[
\chi_{\eta,K}: P_\eta\lra P_{\eta_K},\quad
\chi_{p,K}:P_p\lra P_{p_K},\quad
\chi_{q,K}:P_q\lra P_{q_K}.
\]
Note that the closures of nodes $q_K\in \ul C_K$ may only give a
subset of the nodes of $\ul C_0$; for the other nodes of $\ul C_0$
we only have the composition $\chi_{\eta,K}\circ \chi_{\eta,q}:
P_q\to P_{\eta_K}$. Let us refer to the latter as \emph{isolated
nodes}. 

\begin{proposition}\label{Prop: Uniqueness of type}
There exists a type $\mathbf u=(u_x)_{x\in \ul C_0}$ for stable log
maps with underlying ordinary stable map $(\ul C_0 /\Spec
\kappa,{\mathbf x}_0, \ul f_0)$ such that the central fibre of any
extension to $\Spec R$ of $(C_K/(\Spec K, Q_K), {\mathbf x}_K, f_K)$
as a basic stable log map has type $\mathbf u$. Moreover, $\mathbf
u$ defines the extension on the level of ghost sheaves uniquely up
to unique isomorphism.
\end{proposition}

\proof
Let $(C/(\Spec R,\M_R), {\mathbf x}, f)$ be an extension of
$(C_K/(\Spec K, Q_K), {\mathbf x}_K, f_K)$. We are going to show
that this extension is determined uniquely on the level of ghost
sheaves by data computable without knowing the extension.

By basicness we have a commutative diagram
\[
\begin{CD}
\prod_\eta P_\eta \times\prod_q\NN @>>> Q\\
@VVV @VVV\\
\prod_{\eta_K} P_{\eta_K} \times\prod_{q_K} \NN
@>>> Q_K,
\end{CD}
\]
defining the basic monoids $Q$, $Q_K$ for the central and the
generic fibres, respectively. The left-hand vertical arrow in the
diagram is defined as the product of $\prod\chi_{\eta,K}$ and the
projection to the non-isolated nodes. The horizontal arrows are
quotients by subgroups $R_{\mathbf u}\subset \prod_\eta
P_\eta^\gp\times\prod_q \ZZ$ and $R_{\mathbf u_K}\subset
\prod_{\eta_K} P_{\eta_K}^\gp\times\prod_{q_K} \ZZ$, determined by
the respective types $\mathbf{u}$ and $\mathbf u_K$, followed by
saturation, see Construction~\ref{Q and basic log structure}.

Recall also from Discussion~\ref{Disc: Stalkwise} that the upper
horizontal map together with the types $u_p$ of the marked points
determine the maps of ghost sheaves $\ul\pi^*\ol\M_R \to \ol\M_C$,
$\ul f^*\ol\M_X\to\ol \M_C$  on the central fibre, and similarly for
the lower horizontal arrow and the generic fibre. Commutativity of
the diagram then establishes the generization maps between the
relevant sheaves, and hence determines the extension on the level of
ghost sheaves. Thus it suffices to find the type $\mathbf u =
\big\{(u_p)_p,(u_q)_q\big\}$ with $R_{\mathbf u}$ mapping to
$R_{\mathbf u_K}$.

If $p\in \ul C_0$ is a marked point there exists a unique marked
point $p_K\in \ul C_K$ with $p\in\cl(p_K)$. By the structure of log curves
(Theorem~\ref{Thm: structure of log
curves},ii) $u_p$ must then equal the composition
\[
P_p\stackrel{\chi_{p,K}}{\lra} P_{p_K}\stackrel{u_{p_K}}{\lra} \NN,
\]
where the second map is given by the type of the generic fibre.
Hence $u_p$ depends only on the generic fibre and $\ul f$.

Similarly, for a non-isolated node $q$ let $q_K$ be the
node of $\ul C_K$ with $q\in\cl(q_K)$.
Then the defining equation \eqref{eqn at node} for $u_q$,
\[
\varphi_{\ol\eta_2}\big(\chi_{\eta_2,q}(m)\big)- 
\varphi_{\ol\eta_1}\big(\chi_{\eta_1,q}(m)\big)=
u_q(m)\cdot\rho_q,
\]
with $\varphi:= \ol{f^\flat}$ and $m\in P_q$, generizes to the
defining equation for $u_{q_K}$:
\[
\varphi_{\ol\eta_{2,K}}\big(\chi_{\eta_{2,K},q_K}(m_K)\big)- 
\varphi_{\ol\eta_{1,K}}\big(\chi_{\eta_{1,K},q_K}(m_K)\big)=
u_{q_K}(m_K)\cdot\rho_{q_K},
\]
where $m_K:=\chi_{q,K}(m)$. Hence $u_q= u_{q_K} \circ
\chi_{q,K}$ is also determined a priori by the generic fibre.

So far we have just reversed the reasoning in the proof of
Proposition~\ref{Prop: Basicness is open}. The isolated nodes,
however, require different arguments. Recall that $u_q$ was defined
by the homomorphism
\[
\varphi_{\ol q}: P_q= \ol\M_{X,\ul f(q)}\lra \ol \M_{C,\ol q}=
Q\oplus_\NN \NN^2,\quad m\longmapsto \big(m',(a,b)\big)
\]
as $u_q(m)= b-a$ (for the correct orientation). Now while the
composition with generization to $\eta_K$ maps any $(0,(a,b))$ to
$0\in Q_K$, we can retrieve $b-a$ by working on the level of log
structures. In fact, choosing a compatible chart on an \'etale
neighbourhood of $q$, the generization map has the form
\begin{eqnarray*}
(Q\oplus_\NN \NN^2)\times\O_{\ul C,\ol q}^\times \simeq \M_{C, \ol q}
&\lra& \M_{C,\ol \eta_K}\simeq Q_K\times \O_{\ul C,\ol \eta_K}^\times\\
\big((m',(a,b)),h\big) &\longmapsto& (m'_K, z^a w^b \ul\pi^* g(m')h),
\end{eqnarray*}
with $z=0$, $w=0$ defining the two branches of $\ul C_0$ inside $\ul
C$ at $q$, and $g(m')\in K^\times$. Assume that $\ul C$ has an
$A_{e-1}$-singularity at $q$. Then letting $\ord_{\eta_1}$,
$\ord_{\eta_2}$ be the discrete valuations of these branches it
holds $\ord_{\eta_1}(z)=\ord_{\eta_2}(w)=e$, and hence $e\cdot(b-a)=
(\ord_{\eta_2}-\ord_{\eta_1})(z^aw^b h)$ and
$\ord_{\eta_1}(\ul\pi^*(g(m')))= \ord_{\eta_2}(\ul\pi^*(g(m')))$.
Thus $e\cdot u_q$ equals the following composition:
\begin{equation}\label{u_q at isolated nodes}
P_q\lra \M_{X,\ul f(q)} \stackrel{f^\flat}{\lra} \M_{C,\ol q}
\lra \M_{C,\ol \eta_K}\simeq Q_K\times \O_{\ul C,\ol \eta_K}^\times
\stackrel{\pr_2}{\lra} \O_{\ul C,\ol\eta_K}^\times
\stackrel{\ord_{\eta_2}-\ord_{\eta_1}}{\lra} \ZZ.
\end{equation}
The first map is a choice of chart for $\ul f^*\M_X$ at $q$. Finally
observe that the composition $\M_{X,\ul f(q)}\to \M_{C,\ol q} \to
\M_{C,\ol\eta_K}$ equals $\M_{X,\ul f(q)} \to \M_{X,\ul
f(\eta_K)}\mapright{f_K^\flat}\M_{C,\ol \eta_K}$. Hence $u_q$ is
completely determined also for an isolated node $q$  by the generic
fibre and by the extension $\ul f$ of $\ul f_K$, already known to be
unique.
\qed

\subsection{The log structure on the base}
\label{par: M_R}
We now want to show that the extension of the log structure on the
base is uniquely defined. The decisive tool is the identification of
``fibrewise constant'' subsheaves $\M(\eta)$ of $\ul f^*\M_X$ on
certain open subsets of $\ul C$ as follows. There is one open subset
$U(\eta)$ for each generic point $\eta\in\ul C_0$. To define
$U(\eta)$ let $A\subset \ul C_0$ be the set of non-special points of
$\cl(\eta)\subset\ul C_0$. Then $U(\eta)\subset\ul C$ is the set of
generizations of points in $A$. Said differently, if $D\subset \ul
C$ is the irreducible component of $\ul C$ containing $\eta$ then
$U(\eta)$ is obtained from $D$ by removing all special points and all
irreducible components of $\ul C_0$ not containing $\eta$. Thus
$U(\eta)$ is a smooth open neighbourhood of $\eta$, and $\ul
f^*\ol\M_X|_{U(\eta)}$ has only two interesting stalks, at $\eta$
and at its generization $\eta_K$, see Remark~\ref{Rem: structure of
M}. Hence $\ul f^*\ol\M_X|_{U(\eta)}$ is completely determined by
the generization map $P_\eta \to P_{\eta_K}$. In particular, $\ul
f^*\ol\M_X|_{U(\eta)}$ is globally generated with global section
space equal to $P_\eta$. Note also that each $x\in U(\eta)$ is in
the closure of $\eta_K$, so we may restrict germs of sections to
$\eta_K$. The notation for the restriction of $m$ is $m|_{\eta_K}$.

\begin{definition}
The \emph{fibrewise constant subsheaf} $\M(\eta)\subset \ul
f^*\M_X|_{U(\eta)}$ is the sheaf in the Zariski topology defined by
the condition
\[
m\in \M(\eta)_{x}\quad \Longleftrightarrow\quad
f_K^\flat(m|_{\eta_K})\in \im\big( \Gamma(\M_K)
\mapright{\pi_K^\flat} \M_{C,\ol \eta_K}  \big).
\]
\end{definition}

Since $C_K\to (\Spec K, Q_K)$ is strict at $\ol\eta_K$, the map
$\pi_K^\flat$ is injective and there is a canonical map
\[
\beta_\eta: \M(\eta)_\eta \lra \Gamma(\M_K),\quad
m\longmapsto (\pi_K^\flat)^{-1} \big(f_K^\flat(m|_{\eta_K})\big).
\]
Note the composition of $\beta_\eta$ with the structure
homomorphism $\Gamma(\M_K)\to K$ factors over a map
\[
\alpha_\eta: \M(\eta)_\eta\lra R.
\]
In fact, if $m\in\M(\eta)_\eta$ maps to $a\in K$ then
$\pi^*(a)$ is in the image of the structure homomorphism
$\M(\eta)_\eta\to \O_{\ul C,\eta}$ and hence $a$ is regular at $0$. 

Calling $\M(\eta)$ fibrewise constant is justified by the fact that
the restrictions of $\M(\eta)$ to $\ul C_K\cap U(\eta)$ and to $\ul
C_0\cap U(\eta)$ are constant sheaves with fibres isomorphic to
$P_{\eta_K}\times K^\times$ and to $P_{\eta}\times R^\times$,
respectively. For our purposes only the stalk $\M(\eta)_\eta$ is
really relevant and it suffices to prove the following weaker
statement.

\begin{lemma}\label{Lem: fibrewise constant sections}
There is a (non-canonical) isomorphism
$\M(\eta)_\eta\simeq P_\eta\times R^\times$.
\end{lemma}

\proof
By strictness of $C_K\to (\Spec K,Q_K)$ at $\ol\eta_K$ there exists
an isomorphism $\M_{C,\ol \eta_K}\simeq Q_K\times \O_{\ul
C,\ol\eta_K}^\times$ such that $\pi_K^\flat(\Gamma(\M_K)) =
Q_K\times K^\times$. Let $\theta: P_\eta\to \M_{X,\ul f(\eta)}$ be a
chart, hence inducing an isomorphism $(f^*\M_X)_{\eta}\simeq
P_\eta\times \O_{\ul C,\eta}^\times$. Note that $\O_{\ul
C,\ol\eta_K}^\times\simeq \ZZ\times \O_{\ul C,\ol \eta}^\times$
since $R$ is a discrete valuation ring. We can thus define $g:
P_\eta\to \O_{\ul C,\ol \eta}^\times$ by the composition
\[
P_\eta\stackrel{\theta}{\lra} \M_{X,\ul f(\eta)}
\stackrel{f^\flat(\,.\,|_{\eta_K})}{\lra} \M_{C_K, \ol\eta_K}\simeq 
Q_K\times \O_{\ul C,\ol \eta_K}^\times \stackrel{\pr_2}{\lra}
\O_{\ul C,\ol \eta_K}^\times\simeq \ZZ\times\O_{\ul C,\ol \eta}^\times
\stackrel{\pr_2}{\lra} \O_{\ul C,\ol \eta}^\times.
\]
Note that $g(m)$ measures the defect of $\theta(m)$ to be fibrewise
constant. Thus
\[
g^{-1}\cdot\theta:P_\eta\lra  \M_{X,\ul f(\eta)},\quad
p\longmapsto g^{-1}(p)\cdot\theta(p)
\]
factors over $\M(\eta)_\eta$. From the analogous property of
$\theta$ the constructed homomorphism $P_\eta\to \M(\eta)_\eta$ is
right-inverse to
\[
\M(\eta)_\eta\to \M_{X,\ul f(\eta)}\to \ol \M_{X,\ul
f(\eta)}=P_\eta,
\quad m\longmapsto \ol m.
\]
To prove the conclusion it remains to show that if $m_1,m_2\in
\M(\eta)_\eta$ fulfill $\ol m_1=\ol m_2$ there exists a unique $h\in
R^\times$ with $m_2=h\cdot m_1$. By the definition of the ghost sheaf
there certainly exists $h\in \O_{\ul C, \eta}^\times$ such that
$m_2=h\cdot m_1$ holds in $\M_{X,\ul f(\eta)}$. But then also
$f_K^\flat(m_2|_{\eta_K})= h\cdot f_K^\flat(m_1|_{\eta_K})$ holds in
$\M_{C,\ol \eta_K}$. Since both $m_1,m_2$ are fibrewise constant it
follows $h\in K^\times$. Hence $h\in K^\times\cap \O_{\ul C,
\eta}^\times= R^\times$ as claimed.
\qed
\medskip

To motivate our construction of the log structure on $R$ note that
given $\M_R$, the product
\begin{equation}\label{factorization over M_R}
\M_{R,0}\lra Q\times_{Q_K} \Gamma(\M_K),
\end{equation}
of the quotient $\M_{R,0}\to \ol\M_{R,0}= Q$ and the generization
map exhibits $\M_{R,0}$ canonically as a submonoid of
$Q\times_{Q_K}\Gamma(\M_K)$. The point is that we can determine this
image without knowledge of $\M_{R,0}$. To this end let $\alpha_0:
\M_R^0 \to \O_{\Spec R}$ be the basic log structure on $\Spec R$ for
the pre-stable curve $(\ul C/\Spec R, \mathbf x)$. The universal
property of this log structure produces a homomorphism
\[
\beta_0: \M_{R,0}^0\lra \M_{R,(0)}^0\lra \Gamma(\M_K).
\]
Here $(0)\in \Spec R$ denotes the generic point. From the definition
of $Q$ there are also canonical homomorphisms
\[
\psi_\eta:\M(\eta)_\eta \lra P_\eta\lra Q,\quad
\psi_0: \M^0_{R,0} \lra \prod_q\NN\lra Q.
\]
Define $\beta:=(\prod_\eta \beta_\eta)\cdot \beta_0$ and
$\psi:=(\sum_\eta \psi_\eta)+\psi_0$. Clearly,
$\psi_\eta\times\beta_\eta$ and $\psi_0\times\beta_0$ have image in
$Q\times_{Q_K} \Gamma(\M_K) \subset Q\times\Gamma(\M_K)$. Thus
$\psi\times\beta$ maps $(\prod_\eta \M(\eta)_\eta)\times\M_{R,0}^0$
to $Q\times_{Q_K} \Gamma(\M_K)$. We are now in position to construct
$\shM_R$ canonically.

\begin{construction}
\label{Construction of M_R}
Let $\M'_R$ be the log structure on $R$ extending $\M_K$ by the stalk
\begin{equation}\label{defn of M_R}
\M'_{R,0}:= \im(\psi\times\beta)
\subset Q\times_{Q_K} \Gamma(\M_K)
\end{equation}
at $0$ and by the projection to $\Gamma(\M_K)$ as the generization
map to the generic point. Note that $\M'_{R,0}$ may not be
saturated, so we now take its saturation to define $\M_{R,0}$. Since
$Q$ and $\Gamma(\M_K)$ are saturated, so is $Q\times_{Q_K}
\Gamma(\M_K)$. Thus $\M_{R,0}$ is still canonically a submonoid of
$Q\times_{Q_K} \Gamma(\M_K)$. For the structure homomorphism
$\alpha_{R,0}:\M_{R,0}\to R$ to be compatible with generization it
must be taken as composition
\begin{equation}\label{Structure Map}
\M_{R,0}\lra Q\times_{Q_K} \Gamma(\M_K)\stackrel{\pr_2}{\lra}
\Gamma(\M_K) \lra K,
\end{equation}
which we claim has image in $R$. In fact, since $a\in K$ lies in $R$
iff $a^d\in R$ for some $d>0$ it suffices to check this statement
before saturation. But the restriction of \eqref{Structure Map} to
$\M'_{R,0}$  is compatible with $(\prod_\eta \alpha_\eta)\times
\alpha_0$ and hence indeed has image in $R$.
\end{construction}

Here is the key technical result for establishing that $\M_R$
has the requested properties.

\begin{lemma}\label{preimage of psi^gp}
$\beta^\gp\big( \ker(\psi^\gp)\big)= R^\times$, where
$R^\times\subset \Gamma(\M_K)$ is canonically embedded via the
inverse of the structure homomorphism $\Gamma(\M_K)\to K$.
\end{lemma}

\proof
Each $\beta_\eta$ as well as $\beta_0$ are equivariant with respect
to the multiplication action of $R^\times$, while $\psi_\eta$ and
$\psi_0$ are invariant under this action. Since $1\in R^\times$ lies
in the image it is therefore enough to show $\beta^\gp
(\ker(\psi^\gp)) \subset R^\times$. 

An element $m=\big((m_\eta)_\eta,m_0\big) \in \big(\prod_\eta
\M(\eta)_\eta^\gp\big) \times(\M_{R,0}^0)^\gp$ lies in
$\ker\psi^\gp$ iff $\ol m \in\prod_\eta P_\eta^\gp \times\prod_q\ZZ$
lies in the relation subgroup $R_{\mathbf u}$ defining $Q$. In
particular, $\psi(\ol m)$ generizes to $0\in Q_K$, and hence
$\beta^\gp(\ker (\psi^\gp))\subset K^\times$. It also suffices to
check the statement for $\ol m$ equal to one of the generators
$a_q(\ol m_q)$, $\ol m_q\in P_q$, of $R_{\mathbf u}$, as defined in
Construction~\ref{Q and basic log structure}. Note this holds
regardless of saturation since $a\in K^\times$ lies in $R^\times$
iff $a^d\in R^\times$ for some $d>0$. Thus we have
to show that if $m=\big((m_\eta)_\eta,m_0\big) \in \big(\prod_\eta
\M(\eta)_\eta^\gp\big) \times (\M_{R,0}^0)^\gp$ maps to some
$a_q(\ol m_q)$ then $\ord_0 (\beta^\gp(m))=0$, where $\ord_0$ is the
discrete valuation of $K$. Note from the definition of $a_q(\ol
m_q)$ that $\beta_\eta(m_\eta)\in R^\times$ except possibly for
$\eta=\eta_1,\eta_2$, the generic points of the branches of $\ul
C_0$ at $q$. Similarly, all entries of $\ol m_0\in \prod_q\ZZ$
vanish except at entry $q$, which equals $u_q(\ol m_q)$. By working
on an \'etale neighbourhood of $q$ we may assume $\ul C\to \Spec R$
is locally given by $\Spec \big(R[z,w]/(zw-t_q)\big)$ with $t_q\in
R$. Let $m_q \in \M_{X,\ul f(q)}$ be a lift of $\ol m_q$.

There are two cases, depending on $q$ being an isolated node or not.
Let us first assume $q$ is isolated. Then $t_q\neq 0$ and
$\beta_0(m_0)=h\cdot t_q^{u_q(\ol m_q)}$ for some $h\in R^\times$.
Denote by $\eta_K\in\ul C_K$ the generic point with
$\eta_1,\eta_2\in \cl(\eta_K)$ as before. Let
$\M_{C_K,\ol\eta_K}\simeq Q_K\times \O_{\ul C,\ol\eta_K}^\times$ be
an isomorphism induced by a distinguished chart as in
Theorem~\ref{Thm: structure of log curves}. In particular,
$\pi_K^\flat: \Gamma(\M_K)\to \M_{C_K, \ol\eta_K}$ is then given by
the natural inclusion $Q_K\times K^\times\to Q_K\times \O_{\ul
C_K,\ol \eta_K}^\times$. As in \eqref{u_q at isolated nodes}
consider the composition
\[
\kappa_q: \M_{X,\ul f(q)} \lra \M_{X,\ul f(\eta_K)} 
\stackrel{f_K^\flat}{\lra} \M_{C_K,\ol\eta_K}\simeq
Q_K\times \O_{\ul C,\ol\eta_K}^\times \stackrel{\pr_2}{\lra}
\O_{\ul C,\ol\eta_K}^\times.
\]
Since $\ol m_q$ generizes to $\ol m_{\eta_1}$ and to $-\ol
m_{\eta_2}$, respectively, the generizations of $m_q$ at $\eta_1$,
$\eta_2$ differ from $m_{\eta_1}$ and $m_{\eta_2}^{-1}$ only by
functions invertible at $\eta_i$. Therefore
\begin{eqnarray*}
\ord_0(\beta(m))&=& \ord_0\big(\beta_{\eta_1}(m_{\eta_1})
\cdot\beta_{\eta_2}(m_{\eta_2})\cdot \beta_0(m_0)\big)\\
&=& \ord_{\eta_1}\big(\kappa_q(m_q)\big)
-\ord_{\eta_2}\big(\kappa_q(m_q)\big)
+\ord_0(t_q^{u_q(\ol m_q)}).
\end{eqnarray*}
Now by the discussion of the map~\eqref{u_q at isolated nodes} at
the end of \S\ref{par: Ghost sheaves} we know
\[
e\cdot u_q(\ol m_q) =\ord_{\eta_2}\big(\kappa_q(m_q)\big)
-\ord_{\eta_1}\big(\kappa_q(m_q)\big),
\]
and $e=\ord_0(t_q)$. Thus $\ord_0(\beta(m))=0$ as claimed. Note that
the argument seems to depend on the choice of chart inducing the
isomorphism $\M_{C_K,\ol\eta_K}\simeq Q_K\times \O_{\ul
C,\ol\eta_K}^\times$, but different choices cancel out in the
formula for $\ord_0(\beta(m))$ due to the different signs for
$\eta_1$ and $\eta_2$.

For a non-isolated node $q$ let $q_K$ be the node of $\ul C_K$ with
$q\in\cl(q_K)$. Now $t_q=0$ and a distinguished chart for $C_K$ at
$q_K$ takes the form
\[
\M_{C_K,\ol q_K}\simeq
(Q_K\oplus_\NN \NN^2)\times \O_{\ul C_K,\ol q_K}^\times.
\]
Again denote by $\kappa_q: \M_{X,\ul f(q)}\to \O_{\ul C_K,\ol
q_K}^\times$ the composition of $f^\flat(\,.\,|_{q_K})$ with the
projection induced by the chart, and similarly let $\kappa_K:
\Gamma(\M_K)= Q_K\times K^\times\to K^\times$ be the projection.
Since $\ol m=a_q(\ol m_q)$ the generization of $\ol m_q$ to
$\eta_i$ is $\pm \ol m_{\eta_i}$. Hence the generization of $m_q$ to
$\eta_i$ is $m_{\eta_i}^{\pm 1}$ as element of $\M_{X,\ul
f(\eta_i)}^\gp$ up to an element of $\O^\times_{\ul
C,\eta_i}$. Since $\kappa_K$ is induced from $\kappa_q$ and by
the definition of $\beta_{\eta_i}$, the elements
$\pi_K^*(\kappa_K(\beta_{\eta_i} (m_{\eta_i}^{\pm 1})))$ and $\kappa_q
(m_q)|_{\ol\eta_K}$ of $\O^\times_{\ul C,\ol\eta_K}$ agree up to an
element of $\O^\times_{\ul C,\ol\eta_i}$. Therefore it holds
\begin{equation}\label{ord_0(beta(m)) at non-isolated node}
\pm \ord_0 \big(\kappa_K(\beta_{\eta_i}(m_{\eta_i}))\big) =
\ord_{\eta_i} \big(\kappa_q(m_q)\big)=\ord_{\eta_i}(h),
\end{equation}
with $h:= \kappa_q(m_q)\in \O_{\ul C_K,\ol q_K}^\times$ and opposite
signs for $i=1,2$. Now $h$ is invertible at $\ol q_K$ and, since
$\kappa_q$ factors over a chart for $\M_{C_K}$, the zero and
polar locus of $h$ is contained in $\ul C_0$. Hence we can write
$h=g_q g_K$ with $g_q\in \O_{\ul C,\ol q}^\times$, $g_K\in
K^\times$, and thus $\ord_{\eta_1}(h)=
\ord_0(g_K)=\ord_{\eta_2}(h)$. This shows
\[
\kappa_K (\beta_{\eta_1} (m_{\eta_1}) \cdot
\beta_{\eta_2}(m_{\eta_2}))\in R^\times.
\]
As for $m_0$ note that there is a lift to $\M^0_{R,0}$ of the $q$-th
generator of $\ol \M^0_{R,0}= \prod_q\NN$ mapping to $(0,(1,1),1)\in
(Q_K\oplus_\NN \NN^2) \times \O_{\ul C_K,\ol q_K}^\times$ in our
chart. Hence in our chart $\pi_K^\flat\big( \beta_0(m_0)\big)$ takes
the form $\big((0,(a,a)), h\big)$ for $a= u_q(\ol m_q)\in
\NN\setminus \{0\}$ and $h\in \O_{\ul C,\ol q}^\times$. Therefore
also $\kappa_K(\beta_0(m_0)) \in R^\times$. Taken together we obtain
$\beta(m)\in R^\times$ as claimed. The independence of this argument
of the choice of charts follows as in the case of isolated nodes.
\qed

\begin{proposition}
The extension $\alpha_R:\M_R \to \O_{\Spec R}$ of $\M_K$ defined by
$\M_{R,0} \subset Q\times_{Q_K} \Gamma(\M_K)$
(Construction~\ref{Construction of M_R}) is a log structure on
$R$ with $\ol \M_{R,0}=Q$. Moreover, if $\M'_R$ is the log structure
on $\Spec R$ for an extension of $(C_K/(\Spec K,Q_K), \mathbf x_K,
f_K)$ to a basic stable log map over $R$, then the image of $\M'_{R,0}$
under the canonical map to $Q\times_{Q_K} \Gamma(\M_K)$ equals
$\M_{R,0}$.
\end{proposition}

\proof
By Lemma~\ref{Lem: fibrewise constant sections} the maps
$\M(\eta)_\eta\to \ol\M_{X,\ul f(\eta)}=P_\eta$ are surjective, as
is $\M_{R,0}^0\to \prod_q\NN$. Hence $\psi: \big(\prod_\eta
\M(\eta)_\eta \big)\times\M_{R,0}^0\lra Q$ is surjective up to
saturation by the definition of $Q$. Since
$\M_{R,0}=\im(\psi\times\beta)^\sat$ the projection to the first factor
\[
\kappa:\M_{R,0}\lra Q\times_{Q_K} \Gamma(\M_K)\lra Q
\]
is also surjective up to saturation. Since $\kappa$ is also
$R^\times$-invariant we obtain a homomorphism $\M_{R,0}/R^\times\to
Q$ that is surjective up  to saturation. Conversely, if $m\in
\M_{R,0}$ maps to a non-zero element in $Q$ then $\alpha_{R,0}(m)
\in R$ is not invertible. Indeed, let $\big((m_\eta), m_0\big)\in
\prod_\eta \M(\eta)_\eta\times \M_{R,0}^0$ be a lift of (a power of)
$m$. Then by the definition of $\psi$ we must have $\psi_0(m_0)\neq
0$ or $\psi_\eta(m_\eta)\neq 0$ for some $\eta$. But since
$\alpha_\eta$ and $\alpha_0$ are induced by charts of log
structures, one of $\alpha_0(m_0)$ or $\alpha_\eta(m_\eta)$ is not
invertible. Hence also $\alpha_{R,0}(m)$ is not invertible. This
shows $\alpha_{R,0}^{-1}(R^\times)\subset \M_{R,0}\cap
\big(\{0\}\times \Gamma(\M_K)\big)$. Thus Lemma~\ref{preimage of
psi^gp} says that $\alpha_{R,0}$ maps $\alpha_{R,0}^{-1}(R^\times)$
isomorphically to $R^\times \subset K^\times$. In particular,
$\alpha_{R,0}$ induces an injection of $\M_{R,0}/R^\times$ into $Q$
that is surjective up to saturation. Since both $\M_{R,0}$ and $Q$
are saturated, this map is an isomorphism. We have thus established
that $\M_R\to \O_R$ is a log structure on $R$ with $\ol\M_{R,0} =Q$.

For the uniqueness statement it suffices to produce a factorization
of $\psi\times\beta$ over $\M'_{R,0}$:
\[
\big(\prod_\eta \M(\eta)_\eta\big)\times\M_{R,0}^0\lra \M'_{R,0}
\lra Q\times_{Q_K}\Gamma(\M_K).
\]
In fact, then $\M_{R,0}\subset \im \big(\M'_{R,0}\to
Q\times_{Q_K}\Gamma(\M_K)\big)$ and equality in this inclusion
follows by the equality of ghost sheaves via basicness. The
factorization is immediate on $\M_{R,0}^0$ by the universal property
of this log structure. For $\M(\eta)_\eta$ recall first that
$\Gamma(U(\eta), \ul f^*\ol \M_X)= P_\eta$. Thus for each $\ol m\in
P_\eta$ we obtain an $\O_{U(\eta)}^\times$-torsor $\shL_{\ol m}
\subset (\ul f^* \M_X)|_{U(\eta)}$. Then the extension $f^\flat$ of
$f_K^\flat$ induces an isomorphism of $\shL_{\ol m}$ with an
$\O_{U(\eta)}^\times$-torsor in $\M'_C|_{U(\eta)}$. But $(\ul
C,\M'_C)\to (\Spec R,\M'_R)$ is strict on $U(\eta)$, and hence this
torsor is trivial and equal to the pull-back of an
$\O_R^\times$-torsor in $\M'_R$. Moreover, if $m\in
\Gamma(U(\eta),\shL_{\ol m})$ maps to the pull-back of a section of
$\M'_R$ over $\Spec K$ then there exists a unique
$m'\in\Gamma(\M'_R)$ with $f^\flat(m)= \pi^\flat(m')$. The map
$m\mapsto m'$ defines the desired factorization $\M(\eta)_\eta\to
\M'_{R,0}\to Q\times_{Q_K}\Gamma(\M_K)$.
\qed

\subsection{Extension of the log morphism}
\label{par: Extension of morphism}
Our log structure $\M_R$ comes with a morphism $\M_R^0\to \M_R$,
which by the universal property of $\M_R^0$ extends the given
structure of a log smooth curve on $\ul C_K\to \Spec K$ to $\ul C\to
\Spec R$, in the category $(\mathrm{Log}/S)$. Moreover, by uniqueness of
$\M_{R,0}\subset Q\times_{Q_K} \Gamma(\M_K)$ this extension is
uniquely isomorphic to the log structure on the domain $\ul C$ of
any extension of $(C_K/(\Spec K,Q_K), \mathbf x_K, f_K)$ as a basic
stable log map. Denote by $\varphi: \ul f^*\ol\M_X\to \ol\M_C$ the
extension of $\ol{f_K^\flat}$ constructed in \S\ref{par: Ghost
sheaves}.

Now if $U\subset \ul C$ is an open subset and $\ol m\in\Gamma(U, \ul
f^*\ol\M_X)$ we obtain two $\O_{U}^\times$-torsors in $\ul
f^*\M_X|_U$ and in $\M_C|_U$, respectively. If $\shL_{\ol m}$ and
$\shL_{\varphi({\ol m})}$ are the corresponding line bundles the
question is if the isomorphism
\[
f^\flat_K|_{\shL_{\ol m}}: \shL_{\ol m}|_{\ul C_K\cap U}
\lra \shL_{\varphi({\ol m})}|_{\ul C_K\cap U}
\]
extends to $U$, thus defining $ f^\flat|_{\shL_{\ol m}}$
uniquely. Note that since $\ul C$ is Cohen-Macaulay any such
extension is unique, and an extension exists iff it exists in
codimension one, that is, at the generic points $\eta\in \ul C_0$.
But if $m\in \M_{X,\ul f(\eta)}$ there exists $h\in \O_{\ul
C,\eta}^\times$ with $m':=h\cdot m\in\M(\eta)_\eta$. Moreover,
$\M_{C,\ol\eta}= (\ul\pi^* \M_{R,0})_{\ol\eta}$. Using the definition
of $\M_{R,0}$ we therefore see that
\[
f_K^\flat(m)= h^{-1}\cdot f_K^\flat(m')= h^{-1}\cdot
\pi_K^\flat \big((\psi_\eta\times\beta_{\eta})(m')\big)
\]
extends over $\eta$ as $h^{-1}\cdot \pi^\flat (m'')$ with $m'':=
(\psi_\eta\times \beta_\eta) (m')\in\M_{R,0}$. This
proves the unique existence of an extension $f$ of $f_K$ to $\ul C$,
the last step in our proof of Theorem~\ref{Thm: Stable reduction}.
\qed

\section{Log Gromov-Witten invariants}
\label{Sect: VFC}

We are now in position to define log Gromov-Witten invariants by
constructing a virtual fundamental class on $\MM(X)$. This is quite
standard by now. We follow the method of Behrend and Fantechi
\cite{behrend},\cite{behrendfantechi} with the improvement by Kresch
\cite{kresch}. The method has been adapted by Kim to the log setting
\cite{kim} using Olsson's log cotangent complex \cite{OlssonCC}. A
low tech approach based on a global version of Artin's obstruction
theory \cite{artin}, close to the original approach of Li and Tian
\cite{litian} and avoiding cotangent complexes completely is also
possible \cite{si1},\cite{si2}, but this would be less economic.

Some remarks on logarithmic cotangent complexes are in order. An
argument by W.~Bauer presented in \cite{OlssonCC}, \S7, shows that
there is no theory of logarithmic cotangent complexes with the
following four properties: (1)~For strict morphisms one retrieves
the ordinary cotangent complex. (2)~For a log smooth morphism $X\to
Y$ it is represented by $\Omega^1 _{X/Y}$. (3)~Functoriality.
(4)~Compositions of maps yields distinguished triangles. The point
is that (1), (2) and (4) imply compatibility with base change of the
ordinary cotangent complex $L^\bullet_{\ul Y'/\ul Y}$ for a
morphism of schemes $\ul Y'\to \ul Y$ underlying an arbitrary
log \'etale morphism $Y'\to Y$. Since such morphisms need not be
flat this is not true in general.

Olsson presents two ways out of this. The first method defines
$L^\bullet_{X/Y}$ as the ordinary cotangent complex of the asociated
morphism of algebraic stacks $\ul X\to \Log_{\ul Y}$
(\cite{OlssonCC}, Definition~3.2). This version of the logarithmic
cotangent complex fulfills (1), (2) and (3), but (4) only holds
under an additional assumption (Condition~(T)). The second method,
proposed by Gabber, works by enhancing the simplicial resolution
approach to the cotangent complex by log structures
(\cite{OlssonCC}, \S8). The resulting object $L^{G,\bullet}_{X/Y}$
fulfills (1), (3) and (4), but (2) only provided $f$ is integral.
The two approaches agree in many respects (see \cite{OlssonCC},
Theorem~8.32 and Corollary~8.34). In our case we need to represent
the cotangent complex of $X\to S$ by a locally free sheaf but do not
want to impose integrality of $X\to S$, and hence decided for the
first version.

We should also point out that we need an extension of the theory
developed in \cite{OlssonCC} for morphisms of log schemes to
morphisms of log algebraic stacks with Deligne-Mumford domains. This
should be a straightforward generalization, but since this statement
is not available in the literature we will show in Remark~\ref{etale
cover remark} how to reduce to the case of morphisms of log schemes
by working on an \'etale cover. For the moment let us ignore this
issue.

We now require $X\to S$ to be log smooth and quasi-projective. Let
$\U\to \MM$ be the universal curve over the log stack $\MM$ of (not
necessarily basic) pre-stable curves over $S$, as discussed in
Appendix~\ref{App: Prestable curves}. Let $\MM(X)= \MM(X/S)$ be the
log stack of basic stable log maps over $S$,
\[
\pi:\V:= \MM(X)\times_\MM \U\lra \MM(X)
\]
the universal curve and
\[
f:\V\lra X
\]
the evaluation morphism, both considered as $1$-morphisms of
algebraic log stacks. The commutative diagram
\[\begin{CD}
\V@>f>> X\\
@VVV @VVV\\
\U @>>> S
\end{CD}\]
gives rise to the morphism of log cotangent complexes
(\cite{OlssonCC}, (1.1.2))
\[
Lf^* L^\bullet_{X/S}\lra L^\bullet_{\V/\U}.
\]
Now $\V\to\U$ is strict and $\ul\pi:\ul\V\to\ul{\MM(X)}$ is obtained
from the flat morphism $\ul\U\to\ul\MM$ by a base change. Hence, by
\cite{OlssonCC}, 1.1,(ii) and by the compatibility of the ordinary
cotangent complex with flat base change (\cite{illusie}, II.2.2), it holds
\[
L^\bullet_{\V/\U}\simeq L^\bullet_{\ul\V/\ul\U}\simeq
L\pi^* L^\bullet_{\ul{\MM(X)}/\ul\MM}.
\]
In contrast to previous usage we have to indicate now  by
underlining when we want to view $\MM(X)$ or $\MM$ as ordinary
stacks rather than as log stacks. Tensoring with the relative dualizing
sheaf $\omega_{\ul\pi}$ and using the fact that
$L\pi^!L^\bullet_{\ul{\MM(X)}/\ul\MM} \simeq
L\pi^*L^\bullet_{\ul{\MM(X)}/\ul\MM} \stackrel{L}{\otimes}
\omega_{\ul\pi}$, adjunction now defines a morphism
\[
\phi:R\pi_*(L f^* L^\bullet_{X/S} \stackrel{L}{\otimes} \omega_{\ul\pi})
\lra L^\bullet_{\ul{\MM(X)}/\ul\MM}.
\]
By log smoothness $L^\bullet_{X/S}= [\Omega^1_{X/S}]$
(\cite{OlssonCC}, 1.1,(iii)) is represented by a locally free sheaf.
Hence by duality, the left-hand side equals
\[
E^\bullet:= \big(R\pi_*[f^*\Theta_{X/S}]\big)^\vee,
\]
which is of perfect amplitude contained in $[-1,0]$. We have thus
constructed a morphism
\begin{equation}\label{obstructionmorph}
\phi: E^\bullet \lra L^\bullet_{\ul{\MM(X)}/\ul\MM}.
\end{equation}

\begin{proposition}
The morphism \eqref{obstructionmorph}
is a perfect obstruction theory relative to $\ul\MM$ in the
sense of \cite{behrendfantechi}, Definition~4.4.
\end{proposition}

\proof
We are going to check the obstruction-theoretic criterion of
\cite{behrendfantechi}, Theorem~4.5,3, applied relative to $\ul\MM$. To
this end let $T\to\bar T$ be a square zero extension of schemes over
$\ul\MM$ with ideal sheaf $\shJ$ and $g: T\to \ul{\MM(X)}$. We
consider $T$ and $\bar T$ as endowed with the log structures making
the morphisms to $\MM$ strict. Denote by $\V_T= T\times_{\MM(X)} \V
= T\times_\MM \U$, $\V_{\bar T}= \bar T\times_\MM \U$, and by
$p:\V_T\to T$, $\tilde g:\V_T\to\V$ the projections. Let $f_T:=
f\circ \tilde g: \V_T\to X$ be the induced log morphism. We are thus
led to the following commutative diagram.
\[
\parbox{0pt}{\xymatrix{
&\V_T\ar@/^2pc/[rrrr]^{f_T}\ar[rr]_{\tilde g}\ar[ld]^p\ar[dd]&&
\V\ar[rr]_f\ar[ld]^\pi\ar[dd]&&X\ar[dd]\\
T\ar[rr]^>>>>>g\ar[dd]&&\scrM(X)\ar[dd]\\
&\V_{\bar T}\ar[rr]\ar[ld]&& \U\ar[ld]\ar[rr]&&S\\
\bar T\ar[rr]&&\scrM
}}
\]
Note that all but the front and back faces of the cube are
cartesian. Because $\scrM(X)\to\scrM$ is strict and representable
standard obstruction theory applies to the extension problem of $g$
to $\bar T$ (\cite{illusie}, Ch.III, and \cite{OlssonDefTheory} for
the case of representable morphisms of stacks). Thus the
obstructions class $\omega(g)\in \Ext^1_{\O_T}
(Lg^*L^\bullet_{\ul{\scrM(X)}/\ul\scrM}, \shJ)$ to the existence of
an extension is given by the composition
\[
Lg^* L^\bullet_{\ul{\scrM(X)}/\ul\scrM}
\lra L^\bullet_{T/\bar T} \lra \tau_{\ge -1} L^\bullet_{T/\bar
T}=\shJ[1],
\]
defined by functoriality of the cotangent complex. Now
\cite{behrendfantechi}, Theorem~5.3,3, says that $\phi$ defines an
obstruction theory in the sense of Behrend and Fantechi if the
following hold for any extension problem. (1)~An
extension exists if and only if $\phi^*\omega(g)=0$, and (2)~in this
case the set of isomorphism classes of extensions form a torsor
under $\Hom_{\O_T} (Lg^*E^\bullet,\shJ)$. Here $\phi^*\omega(g)\in
\Ext^1_{\O_T} (Lg^*E^\bullet,\shJ)$ is the image of $\omega(g)$ under
pull-back by $\phi$.

On the other hand, by the definition of $\scrM(X)$, such an
extension exists if and only if $f_T:\V_T\to X$ extends as a log
morphism to $\V_{\bar T}$. \cite{OlssonCC}, Theorem~8.45, provides
the obstruction theory for this situation. From this point of view
there is an obstruction class $o\in \Ext^1_{\O_{\V_T}} (Lf_T^*
L^\bullet_{X/S}, p^*\shJ)$, and the isomorphism classes of 
extensions form a torsor under $\Hom_{\O_{\V_T}} (Lf_T^*
L^\bullet_{X/S}, p^*\shJ)$.  The obstruction class $o$ is defined by
the morphism\footnote{The proof of \cite{OlssonCC}, Theorem~8.45,
proceeds with Gabber's version of the cotangent complex, but in view
of \cite{OlssonCC}, Theorem~8.32, this has no influence on the
obstruction class.}
\begin{equation}\label{log obstruction class}
Lf_T^\bullet L^\bullet_{X/S}\lra L^\bullet_{\V_T/S} \lra
L^\bullet_{\V_T/\V_{\bar T}} \to \tau_{\ge -1}L^\bullet_{\V_T/\V_{\bar T}} 
=\shJ[1].
\end{equation}

To compare the two obstruction situations note that repeated
application of adjunction yields the following sequence of
identifications, for any $k\in\NN$:
\begin{eqnarray*}\label{sequence of isos}
\Ext^k_{\O_T}(Lg^*E^\bullet, \shJ)&=&Ext^k_{\O_T}
\big(Lg^*R\pi_*(Lf^*L^\bullet_{X/S}\stackrel{L}{\otimes} \omega_\pi),
\shJ \big)\\
&=& \Ext^k_{\O_\V} \big(Lf^*L^\bullet_{X/S}\stackrel{L}{\otimes}
\omega_\pi, L\pi^!(Rg_*\shJ) \big)\\
&=& \Ext^k_{\O_\V} \big(Lf^*L^\bullet_{X/S},
R\tilde g_*(p^*\shJ) \big)
\ =\ \Ext^k_{\O_{\V_T}} \big(Lf_T^*L^\bullet_{X/S}, p^*\shJ \big).
\end{eqnarray*}
For the third equality we used $L\pi^!= L\pi^*\otimes\omega_\pi$ and
$L\pi^*\circ R g_*= R\tilde g_*\circ Lp^*$, the latter by
flatness of $\pi$. Tracing through this sequence of isomorphisms
for $k=1$ now indeed maps $\phi^*\omega(g)$ to the obstruction
morphism $o$ in \eqref{log obstruction class}. This proves the
obstruction part of the criterion in \cite{behrendfantechi},
Theorem~5.3,3. The torsor part follows readily from \eqref{sequence
of isos} with $k=0$.
\qed
\medskip

Finally, by \cite{behrendfantechi}, \S5 and \cite{kresch},
Theorem~6.2.1, we now have a \emph{virtual fundamental class} $\lfor
\MM(X/S)\rfor$, as a well-defined rational Chow class on
$\ul{\MM(X/S)}$. Moreover, if $\beta$ is a class of stable log maps
to $X$ fulfilling the maximality condition of Definition~\ref{log
classes},ii then $\MM(X/S,\beta)\subset \MM(X/S)$ is an open
substack and hence also carries a virtual fundamental class.

\begin{remark}\label{etale cover remark}
To reduce to an \'etale cover let us recall the construction of the
virtual fundamental class in \cite{behrendfantechi}.  Behrend and
Fantechi first construct a canonical cone stack $\mathfrak C_{Z/S}$
for any morphism of Artin stacks $Z\to S$ that is relatively
Deligne-Mumford, the \emph{intrinsic normal cone}. It is a closed
substack of the \emph{intrinsic normal sheaf} $\mathfrak N_{Z/S}$, a
Picard stack over $Z$ (that is, it has an additive structure
relative $Z$). Both stacks are functorial for \'etale morphisms by
\cite{behrendfantechi}, Proposition~3.14. A perfect obstruction
theory $E^\bullet\to L^\bullet_{Z/S}$ then provides a closed
embedding $\mathfrak C_{Z/S}\to h^1/h^0(E^\vee)$ into a vector
bundle stack, where we use the notation of \cite{behrendfantechi}.
Again, this embedding is compatible with \'etale morphisms, so can
be constructed on an \'etale cover. One then obtains a cone stack of
pure dimension $0$ inside the vector bundle stack $h^1/h^0(E^\vee)$
over $\MM(X)$. Intersecting with the zero section as defined in
\cite{kresch} then defines the virtual fundamental class.

The theory of log cotangent complexes only enters in the
construction of the perfect obstruction theory. Since the
obstruction theory is compatible with \'etale morphisms
we can go over to an \'etale cover and work with diagrams of
log schemes rather than stacks.
\end{remark}

\section{The relationship with expanded degenerations}
\label{Sect: Expanded degenerations}

Let $\ul X$ be a non-singular variety and $D\subset \ul X$ a smooth
divisor, and define $X$ as the log scheme with the divisorial log
structure $\M_{(X,D)}$.  This is the case of relative Gromov-Witten
invariants considered by \cite{LiRuan}, \cite{IonelParker},
\cite{Gathmann}, \cite{JunLi1}. It is insightful to compare the
moduli space of stable log maps in this context with the moduli space
constructed by Jun Li.

Conjecturally, the Gromov-Witten invariants defined using these
moduli spaces will coincide, though we will make no attempt to prove
this here\footnote{Note added in final revision: This statement has
now been proved \cite{AMW}.}. On the other hand, the moduli
spaces themselves are demonstrably not isomorphic.

We will sketch here the relationship between these
moduli spaces, assuming familiarity with Li's notion of 
stable relative maps. 

Let $\M_X$ be the divisorial log structure on $\ul X$ defined by
$D$, $\beta$ a class of stable log maps, $X=(\ul X,\M_X)$, and let
$\MM(X,\beta)$ be as usual the stack of basic stable log maps of
class $\beta$. Let $\mathbf{M}( \ul X/D,\beta)$ be Jun Li's moduli
space of stable relative maps. In fact, Jun Li constructed a log
structure on this stack, but it is not saturated, and as a
consequence, we cannot get a morphism $\mathbf{M}( \ul
X/D,\beta)\to\MM(X,\beta)$. Rather, one must pass to the saturation
$\mathbf{M}( \ul X/D,\beta)^{\sat}$, constructed using \cite{Og},
II, Prop.\ 2.4.5. We will then obtain a morphism of fs log stacks
$\Psi:\mathbf{M}( \ul X/D,\beta)^{\sat}\to\MM(X,\beta)$.

Recall a family of relative stable maps over a base scheme $\ul W$
is given by the following data. For each $n\ge 0$, there is a pair
$(\ul X[n],D[n])$  constructed from $(\ul X,D)$ where $\ul X[n]$ is
defined over $\AA^n$ and $D[n]\subset \ul X[n]$ is a divisor.
Then a family of relative stable maps is a diagram
\begin{equation}
\label{JLicurve}
\xymatrix@C=30pt
{ (\ul  C,{\bf x})\ar[d] \ar[r]^{\ul f} &\ul X[n]\ar[d]^g\\
\ul W\ar[r]&\AA^n}
\end{equation}
which is a family of ordinary stable maps in $\ul X[n]$ satisfying
certain conditions (predeformability and finiteness of a certain
notion of automorphism group), and has specified tangencies with
$D[n]$. In \S1.3 of \cite{JunLi2}, Li observes that $\ul X[n]$ and
$\AA^n$ carry canonical log structures. The log structure on
$\AA^n$ is the divisorial log structure associated to the
divisor $B$ given by $t_1\cdots t_n=0$, and the one on $\ul X[n]$ is
induced by $g^{-1}(B)\cup D[n]$. Furthermore, there is a natural
projection $\ul \Theta:\ul X[n]\to\ul X$ such that
$\ul\Theta^{-1}(D)\subset g^{-1}(B)\cup D[n]$, hence giving rise to a
log morphism $\Theta: X[n]\to X$. Li then constructs log structures
on $\ul C$ and $\ul W$ making \eqref{JLicurve} into a commutative
diagram of log schemes. Now $C\to W$ is a log smooth curve with
marked points along $\ul f^{-1}(D[n])$.

To pass between relative stable maps and stable log maps, we need
the following proposition:

\begin{proposition}
\label{junliprop}
There is a commutative diagram of fine log schemes
\[
\xymatrix@C=30pt
{C\ar[d]_{\pi}\ar[r]^f\ar[dr]^{\psi}
&X[n]\ar[dr]^{\Theta}&\\
W&\bar C\ar[l]^{\ol\pi}\ar[r]_{\bar f}&X}
\]
so that $(\bar C/W,\mathbf x, \bar f)$ satisfies all conditions of
being a stable log map except that the log structures on $\bar C$
and $W$ need not be saturated.
\end{proposition}

\proof
It is standard that such a diagram exists at the level of schemes,
with $\bar f$ a stabilization of $\ul\Theta\circ \ul f$
(\cite{behrendmanin}, Theorem~3.6).  The map $\ul\psi$ contracts in
every fibre of $\ul\pi$ every $\PP^1$ component with only two special
points on which $\ul\Theta\circ \ul f$ is constant.  We need to understand
these morphisms at the log level.
\medskip

\textbf{Step 1}. \emph{Review of $g:X[n]\to X$}.
The pair $(\ul X[n],D[n])$ is constructed inductively: $(\ul
X[0],D[0])$ is the pair $(\ul X,D)$. Then the pair $(\ul X[n],D[n])$
is obtained by blowing up $\ul X[n-1]\times\AA^1$ along
$D[n-1]\times \{0\}$, and $D[n]$ is the proper transform of
$D[n-1]\times\AA^1$. If $D$ is given by the vanishing of a regular
function $w_0=0$, (which can always be accomplished locally on $\ul
X$) then in fact $\ul X[n]$ can be described as the subscheme of
$\ul X\times (\PP^1)^n \times \AA^n$ given by the equations
\[
w_0 z_1=w_1t_1,\quad w_1z_2=w_2z_1t_2,\quad\ldots,\quad w_{n-1}z_n
=z_{n-1}w_nt_n,
\]
where $z_i,w_i$ are homogeneous coordinates on the $i$-th copy of
$\PP^1$ and $t_1,\ldots,t_n$ are coordinates on $\AA^n$.
This is covered by affine open subsets $\ul X_i$, $1\le i\le n+1$,
where $\ul X_i$ is given by $z_1=\cdots=z_{i-1}=1=w_i=\cdots=w_n$.
Thus $w_{i-1}z_i=t_i$ on $\ul X_i$ for $i\le n$. The log structure
of $X[n]$ restricted to $\ul X_i$ for $i\le n$ has a chart
\[
\NN^n\oplus_{\NN}\NN^2\to \O_{X_i}
\]
where the map $\NN\to \NN^n$ is $1\mapsto e_i$, the 
$i$-th generator of $\NN^n$, and the chart is
\[
\left(\sum_{i=1}^n a_ie_i, (a,b)\right)\longmapsto 
\left(\prod_{i=1}^n t_i^{a_i}\right)
w_{i-1}^az_i^b.
\]
The chart on $X_{n+1}$ is just $\NN^n\oplus \NN\to \O_{X_{n+1}}$ 
given by $(\sum a_i e_i,a)\mapsto w_n^a\prod t_i^{a_i}$.

Recall that $X[n]$ fibres over $\AA^n$ with coordinates $t_1,\ldots,t_n$,
and that the critical locus of this map consists of $n$ distinct
subvarieties $D_1,\ldots,D_n$, with $D_i$ sitting over $t_i=0$.
In the above description, $D_i\subset X_i$ is given by
$w_{i-1}=z_i=0$. Also $D[n]\subset X_{n+1}$ is given by $w_n=0$.

The map $\Theta:X[n]\to X$ is the projection to $X$.
We describe this as a log morphism as follows.
Let $s_{w_0}$ be the section of $\M_X$ corresponding
to the function $w_0$ vanishing only along $D$. In general, if
we have a chart $P\to\O_Y$ for a log structure, for $m\in P$
we denote by $s_m$ the corresponding section of the associated log
structure, so that any section of the associated log structure
is of the form $h\cdot s_m$ for some $m\in P$, $h\in \O_Y^{\times}$.
Then $\Theta:
X_i\to X$ is given by (for $i\le n$)
\[
s_{w_0}\longmapsto s_{(e_1+\cdots+e_{i-1},(1,0))}
\]
since on $\ul X_i$, $w_0=w_{i-1}t_{1}\cdots t_{i-1}$,
while $\Theta:X_{n+1}\to X$
is given by 
\[
s_{w_0}\longmapsto s_{(e_1+\cdots+e_n,1)}.
\]

\medskip

\textbf{Step 2.} \emph{Review of $C\to W$}.
Fix a geometric point $\ol w\in |W|$. Let ${\bf D}$ be the set of 
\emph{distinguished} double  points of $\ul C_{\ol w}$, that is,
double points mapping to $\bigcup_{\ell=1}^n D_\ell$, with ${\bf
D}_{\ell} \subset {\bf D}$ the subset mapping to $D_{\ell}$. Let
${\bf U}$ be the set of undistinguished double points. For each
$q\in {\bf D}_{\ell}$, let $\mu_q$ denote the order of tangency of
either branch of $C_{\ol w}$ at $q$ with $D_{\ell}$. Let $N_{\ell}$
be the free monoid $\NN\rho_{\ell}$  generated by $\rho_{\ell}$ if
${\bf D}_{\ell}=\emptyset$; otherwise $N_\ell$ is the monoid
generated by  $\{\rho_q\,|\,q\in {\bf D}_{\ell}\}$ modulo the
relations $\mu_q\rho_q=\mu_{q'}\rho_{q'}$ for each $q,q'\in {\bf
D}_{\ell}$. Then\footnote{Jun Li does not include the contributions
from the undistinguished nodes; these must be included in order for
$C\to W$ to be log smooth. However, these contributions will play no
further role in the discussion.}
\[
\ol\M_{W,\ol w}=\bigoplus_{\ell=1}^n N_{\ell}\oplus
\bigoplus_{q\in {\bf U}} \NN\rho_q.
\]
Note the monoids $N_{\ell}$ need not be saturated; as a consequence,
the log structures on $W, C$ and $\bar C$ we construct below need
not be saturated.

For every point $q\in {\bf D}_{\ell}$, one can find an \'etale open
neighbourhood $U_q$ of $q\in  \ul C$ such that $\ul f$ maps $U_q$ into
$\ul X_{\ell}$. Furthermore, $U_q$ is of the form $\Spec
A[x_q,y_q]/(x_qy_q-t_q)$ where \'etale locally $W=\Spec A$, and
$t_q\in A$. As observed in \cite{JunLi2}, Simplification 1.7, we can
choose $U_q$ and coordinates $x_q,y_q$ so that
\[
\ul f^*(w_{\ell-1})=x_q^{\mu_q},
\quad \ul f^*(z_{\ell})=y_q^{\mu_q}.
\]
Then $\ul f^*(t_\ell)=t_q^{\mu_q}$.

Similarly, for each $q\in {\bf U}$, we can describe a neighbourhood
$U_q$ of $q$ in $\ul C$ as 
\[
U_q\cong\Spec A[x_q,y_q]/(x_qy_q-t_q)
\]
for some $t_q\in A$.

Given these choices, Li puts log structures on $C$ and $W$ as follows.
There is a chart for the log structure of $W$, 
\[
\sigma:\ol\M_{W,\ol w} \lra A
\]
given by
\[
\rho_q\longmapsto t_q, \quad \rho_{\ell}\longmapsto \ul h^*t_{\ell},
\]
where $\ul h:\ul W\to \AA^n$. A chart for the log structure on $C$
for the neighbourhood $U_q$ is
\begin{align*}
\psi_q:\ol\M_{W,\ol w}\oplus_{\NN}\NN^2&\lra \O_{U_q}\\
(\alpha,(a,b))&\longmapsto \ul\pi^*(\sigma(\alpha))x_q^ay_q^b.
\end{align*}
Here the map $\NN\to\ol\M_{W,\ol w}$ is $1\mapsto \rho_q$. 

If $p$ is any marked point of $\ul C$ with $\ul f(p)\in D[n]$, then
in a suitable neighbourhood $U_p=\Spec A[x_p]$ of $p$, with chart
$\ol\M_{W,\ol w}\oplus\NN\to A[x_p]$ given as usual as
$(\alpha,a)\mapsto \pi^*(\sigma(\alpha))x_p^a$. We have $\ul
f^*(w_n)= x^{\mu_p}_p \cdot h_p$ for a unit $h_p$ and some
$\mu_p\in\NN$, which will be assumed to be $1$ via a suitable choice
of $U_p$ and $x_p$. 

One checks these induced log structures do not depend on any choices
and glue uniquely, yielding a well-defined $C \to W$.

\textbf{Step 3}. \emph{The log morphism $f:C \to
X[n]$.} For $q\in {\bf D}_i$, $f: (U_q, \M_C|_{U_q}) \to X_i$ is given by
\[
s_{(\sum a_{\ell}e_{\ell},(a,b))}\longmapsto s_{(\sum a_{\ell}\rho_{\ell}',
(\mu_q a,\mu_q b))}
\]
where $\rho_{\ell}'=\rho_{\ell}$ if ${\bf D}_{\ell}$ is empty and is
equal to $\mu_q\rho_q$ for any $q\in {\bf D}_{\ell}$ otherwise. If
$p$ is a marked point of $\ul C$ with $\ul f(p)\in D[n]$, then  $f$ is
given on $(U_p,\M_C|_{U_p})\to X_{n+1}$ by
\[
s_{(\sum a_{\ell}e_{\ell},a)}\longmapsto
s_{(\sum a_{\ell}\rho_{\ell}',\mu_pa)}.
\]
This is sufficient to completely specify $f$.

\textbf{Step 4}. \emph{The map $C\to \bar C$.}
The log structure on $\ul{\bar C}$ can be described as follows.
First, if $\ul{\bar C}^o\subset \ul{\bar C}$ is the largest open set for
which $\ul\psi^{-1}(\ul{\bar C}^o)\to \ul{\bar C}^o$ is an isomorphism, 
$\bar C^o$ has the same log structure as $\psi^{-1}(\bar C^o)$.

Next, let $\ol y\in |\ul{\bar C}|$ be a geometric point not in $\ul{\bar
C}^o$, with $\ol w=\ul{\ol\pi}(\ol y)$. There are two cases: $\ol y$
either  is not, or is, a double point of $\ul{\bar C}_{\ol w}$. In
both cases, $\psi^{-1}(\ol y)$ is a chain $\ul C_1\cup\cdots\cup\ul
C_m$ of rational curves, with $q_i=\ul C_i\cap\ul C_{i+1}$ a double
point, and $\ul C_{\ol w}$ has another component $\ul C_0$ with $\ul
C_0\cap\ul  C_1=q_0$ also a double point. However, if $\ol y$ is not
a double point, then there is a point $p_m\in \ul C_m$ which is a
log marked point, while if $\ol y$ is a double point, $\ul C_{\ol w}$
has a component $\ul C_{m+1}$ and a double point $q_{m}=\ul C_m\cap
\ul C_{m+1}$.

If $\ol y$ is not a double point, then there is an \'etale open
neighbourhood of $\ol y$ of the  form $U:=\Spec A[x]$, where $x=0$
is the image of the section $\ul\psi\circ x_i$, with $x_i$ the marked
point of $\ul C$ corresponding to $p_m$.

If $\ol y$ is a double point, then one can show that there is an 
\'etale open neighbourhood of $\ol y$
of the form $U:=\Spec A[x,y]/(xy-\prod_{i=0}^m t_{q_i})$. 

In both cases, one can describe $U\times_{\ul{\bar C}} \ul C$.  In
the case that $\ol y$ is not a double point, then this scheme is
given by the equations
\[
xu_1=v_1t_{q_0},\quad v_1u_2=u_1v_2t_{q_1},\quad \ldots,
\quad v_{m-1}u_{m}=u_{m-1}v_mt_{q_{m-1}}
\]
in $\AA^1\times (\PP^1)^{m}\times \Spec A$, with homogeneous
coordinates on the $i$-th $\PP^1$ being $u_i,v_i$. If $\ol y$ is a
double point, then this scheme is given by the equations
\[
xu_1=v_1t_{q_0},\quad v_1u_2=u_1v_2t_{q_1},\quad \ldots,
\quad v_{m-1}u_{m}=u_{m-1}v_{n}t_{q_{m-1}}, \quad yv_{m}=u_{m}t_{q_m}
\]
in $\AA^2\times (\PP^1)^{m}\times\Spec A$.

We can cover $U\times_{\ul{\bar C}}\ul C$ with Zariski open
subsets $U_i$ for $0\le i\le m$, where $U_i$ is the
set where $u_1=\cdots=u_i=1=v_{i+1}=\cdots=v_m$. Note $U_i$ is
an open neighbourhood of $q_i$, $U_i\cong \Spec
A[v_i,u_{i+1}]/(v_iu_{i+1}-t_{q_i})$ if $i\le m$,
and $U_m\cong\Spec A[v_m]$ if $\ol y$ is not a double point.

We have to be slightly careful with the coordinates $x$ and $y$;
these can't be chosen arbitrarily. To do so we have to relate these
to the map $\ul f$. At this point, we shall make a simplifying assumption
that $\ol w\in W$ maps to $0\in \AA^n$. This can always be achieved
locally on $W$ by decreasing $n$, and the general case can be dealt
with by the reader with some extra book-keeping.

Possibly after reversing the order of $\ul C_0,\ldots,\ul C_{m+1}$
in the case that $\ol y$ is a double point, 
we can assume that there is some positive integer $\ell_0$
such that $\ul f(q_i)\in D_{\ell_0+i}$ for 
$i\le m$, and $\ul f(p_m)\in \ul X_{n+1}\cap D[n]$ if
$\ol y$ is not a double point. In the latter case, $\ell_0=n-m+1$.
The map $\ul f:U_i\to \ul X[n]$ factors through $\ul X_{\ell_0+i}$, and
necessarily yields 
\begin{align*}
i=0:\quad\quad & \ul f^*(w_{\ell_0-1})=h_{0,1} x^{\mu}, \quad 
\ul f^*(z_{\ell_0})=h_{0,2}u_1^{\mu}\\
0<i<m: \quad\quad &
\ul f^*(w_{\ell_0+i-1})=h_{i,1} v_{i}^{\mu},\quad \ul f^*(z_{\ell_0+i})=
h_{i,2} u_{i+1}^{\mu}\\
i=m: \quad\quad &
\ul f^*(w_{\ell_0+m-1})=h_{m,1} v_m^{\mu},\quad \ul f^*(z_{\ell_0+m})=
h_{m,2} y^{\mu},
\end{align*}
the last line if $\ol y$ is a double point. Here $\mu=\mu_{q_i}$
for any $i$.
If $\ol y$ is not a double point, then $U_m\to \ul X[n]$
factors through $X_{n+1}$, with
\[
f^*(w_n)=h_{m,1}v_m^{\mu}.
\]
Now in fact we can assume $h_{0,2}=h_{m,1}=1$, $h_{i,k}=1$ for $0< i
<m$, $k=1,2$. This is because $\ul f$ on each component $\ul C_i$, $1\le
i\le m$,  is just a $\mu$-fold cover of $\PP^1$ totally ramified at
$0$ and $\infty$, so that the above listed $h_{i,k}$'s are constant,
and then after applying a suitable change of coordinates on $\ul C$, one
can assume these constants are $1$. Further, by making a change of
coordinates for $x$ (and $y$), one can assume $h_{0,1}=1$ (and
$h_{m,2}=1$ in the double point case). 

In particular, one can then assume that the open sets
$U_{q_0},\ldots,U_{q_{m-1}},U_{p_m}$ if $\ol y$ is not a double
point ($U_{q_0},\ldots,U_{q_m}$ if $\ol y $ is a double point) are
taken to be the open sets $U_0,\ldots, U_m$. Thus we know what the
chart for the log structure on $U_i$ is, by Step 2. This is given on
$U_i$ for $i\le m-1$ (and for $i=m$ in the case of a double point)
by charts  $\ol\M_{W,\ol w}\oplus_{\NN}\NN^2\to \O_{U_i}$ given by 
\[
(\alpha,(a,b))\longmapsto \ul\pi^*(\sigma(\alpha))v_i^au_{i+1}^b,
\]
and on $U_m$ for $\ol y$ not a double point by a chart
$\ol\M_{W,\ol w}\oplus\NN\to\O_{U_m}$ given by
\[
(\alpha,a)\longmapsto \ul\pi^*(\sigma(\alpha))v_m^a.
\]

Using this particular choice of the coordinate $x$ (and $y$ if
$\ol y$ is a double point), we can put a log structure on $U$
as follows. If $\ol y$ is not a double point, we have a chart
$\ol\M_{W,\ol w}\oplus\NN\to \O_{ U}$ given by
\[
(\alpha,a)\longmapsto \ul{\ol\pi}^*(\sigma(\alpha))\cdot x^a.
\]
If $\ol y$ is a double point, we have a chart
\begin{align*}
\ol\M_{W,\ol w}\oplus_{\NN}
\NN^2&\lra \O_{ U}\\
(\alpha,(a,b))&\longmapsto
\ul{\ol\pi}^*(\sigma(\alpha)) x^ay^b,
\end{align*}
with the map $\NN\to  \ol\M_{W,\ol w}$ given by $1\mapsto
\sum_{i=1}^m\rho_{q_i}$. It is straightforward to check that this
log structure is independent of choices and this description gives
compatible log structures for different choices of the point $\ol
y$. Hence one obtains a log structure on $\ul{\bar C}$ clearly making
$\ol\pi$ log smooth. 

We can specify the log morphism $\psi:U_i\to U$ as
follows. If $\ol y$ is not a double point, this map is defined as
\[
s_{(\alpha,a)}\longmapsto 
\begin{cases}s_{(\alpha+a\sum_{j<i}\rho_{q_j},(a,0))}&i<m;\\
s_{(\alpha+a\sum_{j<m}\rho_{q_j},a)}&i=m. 
\end{cases}
\]
On the other hand, if $\ol y$ is a double point, then for every $i$,
we take the map
\[
s_{(\alpha,(a,b))}
\longmapsto s_{(\alpha+a\sum_{j<i}\rho_{q_j}+b\sum_{j>i}\rho_{q_j},(a,b))}.
\]
Using the explicit description for the charts on the sets $U_i$, one
checks that the morphisms $U_i\to U$ agree on overlaps and
hence give a morphism $U\times_{\bar C} C \to U$.
Furthermore, these morphisms are compatible for different choices of
open neighbourhoods $U$ of different points $\ol y$, giving the log
morphism $\psi:C\to\bar C$.

\textbf{Step 5}. \emph{The map $\bar f:\bar C
\to X$}. We now define the log map
$\bar f:U\to X$
by
\[
s_{w_0}\longmapsto 
\begin{cases}s_{(\sum_{\ell=1}^{\ell_0-1}\rho_\ell',\mu)}\in \M_U&
\hbox{if $\ol y$ is not a double point;}\\
s_{(\sum_{\ell=1}^{\ell_0-1}\rho_\ell',(\mu,0))}\in\M_U&\hbox
{if $\ol y$ is a double point,}
\end{cases}
\]
where $\rho'_{\ell}$ is defined in Step 3.
One checks from the description of all the maps above that
this yields the desired commutative diagram, at least where
these maps are now defined. Further, one checks that all maps
are independent of choices and glue, hence giving the desired
global commutative diagram.
\qed

\medskip

By \cite{Og}, II 2.4.5, for any fine log scheme $W$ there
is a finite surjective morphism $W^{\sat}\rightarrow W$
from an fs log scheme $W^{\sat}$, such that every morphism
$W'\rightarrow W$ with $W'$ fine and saturated has a unique
factorization through $W^{\sat}\rightarrow W$. Thus
we can define an fs log stack $\mathbf{M}(\ul X/D,\beta)^{\sat}$
by defining an object over a scheme $\ul{W}$ to be a choice
of relative stable map $f:C/W\rightarrow X[n]$ (with log structures
as defined by Jun Li), yielding
a log scheme $W^{\sat}$, and a choice of section $\ul{W}\rightarrow
\ul{W}^{\sat}$ of $\ul{W}^{\sat}\rightarrow \ul{W}$. Denote by
$W'$ the pull-back fs log structure on $\ul{W}$ from $W^{\sat}$
under this morphism. Then from this data, using Proposition \ref{junliprop},
we obtain a stable log map $\bar C\times_W W'\rightarrow X$ over $W'$.
Thus we obtain:

\begin{corollary}
\label{Cor: Comparison with Jun Li's}
There is a morphism of stacks
$\Psi:\mathbf{M}(\ul X/D,\beta)^{\sat}\to \MM(X,\beta)$.
\end{corollary}

We shall see in the next section that this morphism is not in general
expected to be an isomorphism.

\section{Examples}

\begin{example}
\label{basicrelativeexample}
Let $\ul{X}$ be a non-singular variety with a smooth divisor $D\subset
\ul{X}$ as in the previous section, yielding the log scheme $X$
with the divisorial log structure.
Consider first the case where $\ul C$ is a smooth curve and $\ul f: \ul
C\to \ul X$ is an ordinary stable map such that $\ul f^{-1}(D)$
consists of a finite number of points. We would like to understand
when this can be lifted to a stable log map. We know that if $\ul f$
lifts to a stable log map, then $\ul f^*\ol\M_X$ can only jump at
marked and double points, by Remark \ref{Rem: structure of M}, 
and hence every point of $\ul f^{-1}(D)$ must
be marked. So we should consider a situation
\[
\ul f:\big(\ul C,(x_1,\ldots,x_d,y_1,\ldots,y_p)\big )\lra \ul X
\]
with $\ul f^{-1}(D)=\{x_1,\ldots,x_d\}$.

What is the possible type of a log lifting? Since $\ul C$ has no
double points, the only relevant information is the choice of
$u_{x_i}\in\NN^{\vee}=\NN$, where $\NN=\ol\M_{X,f(x_i)}$.  (The
$u_{y_i}$'s are necessarily zero as they lie in the zero monoid.) 
Furthermore, $P_{\eta}=0$ for $\eta$ the generic point of $\ul C$,
so once the type is chosen, necessarily $Q=0$. Thus $\M_C$ is just
the divisorial log structure on $\ul C$ associated with the divisor
$x_1+\cdots+y_p$.

The only constraint on the choice of $u_{x_i}$ comes from the
balancing condition. More precisely, note that the torsor coming
from $n\in\Gamma(\ul X, \ol\M_X)=\NN$ corresponds to the line bundle
$\O_{\ul X}(-nD)$. Thus in the notation of \eqref{tau_eta},
$\tau_\eta^X:\Gamma(\ul C,\ul f^*\ol{\M}_X)=\NN^d\to \ZZ$ is given
by  $\tau_\eta^X(n_1,\ldots,n_d)=-\sum \mu_in_i$, where $\mu_i>0$ is
the order of tangency of $D$ with $\ul C$ at the point $x_i$. On the
other hand, the map $\tau_\eta^C:\Gamma(\ul C,\ol\M_C)=
\NN^{d+p}\to\ZZ$ is given by $\tau_\eta^C(n_1,\ldots,n_d,m_1,
\ldots,m_p)=-\sum n_i-\sum m_j$. Since the map $\varphi: \Gamma(\ul
C,\ul f^*\ol\M_X)\to\Gamma(\ul C,\ol\M_C)$ induced by $\ol{f^\flat}$
is given by
\[
\varphi(n_1,\ldots,n_p)=(u_{x_1} n_1,
\ldots,u_{x_d} n_d,0,\ldots,0),
\]
the only way that $\tau_\eta^C\circ \varphi$ can coincide with
$\tau_\eta^X$ is if $u_{x_i}=\mu_i$ for $1\le i\le d$. Thus we
see that the elements $u_{x_i}\in\NN^{\vee}$ can be interpreted
as imposing the orders of tangency.

Note that once $u_{x_i}=\mu_i$ for each $i$, there is a unique
map $\ul f^*\M_X\to \M_C$ induced by $\ul f^*:f^{-1}\O_{\ul X}\to
\O_{\ul C}$, as in this case the structure maps for the log structures
embed $\ul f^*\M_X$ and $\M_C$ in $\O_{\ul C}$. Thus the open substack
of the moduli space of stable log maps corresponding to
curves considered in this example coincides with the corresponding
substack (not necessarily open now) of the moduli of ordinary
stable maps consisting of stable maps as above with the correct
orders of tangencies, and the log structure is trivial.
\end{example}

\begin{example}
\label{squaremonoideg}
Continuing with the case of a pair $(\ul X,D)$ as in the previous
example, it is not difficult to obtain interesting examples for the monoid
$Q$. Consider the case that $\ul X=\PP^1$, $D=\{0\}$, and consider the
limiting situation in Figure~\ref{quadricconeegfigure}. The figure
on the left shows a rational curve with marked points
$p_1,p_2,p_3$ mapping $2:1$ to $\PP^1$, with order of tangencies
to $D$ at the three marked points being $0,0$ and $2$ respectively.
Such a situation can degenerate to the stable marked curve $\ul f$ pictured 
to the right, in which all but the two horizontal components map
to $0$. By semi-stable reduction for stable log maps, there is
a limiting stable log map, with some associated monoid $Q$. Let us
determine what the type of this curve is and what $Q$ is.

\begin{figure}
\input{Figure1.pstex_t}
\caption{}
\label{quadricconeegfigure}
\end{figure}

Let $D_i$, $i=1,2,3$ be the irreducible components of the degenerate
domain curve containing the points $p_i$, $i=1,2,3$, respectively.
Let $D_4$ and $D_5$ be the two horizontal components, with $D_4\cap
D_1\not=\emptyset$ and $D_5\cap D_2\not=\emptyset$. Denote by
$\eta_i$ the generic point of $D_i$. Because $\ul f(D_i)=\{0\}$,
$i=1,2,3$, we have $\tau^X_{\eta_i}=0$ for these $i$, but
$\tau^X_{\eta_4}=\tau^X_{\eta_5}: \NN\to \ZZ$ is given by $n\mapsto
-n$, as in Example \ref{basicrelativeexample}. We have the four
double points 
\[
q_1= D_4\cap D_1,\quad q_2=D_1\cap D_3,\quad q_3=D_2\cap D_3, \quad
q_4=D_5\cap D_2.
\]

To determine the type of the central fibre, one can use the method
of the proof of Proposition \ref{Prop: Uniqueness of type}. In
particular, necessarily $u_{p_1}=u_{p_2}=0\in\NN$ and
$u_{p_3}=2\in\NN$. However, to determine the $u_{q_i}$, it is easier
to use the balancing condition Proposition~\ref{Prop: Balancing
condition} to observe that there is only one choice of type for the
central fibre with these $u_{p_i}$. Indeed, suppose we give the
central fibre the structure of a stable log map over the standard
log point as in \S\ref{par: Log maps over std log pt}. This means in
particular we have data $V_{\eta_1},\ldots,V_{\eta_5}$ with
$V_{\eta_1}, V_{\eta_2},V_{\eta_3}\in \NN^{\vee}$ and
$V_{\eta_4},V_{\eta_5} =0$, as well as positive integers
$e_{q_1},\ldots,e_{q_4}$.  Suppose the $u_{p_i}$'s are as given
above. If we define the $u_{q_i}$'s by the choice of signs
\begin{align*}
u_{q_1}=(V_{\eta_1}-V_{\eta_4})/e_{q_1},\quad
&u_{q_2}=(V_{\eta_3}-V_{\eta_1})/ e_{q_2}\\
u_{q_3}=(V_{\eta_3}-V_{\eta_2})/e_{q_3},\quad
&u_{q_4}=(V_{\eta_2}-V_{\eta_5})/ e_{q_4}
\end{align*}
then by Proposition~\ref{Prop: Balancing condition} 
necessarily $u_{q_i}=1\in\NN^{\gp}$ for all $i$. 
Then the monoid $Q$ associated to this type is 
$(P_{\eta_1}\oplus P_{\eta_2}\oplus P_{\eta_3}\oplus \NN^4)/R_{\bf u}$,
where $R_{\bf u}$ is generated by
\begin{align*}
a_{q_1}(1)= {} &(-1,0,0,1,0,0,0),\\
a_{q_2}(1)= {} &(1,0,-1,0,1,0,0),\\
a_{q_3}(1)= {} &(0,1,-1,0,0,1,0),\\
a_{q_4}(1)= {} &(0,-1,0,0,0,0,1).
\end{align*}
By eliminating the last four components using these relations,
$\ZZ^7/R_{\bf u}$ can be identified with $\ZZ^3$, with the monoid
$\NN^7$ having image in $\ZZ^3$ being generated by $e_1=(1,0,0)$,
$e_2=(-1,0,1)$, $e_3=(0,-1,1)$ and $e_4=(0,1,0)$.
Thus $Q$ is the monoid generated by $e_1,e_2,e_3,e_4$ subject to the
relation $e_1+e_2=e_3+e_4$. This is the monoid defining the quadric
cone in $\AA^4$. 

Note the choice of a map $Q\to\NN$
yields a tropical curve in $(\Gamma(\ul C,\ul f^*\ol\M_X)^{\gp})^*$,
as depicted in Figure~\ref{tropcurve1}. \qed

\begin{figure}
\input{tropcurve1.pstex_t}
\caption{The tropical curve associated to $f$. Here the range is
$\Hom(\Gamma(\ul C,\ul f^*\ol\M_X)^{\gp},\RR)$.}
\label{tropcurve1}
\end{figure}

\begin{example}
Let us return to the situation of Example~\ref{squaremonoideg},
comparing our picture with the expanded degeneration picture as
described in Section~\ref{Sect: Expanded degenerations}. In our
picture, the degenerate curve $f$ corresponds to one point in the
moduli space. On the other hand, there are three combinatorially
distinct possible limits in the expanded degeneration picture, as
depicted in Figure~\ref{expandedpicture}. In Cases I and II, one
must insert additional irreducible components into the limit  curve
$C$: these components are contracted in our picture, as they are not
stable components. Roughly speaking, these two cases correspond, in
our language, to curves over $(\Spec \kappa,\NN)$ with
$V_{\eta_1}<V_{\eta_2}$ and $V_{\eta_1}>V_{\eta_2}$ respectively. In
Case III, however, the limit domain is the same as ours, but note
there is still some actual moduli. Indeed, the first bubble
component of the range has four special points: the two double
points of the range contained in this component, and the images of
the marked points $P_1$ and $P_2$. The cross-ratio of these four
points provides a one-parameter moduli space, a copy of $\Gm$. Cases
I and II can be viewed as limit points in this one-parameter moduli
space. So the map $\mathbf{M}(\ul X/D,\beta)\to \MM(X,\beta)$
contracts a $\PP^1$ to a point.

This suggests that in fact $\mathbf{M}(\ul X/D,\beta)$ is, in this
case, only birational onto its image in $\MM(X,\beta)$. We
conjecture this is the case in general. In this case, this morphism
gives a small resolution of an ordinary double point in the
three-dimensional moduli space of basic stable log maps of the type being
considered here. Nevertheless, one may conjecture that the log
Gromov-Witten invariants defined here coincide with relative
Gromov-Witten invariants as defined by Jun Li. We shall leave it to
others to prove such a conjecture\footnote{Note added in final
revision: This conjecture has been verified in \cite{AMW}.}.
\end{example}

\begin{figure}
\input{expandedpicture.pstex_t}
\caption{}
\label{expandedpicture}
\end{figure}
\end{example}

\begin{example} 
\label{twocomponentexample} 
Next let us consider $\ul X=\PP^2$ with $D\subset \ul X$ a union of
two distinct lines $L_1,L_2$. This is a normal crossings divisor, and is
the first case that the classical form of relative Gromov-Witten
invariants does not cover. Consider the moduli space of degree
two stable log rational curves in $\PP^2$ with four marked points
$p_{ij}$, $1\le i,j\le 2$, with $p_{ij}$ a point with tangency of
order $1$ with $L_i$ and tangency of order $0$ with $L_{3-i}$. 
The generic case just consists of a conic intersecting $D$
at four distinct smooth points of $D$. For such curves,
the existence of the double point in $D$ is irrelevant,
and the analysis of Example~\ref{basicrelativeexample} still applies.

Now consider a limiting curve in which the image of the conic
degenerates to a reducible conic $F$ with $F\cap D
=L_1\cap L_2$, that is, the double points of $F$ and $D$ coincide.
By stable reduction, there is a limiting stable map $\ul f$, which
in the simplest case is given in Figure~\ref{conicdegen}.
Furthermore, by stable log reduction, this gives a stable log map.
Let us determine its type.

\begin{figure}
\input{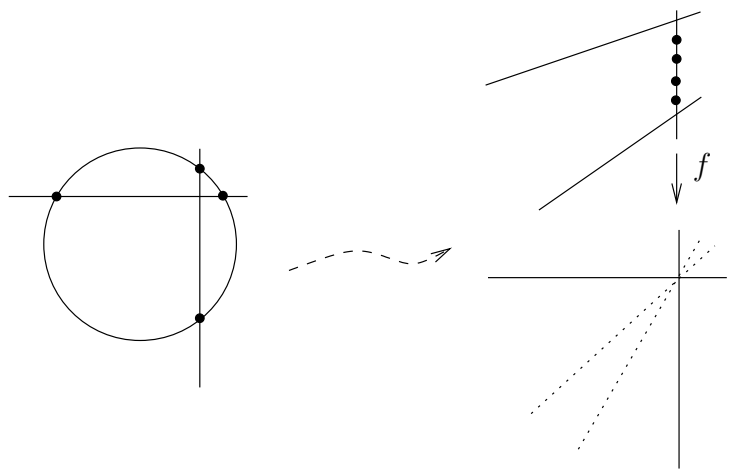}
\caption{The conic on the left degenerates to the reducible curve
on the right; the dotted lines are the image of the stable map.}
\label{conicdegen}
\end{figure}

First note that $\Gamma(\ul C, \ul f^*\ol\M_X)=\NN^2$. Let the irreducible
components of $\ul C$ be $D_1,D_2,D_3$ with $D_3$ the contracted
component. Then necessarily $u_{p_{1j}}=(1,0)$ and
$u_{p_{2j}}=(0,1)$ for $j=1,2$. One can think of this as saying that
even though the $p_{ij}$'s map into the intersection of the
irreducible components, $u_{p_{ij}}$ remembers which component the
curve should be tangent to at $p_{ij}$.

As in Example \ref{squaremonoideg}, we can determine the $u_q$'s by
considering the structure of a stable log map over the standard log
point. Necessarily $V_{\eta_1}=V_{\eta_2}=0$. Letting $q_i=D_i\cap
D_3$, $i=1,2$,   there is no choice but for $u_{q_i}=(1,1)$ (with
the proper choice of order) in order to achieve the balancing
condition at $V_{\eta_3}$. One sees in this case that $Q=\NN$, with
an element of $Q^{\vee}$ specifying $e_{q_1}=e_{q_2}$.
\end{example}

\begin{example}
Let $\ul X$ be a complete toric variety and let $\partial \ul
X\subset \ul X$ be the toric boundary, the union of toric divisors
on $\ul X$. Let $\ul X$ be given the divisorial log structure
induced by this divisor. In general, we do not expect that $\ol\M_X$
will be generated by global sections; however, it is easy to  see
that if $M$ is the character lattice of the torus acting on $X$,
then there is a surjection $\ul M\to \ol\M_X^{\gp}$, as each stalk
of $\ol\M_X^{\gp}$ is a quotient of $M$. In particular, given a
stable log map $(C/(\Spec\kappa,\NN), {\bf x}, f)$ with target $X$,
by Discussion~\ref{Disc: Tropical curve} one obtains a balanced
tropical curve $h:\Gamma_{\ul C}\to N_{\RR}=\Hom(M,\RR)$.  In
general, this tropical curve contains one unbounded edge for each
irreducible component $D =\cl(\eta)$ of $\ul C$ in the direction of
$\tau^X_\eta$, or rather, the image of $\tau^X_\eta$ in
$N=\Hom(M,\ZZ)$. More specifically, given the component $D$ and its
normalization $g:\tilde D\to D$, one has a composition
\[
M\lra \Gamma(\tilde D,g^*\ul f^*\ol\M_X^\gp)\mapright{\tau^X_\eta}\ZZ,
\]
which yields the image of $\tau^X_\eta$ in $N$. Note, however, that 
in the present toric case the composed map $M\to \Gamma(\ul
X,\ol\M_X^{\gp})\to \Pic X$ is the zero map, as the divisor of
zeroes and poles of a monomial $z^n$ is linearly equivalent to zero.
Hence the image of $\tau^X_\eta$ in $N$ is in fact zero.

Thus the construction of Discussion~\ref{Disc: Tropical curve}
yields an ordinary tropical curve in $N_{\RR}$ whose only unbounded
edges necessarily correspond to marked points of the curve $\ul C$, and
the balancing condition of Proposition~\ref{Prop: Balancing
condition} gives the usual tropical balancing condition in
$N_{\RR}$. Furthermore, there is a relationship between this
tropical curve and the fan $\Sigma$ in $N_{\RR}$ for $X$. Indeed,
suppose that for an irreducible component $D=\cl(\eta)$ of $\ul C$, the
toric stratum of $X$ containing  $\ul f(\eta)$ corresponds to a cone
$\tau\in\Sigma$. Then $h(v_{\eta})\in\tau$. Indeed, we have the
composition
\[
M\lra \Gamma(\ul C,\ul f^*\ol\M_X^{\gp})
\twoheadrightarrow P_{\eta}^{\gp},
\]
with $P_{\eta}=(\tau^{\vee}\cap M)/(\tau^{\vee}\cap M)^{\times}$. 
Thus we have the dual map
\[
P_{\eta}^\vee\hookrightarrow (P_{\eta}^{\gp})^* \lra N,
\]
which identifies $P_{\eta}^{\vee}$ with $\tau\cap N$. As $h(v_{\eta})$
is the image of $V_{\eta}\in P_{\eta}^\vee$ under this map, one sees
$h(v_{\eta})\in\tau$.

\begin{figure}
\input{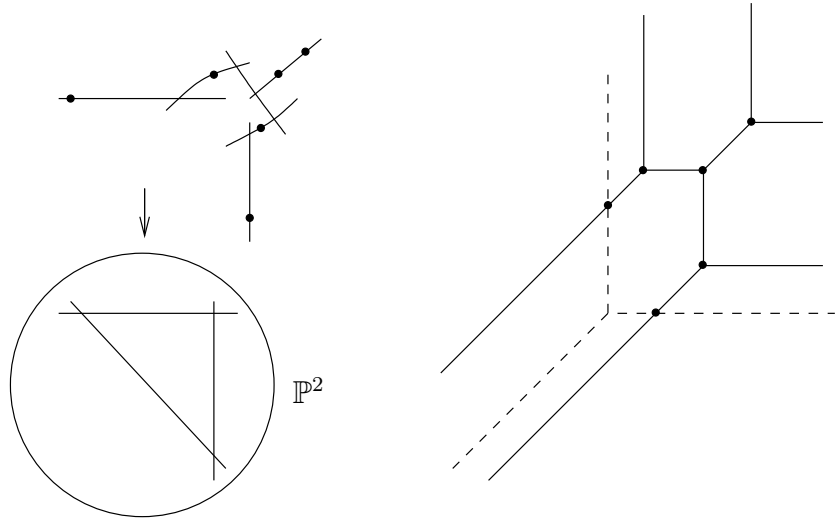}
\caption{The figure on the left depicts the underlying
stable map, with marked points as indicated. All components
are contracted but the left-hand and lower components, which
map to coordinate lines. The right-hand figure shows a possible
corresponding tropical curve, with the dotted lines indicating
the fan for $\PP^2$.}
\label{tropicalconic}
\end{figure}

Note further that the tropical curve $h:\Gamma_{\ul C}\to N_{\RR}$ in
fact determines the type of the curve, as the maps $P_x^{\vee}\to N$
are injective for any $x\in |\ul C|$. However, given an ordinary
stable map $\ul f:\ul C\to \ul X$, it is not true that any type of tropical
curve $h:\Gamma_{\ul C}\to N_{\RR}$ with $h(v_\eta)$ lying in the
correct cone of $\Sigma$ is an allowable type of log map. Indeed,
the balancing condition holding in $N_{\RR}$ does not imply the
balancing condition of Proposition~\ref{Prop: Balancing condition}.

Figure~\ref{tropicalconic} gives an example of a stable map
$\ul f:\ul C\to \PP^2$ and a corresponding possible tropical
curve. It is not difficult to see that $Q=\NN^5$ here.
\end{example}

\begin{appendix}
\section{The log stack of prestable curves}
\label{App: Prestable curves}

In this appendix we will sketch the argument that the stack $\scrM=
\scrM_S$ of pre-stable log smooth curves defined over $S=(\ul
S,\M_S)$ is an algebraic log stack locally of finite type over $S$.
We note that this will be an algebraic stack only in the weaker
sense of \cite{OlssonENS}, in the sense that the diagonal morphism
will not be separated, due to the phenomenon of Example~\ref{Rem:
non-sparatedness issue} inherent in log moduli problems.

We begin by considering the stack $\mathbf{M}$ of ordinary
pre-stable curves over $\ul S$, that is, for a scheme
$V$,
\[
\mathbf{M}(V)=\{\hbox{$(C/V, \,{\bf x})$ is a pre-stable curve}\}.
\]
Denote by $\mathbf{M}_{g,k}$ the Deligne-Mumford stack of stable
curves of genus $g$ with $k$ marked points \cite{delignemumford},
\cite{knudsen}.

\begin{lemma} $\mathbf{M}$ is an algebraic stack, locally of finite
type over $\ul S$.
\end{lemma}

\proof The fact that the diagonal $\Delta:\mathbf{M} \to
\mathbf{M}\times_{\ul S}\mathbf{M}$ is representable follows
from the fact that given two pre-stable curves $(C_1/V,{\bf x}_1)$,
$(C_2/V,{\bf x}_2)$, the isomorphism functor 
\[
\Isom_V\big((C_1/V,{\bf x}_1),(C_2/V,{\bf x}_2)\big)
\]
is representable by a closed subscheme of
$\Isom_V(C_1,C_2)$, which in turn is representable by an open
subscheme of the Hilbert scheme of $C_1\times_V C_2$ by \cite{Gr}.
In particular, $\Delta$ is separated and quasi-compact.

To see that $\mathbf{M}$ has a smooth cover, denote by
$\mathbf{M}_{g,k}^o$ the open sub\emph{scheme} of $\mathbf{M}_{g,k}$
whose points parameterize stable marked curves with no non-trivial
automorphisms. Consider the map
\[
p:\coprod_{g,k,\ell\ge 0} \mathbf{M}_{g,k+\ell}^o\to\mathbf{M}
\]
which takes a curve $(C/V,x_1,\ldots,x_{k+\ell})$ to the curve
$(C/V,x_1,\ldots,x_k)$. Let $C^o$ be the complement of the critical
points of $C\to V$ and let $(C^o)^{\ell}$ denote
$C^o\times_V\cdots\times_V C^o$ ($\ell$ copies).  Given any map
$V\to\mathbf{M}$  corresponding to a pre-stable curve
$(C/V,x_1,\ldots,x_k)$, one sees easily that $V\times_{\mathbf{M}}
\coprod_{g,k,\ell} \mathbf{M}_{g,k+\ell}^o$ is represented by an
open subscheme of $\coprod_{\ell\ge 0} (C^o)^\ell$, and hence $p$ is
smooth. On the other hand, $p$ is clearly surjective, since given
any pre-stable curve over a separably closed field, one can always
add enough marked points so that its automorphism group becomes
trivial.

Since $\mathbf{M}_{g,k}$ is of finite type for each $g$ and $k$,
this shows $\mathbf{M}$ is an algebraic stack locally of finite
type.
\qed

\medskip

Now let us incorporate log structures. Given a pre-stable curve
$(\ul C/ \ul W,{\bf x})$, the argument of \cite{Kato
2000}, p.~227ff, constructs canonical log structures and a log
morphism $(\ul C,\M_C)\to (\ul W,\M_W)$,
with the log structure on $\ul W$ called \emph{basic}. While
the result is stated there for stable curves, stability is not used.
The point of the basicness property is the following:

\begin{proposition}\label{Prop: Basic log curves}
Given $(C/W, {\bf x})$ basic in the above sense, and a pre-stable
marked log curve $(D/Z,{\bf y})$ with maps $\ul \alpha:
\ul Z\to \ul W$ and $\ul \beta:
\ul D\to \ul C$  inducing an isomorphism
$\ul D\to \ul C\times_{\ul W}
\ul Z$, there exist unique maps $\alpha:Z\to W$ and
$\beta:D\to C$ with underlying scheme morphisms $\ul
\alpha$, $\ul \beta$ such that the diagram
\[
\xymatrix@C=30pt
{D\ar[r]^{\beta}\ar[d]&C\ar[d]\\
Z\ar[r]_{\alpha}&W
}
\]
is cartesian.
\end{proposition}

\proof This is the content of Proposition 2.1 and Theorem 2.1 of
\cite{Kato 2000}. The result as stated does not rely on stability of the
curves.
\qed
\medskip

In particular, basicness is stable under base change, that is, given
$\ul W'\to \ul W$, the basic log structure on
$(\ul C\times_{\ul W} \ul W'/\ul W', {\bf
x})$ is the pull-back of the basic log structure on $(\ul
C/\ul W,{\bf x})$. This endows $\mathbf{M}$ with a log
structure,
\[
\mathbf{M}\lra (\mathrm{Log}/S)\lra (\mathrm{Sch}/\ul S),
\]
thus generalizing \eqref{Mgk is a log stack}.

Finally we want to allow arbitrary log structures on the base. For
this we use Olsson's deep result that the stack $\Log_S$ of fine log
schemes over $S$ is algebraic (but not with a separated diagonal)
and locally of finite type over $\ul S$ (\cite{OlssonENS},
Theorem~1.1). Similarly, there is an algebraic stack
$\Log^{\bullet\to\bullet}_S$ of schemes $\ul T$ over $\ul S$
together with a morphism of fine log structures $\M_1\to\M_2$ on $T$
(\cite{OlssonCC}, Example~2.1; this is a direct consequence of
\cite{OlssonENS}, Proposition~5.9). The forgetful functor to the
first log structure $\M_1$ defines a morphism
\[
\Log^{\bullet\to\bullet}_S \lra \Log_S
\]
of algebraic stacks. Clearly, $\mathbf M\to (\mathrm{Log}/S)$
in fact defines a morphism of algebraic stacks
$\mathbf M\to\Log_S$. Now the fibre product
\[
\Log^{\bullet\to\bullet}_S \times_{\Log_S} \mathbf M
\]
is isomorphic to the stack $\MM$ of pre-stable log curves (with
arbitrary log structures on the base). In fact, an object over $\ul
W\in (\mathrm{Sch}/\ul S)$ consists of a \emph{basic} pre-stable log
curve $(C/(\ul W,\M^0_W), \mathbf x)$ and a morphism of fine log
structures $\M_W^0\to \M_W$. Pulling back thus defines a log curve
over $(\ul W,\M_W)$, and Proposition~\ref{Prop: Basic log curves} says
that this functor defines an isomorphism of stacks from the fibre
product to $\MM$. Summarizing, we have the following result.

\begin{proposition}\label{Prop: MM is algebraic}
The stack $\MM$ of pre-stable marked log curves over fine log schemes
over $S$ is an algebraic log stack locally of finite type over $\ul
S$.
\qed 
\end{proposition}

\section{Tropicalization of a log space}
\label{App: Tropicalization}

Traditionally, tropical geometry provides a discrete version of
algebraic geometry in $\GG_m^n$ over some valuation ring. An almost
equivalent point of view is to work in an equivariant
compactification, thus replacing $\GG_m^n$ by a toric variety. The
authors have emphasized at various places that tropical geometry can
also be viewed as providing the discrete information captured in the
ghost sheaf of a log structure. In the present context we have seen
traditional tropical curves to arise from a stable log map
$(C/(\Spec\kappa,\NN), \mathbf{x}, f)$ over the standard log point
provided $\ul f^* \ol\M^\gp$ is globally generated
(Discussion~\ref{Disc: Tropical curve}). Without global
generatedness, a natural target space for a tropical curve might be a
space we call the \emph{tropicalization} of $X$, (analogous to the
notion of the tropical part of an exploded manifold of Parker's work
\cite{parker}), defined as follows.

Given a log scheme $X$ with log structure in the Zariski topology,
we set
\[
\Trop(X):=\Bigg(\coprod_{ x\in \ul X} \Hom(\ol\M_{X, x},
\RR_{\ge 0})\Bigg)\bigg/\sim,
\]
where the disjoint union is over all scheme-theoretic points of
$\ul X$ and the equivalence relation is generated by the
identifications of faces given by dualizing generization maps
$\ol\M_{X, x} \to\ol\M_{X, x'}$ when $x$ is specialization of $x'$.
One then obtains for each $ x$ a map
\[
i_x:\Hom(\ol\M_{X, x},\RR_{\ge 0})\to \Trop(X).
\]
For general $X$, $\Trop(X)$ may not be particularly well-behaved, as
the equivalence relation might yield strange self-identifications of
faces.

\begin{example}
Let $X$ be a log scheme whose underlying scheme is a union $C_1\cup C_2$
of two copies of $\PP^1$ with generic points $\eta_1,\eta_2$, 
meeting at two points $q_1,q_2$ (hence a degenerate 
elliptic curve). It is not difficult to find a Zariski log structure such
that $\overline\M_{X,q_i}\cong \NN^4$, $\overline\M_{X,\eta_i}\cong \NN^3$,
with generization maps: 
\begin{align*}
\overline\M_{X,q_1}\lra\overline\M_{X,\eta_1},& \quad
\sum_i a_ie_i\longmapsto a_1e_1+a_2e_2+a_3e_3,\\
\overline\M_{X,q_1}\lra\overline\M_{X,\eta_2},& \quad
\sum_i a_ie_i\longmapsto a_1e_1+a_2e_2+a_4e_4,\\
\overline\M_{X,q_2}\lra\overline\M_{X,\eta_1},& \quad
\sum_i a_ie_i\longmapsto a_1e_1+a_2e_2+a_3e_3,\\
\overline\M_{X,q_2}\lra\overline\M_{X,\eta_2},& \quad
\sum_i a_ie_i\longmapsto a_2e_1+a_1e_2+a_4e_4.
\end{align*}
Then one sees that $\Trop(X)$ is not particularly well-behaved. It is a quotient
of a disjoint union of two copies of the orthant generated by $e_1^*,\ldots,
e_4^*$, with the face generated by $e_1^*,e_2^*$ of one of these orthants
identified with the same face of the other orthant using two different
identifications: the identity and the identification swapping $e_1^*$ and
$e_2^*$. Hence $\Trop(X)$ cannot be viewed as a polyhedral complex in this
case.

In general, it is helpful to view the tropicalization more abstractly
as a collection of cones with maps between them. This might provide
the correct notion of a tropical stack. This point of view will be
explored in more detail elsewhere.
\end{example}

It is sometimes useful to impose a condition which allows us
to avoid such a possibility:

\begin{definition}
\label{Def: monodromy free}
Let $X$ be an fs log scheme. We say $X$ is \emph{monodromy free} if
for any geometric point $ x\in \ul X$, $i_{ x}$ is injective
on the interior of any face  of $\Hom(\ol\M_{X, x},\RR_{\ge 0})$.
\end{definition}

In general, $\Trop(X)$ is a good range space for the tropical curves
arising from log maps. Specifically, given a stable log map
$(C/(\Spec\kappa,\NN), \mathbf x,f)$ over a standard log point, we
map $v_{\eta}$ to $i_{f(\ol\eta)}(V_{\eta})$, where now $V_{\eta}$
is viewed as an element of $\Hom(P_{\eta},\RR_{\ge 0})$. We map an
edge $E_q$ to the image of the line segment in $\Hom(P_q,\RR_{\ge
0})$ joining the images of the endpoints of $E_q$, as usual. We map
$E_p$ to a ray with endpoint  $i_{f(\ol\eta)}(V_{\eta})$ in the
direction $i_{f(\ol\eta)}(u_p)$, if $p$ is in the closure of
$\eta$. 

\begin{remark}
A more conceptual way to view this construction is to observe that
$\Trop$ is a covariant functor: given a morphism of Zariski
log schemes $f:X\rightarrow Y$, for $x\in X$ we obtain a map
$f^{\flat}:\overline\M_{Y,f(x)} \rightarrow\overline\M_{X,x}$ and
hence a map $\Hom(\overline\M_{X,x},\RR_{\ge 0}) \rightarrow
\Hom(\overline\M_{Y,f(x)},\RR_{\ge 0})$. This is compatible with the
equivalence relations defining $\Trop(X)$ and $\Trop(Y)$. Indeed, if
$x$ is a specialization of $x'$, $f(x)$ is a specialization of
$f(x')$, and hence inclusions of faces are compatible with the
induced maps on cones.

The above description of the tropical curve associated to
$(C/(\Spec\kappa,\NN), {\bf x},f)$ can be viewed as follows. We have
a map $\Trop(\pi):\Trop(C) \rightarrow\Trop(\Spec\kappa,\NN)
=\RR_{\ge 0}$, and $\Trop(f): \Trop(\pi)^{-1}(1)\rightarrow\Trop(X)$
is easily seen to coincide with the tropical curve described above.

More generally, suppose we have a basic stable log map $(C/(\Spec\kappa,Q),
{\bf x},f)$ over a point. Then we obtain a family of tropical curves
parameterized by $\Hom(Q,\RR_{\ge 0})$. In particular, the tropical curves
corresponding to pull-backs of this stable log map to standard log points
are given by the restriction of $\Trop(f)$ to $\Trop(\pi)^{-1}(q)$,
for $q\in \Int(Q^\vee)$.
\end{remark}

\end{appendix}



\begin{thebibliography}{cccccccc}
\bibitem[AbCh]{AbramovichChen} D.~Abramovich, Q.~Chen:
  \emph{Stable logarithmic maps to Deligne-Faltings pairs~II},
  preprint \texttt{arXiv:1102.4531 [math.AG]}, 19pp.  
\bibitem[ACGM]{abramovich etal} D.~Abramovich, Q.~Chen,
  W.D.~Gillam, S.~Marcus:
  \emph{The evaluation space of logarithmic stable maps},
  preprint \texttt{arXiv:1012.5416 [math.AG]}, 19pp.
\bibitem[AMW]{AMW} D.\ Abramovich, S.\ Marcus, J.\ Wise:
  \emph{Comparison theorems for Gromov-Witten invariants of smooth
  pairs and of degenerations},
  preprint \texttt{arXiv:1207.2085 [math.AG]}, 43pp.
\bibitem[Ar]{artin} M.\ Artin:
  \emph{Versal deformations and algebraic stacks},
  Inv.\ Math~\textbf{27} (1974) 165--189.
\bibitem[Be]{behrend} K.\ Behrend:
  \emph{GW-invariants in algebraic geometry},
  Inv.\ Math~\textbf{127} (1997) 601--617.
\bibitem[BeFa]{behrendfantechi} K.\ Behrend, B.\ Fantechi:
  \emph{The intrinsic normal cone},
  Inv.\ Math~\textbf{128} (1997) 45--88.
\bibitem[BeMa]{behrendmanin} K.\ Behrend, Y.\ Manin:
  \emph{Stacks of stable maps and Gromov-Witten invariants},
  Duke Math.\ Journ~\textbf{85} (1996) 1--60.
\bibitem[Ch]{chen} Q.\ Chen:
	\emph{Stable logarithmic maps to Deligne-Faltings pairs I},
	preprint \texttt{arXiv:1008.3090 [math.AG]}, 48pp.
\bibitem[DeMu]{delignemumford} P.\ Deligne, D.\ Mumford:
  \emph{The irreducibility of the space of curves of given genus},
  Publ.\ Math.\ IHES \textbf{36} (1996) 75--110.
\bibitem[FuPa]{fultonpandh} W.~Fulton, R.~Pandharipande:
  \emph{Notes on stable maps and quantum cohomology},
  in: Proceedings of the Algebraic Geometry Conference in Santa Cruz
  1995, Proc.\ Symp.\ Pure Math~\textbf{62}, Part~2, Amer.\ Math.\
  Soc.\ 1997, 45--96.
\bibitem[Ga]{Gathmann} A.~Gathmann:
  \emph{Absolute and relative Gromov-Witten invariants of
  very ample hypersurfaces},
  Duke Math.\ J.~\textbf{115}  (2002) 171--203.
\bibitem[Gr]{Gr} A.~Grothendieck:
  \emph{Techniques de construction et th\'eor\`emes
  d'existence en g\'eom\'etrie alg\'ebrique.~IV,
  Les sch\'emas de Hilbert}, S\'eminaire Bourbaki, Vol.~\textbf{6},
  Exp.\ No.~221 (1960/61), 249--276, Soc.\ Math.\ France, Paris, 1995.
\bibitem[Gr]{gross_seattle} M.\ Gross:
  \emph{The Strominger-Yau-Zaslow conjecture: from torus fibrations
  to degenerations},
  in: Algebraic geometry---Seattle 2005, 149--192,
  Proc.\ Sympos.\ Pure Math., 80, Part 1, Amer.\ Math.\ Soc.~2009.
\bibitem[Il]{illusie} L.~Illusie:
	\emph{Complexe cotangent et d\'eformations I}.
	Lecture Notes Math.~239, Springer~1971.
\bibitem[IoPa]{IonelParker} E.-N.\ Ionel, T.\ Parker:
  \emph{Relative Gromov-Witten invariants},
  Ann.\ of Math.~\textbf{157} (2003), 45--96. 
\bibitem[Kf]{Kato 2000} F.~Kato:
  \emph{Log smooth deformation and moduli of log smooth curves},
  International J.\ Math~\textbf{11} (2000) 215--232.
\bibitem[Kk]{Kato 1989} K.~Kato:
  \emph{Logarithmic structures of
  Fontaine--Illusie}, in: J.-I.~Igusa (ed.) et.~al.,
  Algebraic analysis, geometry, and number theory, pp.~191--224.
  Johns Hopkins Univ.~Pr., Baltimore, 1989
\bibitem[Ki]{kim} B.~Kim: \emph{Logarithmic stable maps},
  preprint \texttt{arXiv:0807.3611v2 [math.AG]}, 30pp.
\bibitem[Kd]{knudsen} F.\ Knudsen:
  \emph{The projectivity of the moduli space of stable curves.\ II.~The
  stacks $M_{g,n}$},
  Math.\ Scand.~\textbf{52} (1983), 161--199.
\bibitem[Kt]{knutson} D.\ Knutson:
  \emph{Algebraic spaces},
  Lect.\ Notes Math.~\textbf{203}, Springer~1971.
\bibitem[Kr]{kresch} A.\ Kresch:
  \emph{Cycle groups for Artin stacks},
  Inv.\ Math.~\textbf{138} (1999), 495--536.
\bibitem[LaMB]{champs} G.\ Laumon, L.\ Moret-Bailly: \emph{Champs
  alg\'ebriques}, Springer 2000.
\bibitem[LiRu]{LiRuan} A.-M.\ Li and Y.\ Ruan:
  \emph{Symplectic surgery and Gromov-Witten invariants of Calabi-Yau
  3-folds, I},
  Invent.\ Math.~\textbf{145} (2001), 151--218. 
\bibitem[Li1]{JunLi1} J.~Li, \emph{Stable morphisms to singular schemes
  and relative stable morphisms},
  J.\ Differential Geom~\textbf{57} (2000) 509--578.
\bibitem[Li2]{JunLi2} J.~Li, \emph{A degeneration formula of
  GW-invariants},
  J.\ Differential Geom.~\textbf{60} (2002) 199--293. 
\bibitem[LiTi]{litian} J.\ Li, G.\ Tian: \emph{Virtual moduli cycles and
  GW-invariants of algebraic varieties},
  Journal Amer.\ Math.\ Soc~\textbf{11} (1998) 119--174.
\bibitem[NiSi]{nisi} T.\ Nishinou, B.\ Siebert:
  \emph{Toric degenerations of toric varieties and tropical curves.},
  Duke Math.\ J.~\textbf{135} (2006), 1--51.
\bibitem[Nt]{Nt} N.\ Nitsure:
	\emph{Construction of Hilbert and Quot schemes},
	Fundamental algebraic geometry, 105--137, Math.\ Surveys
	Monogr.~\textbf{123}, Amer.\ Math.\ Soc.~2005.
\bibitem[Og]{Og} A.\ Ogus: \emph{Lectures on logarithmic algebraic
  geometry.} TeXed notes (2006).
\bibitem[Ol1]{OlssonThesis} M.\ Olsson:
  \emph{Log algebraic stacks and moduli of log schemes},
  Thesis, Berkeley 2002.
\bibitem[Ol2]{OlssonENS} M.\ Olsson:
  \emph{Logarithmic geometry and algebraic stacks}
   Ann.\ Sci.\ École Norm.\ Sup~\textbf{36} (2003), 747--791.
\bibitem[Ol3]{OlssonCC} M.\ Olsson:
  \emph{The logarithmic cotangent complex},
  Math.\ Ann.~\textbf{333}  (2005), 859--931.
\bibitem[Ol4]{OlssonDefTheory} M.\ Olsson:
	\emph{Deformation theory of representable morphisms of algebraic
	stacks},
	Math.\ Z.~\textbf{253} (2006), 25--62.
\bibitem[Pa]{parker} B.\ Parker:
  \emph{Exploded manifolds},
  preprint \texttt{arXiv:0910.4201 [math.SG]}, 63pp.
\bibitem[SGA1]{SGA1} A.\ Grothendieck:
  \emph{Rev\^etement \'etale et groupe fondamental},
  Lecture Notes Math~\textbf{224}, Springer 1971
\bibitem[Si1]{si1} B.\ Siebert:
  \emph{Virtual fundamental classes,
  global normal cones and Fulton's canonical classes},
  in: Frobenius manifolds (C.~Hertling, M.~Marcolli eds.), 341--358,
  Vieweg~2004.
\bibitem[Si2]{Talk} B.\ Siebert:
  \emph{Gromov-Witten invariants in relative and singular cases},
  Lecture given at the workshop ``Algebraic aspects of mirror
  symmetry'', Univ.\ Kaiserslautern, Germany, June~2001.
\bibitem[Si3]{si2} B.\ Siebert:
  \emph{Obstruction theories revisited},
  manuscript~2002.
\end{thebibliography}
\end{document}